\documentclass{article}

\usepackage[titletoc,title]{appendix}
\usepackage[cp1250]{inputenc}
\usepackage[IL2]{fontenc}
\usepackage{yfonts,fancyhdr}
\usepackage{a4wide}
\usepackage[english]{babel}
\usepackage{euscript}
\usepackage{amstext,amsbsy,amscd,amssymb}
\usepackage{amsmath,enumerate}
\usepackage{amsfonts}
\usepackage{mathrsfs,tensor,xy}
\usepackage{graphics}
\usepackage{graphicx}
\usepackage{microtype}
\include{diagxy}

\let\rarr=\rightarrow

\let\larr=\leftarrow

\let\veps=\varepsilon
\let\mcal=\mathcal
\let\mfrak=\mathfrak
\let\eus=\EuScript

\def\N{\mathbb{N}}
\def\Z{\mathbb{Z}}
\def\R{\mathbb{R}}
\def\C{\mathbb{C}}

\def\Con{\mathop {\rm Con} \nolimits}
\def\Mod{\mathop {\rm Mod} \nolimits}

\def\Hom{\mathop {\rm Hom} \nolimits}

\def\ad{\mathop {\rm ad} \nolimits}
\def\gr{\mathop {\rm gr} \nolimits}
\def\im{\mathop {\rm im} \nolimits}

\def\Diff{\mathop {\rm Diff} \nolimits}

\def\SL{\mathop {\rm SL} \nolimits}
\def\diag{\mathop {\rm diag} \nolimits}

\def\Ad{\mathop {\rm Ad} \nolimits}

\def\Sol{\mathop {\rm Sol} \nolimits}

\def\tr{\mathop {\rm tr} \nolimits}
\def\htt{\mathop {\rm ht} \nolimits}

\def\Ind{\mathop {\rm Ind} \nolimits}

\newsymbol\squares 1003

\long\def\proof #1{\noindent \emph{Proof.}\ #1 \hfill $\squares$
\medskip}

\newcounter{num}[section]
\numberwithin{equation}{section}
\numberwithin{num}{section}

\long\def\definition #1 {\refstepcounter{num} \noindent {\bf
Definition \thenum.} #1

\medskip}

\long\def\theorem #1{\refstepcounter{num} \noindent {\bf Theorem
\thenum.} #1

\medskip}

\long\def\lemma #1{\refstepcounter{num}  \noindent {\bf Lemma
\thenum.} #1

\medskip}

\long\def\proposition #1{\refstepcounter{num}  \noindent {\bf
Proposition \thenum.} #1

\medskip}

\newenvironment{enum}{\begin{list}{}{\topsep=2pt \itemsep=0pt
\parsep=0pt}}{\end{list}}

\newcommand*\riso{%
  \xrightarrow[]{\raisebox{-0.25em}{\smash{\ensuremath{\sim}}}}%
}

\makeatletter
\newcommand*\if@single[3]{%
  \setbox0\hbox{${\mathaccent"0362{#1}}^H$}%
  \setbox2\hbox{${\mathaccent"0362{\kern0pt#1}}^H$}%
  \ifdim\ht0=\ht2 #3\else #2\fi
  }
\newcommand*\rel@kern[1]{\kern#1\dimexpr\macc@kerna}
\newcommand*\widebar[1]{\@ifnextchar^{{\wide@bar{#1}{0}}}{\wide@bar{#1}{1}}}
\newcommand*\wide@bar[2]{\if@single{#1}{\wide@bar@{#1}{#2}{1}}{\wide@bar@{#1}{#2}{2}}}
\newcommand*\wide@bar@[3]{%
  \begingroup
  \def\mathaccent##1##2{%
    \if#32 \let\macc@nucleus\first@char \fi
    \setbox\z@\hbox{$\macc@style{\macc@nucleus}_{}$}%
    \setbox\tw@\hbox{$\macc@style{\macc@nucleus}{}_{}$}%
    \dimen@\wd\tw@
    \advance\dimen@-\wd\z@
    \divide\dimen@ 3
    \@tempdima\wd\tw@
    \advance\@tempdima-\scriptspace
    \divide\@tempdima 10
    \advance\dimen@-\@tempdima
    \ifdim\dimen@>\z@ \dimen@0pt\fi
    \rel@kern{0.6}\kern-\dimen@
    \if#31
      \overline{\rel@kern{-0.6}\kern\dimen@\macc@nucleus\rel@kern{0.4}\kern\dimen@}%
      \advance\dimen@0.4\dimexpr\macc@kerna
      \let\final@kern#2%
      \ifdim\dimen@<\z@ \let\final@kern1\fi
      \if\final@kern1 \kern-\dimen@\fi
    \else
      \overline{\rel@kern{-0.6}\kern\dimen@#1}%
    \fi
  }%
  \macc@depth\@ne
  \let\math@bgroup\@empty \let\math@egroup\macc@set@skewchar
  \mathsurround\z@ \frozen@everymath{\mathgroup\macc@group\relax}%
  \macc@set@skewchar\relax
  \let\mathaccentV\macc@nested@a
  \if#31
    \macc@nested@a\relax111{#1}%
  \else
    \def\gobble@till@marker##1\endmarker{}%
    \futurelet\first@char\gobble@till@marker#1\endmarker
    \ifcat\noexpand\first@char A\else
      \def\first@char{}%
    \fi
    \macc@nested@a\relax111{\first@char}%
  \fi
  \endgroup
}
\makeatother

\newcommand\rsetminus{\mathbin{\mathpalette\rsetminusaux\relax}}
\newcommand\rsetminusaux[2]{\mspace{-4mu}
  \raisebox{\rsmraise{#1}\depth}{\rotatebox[origin=c]{-20}{$#1\smallsetminus$}}
 \mspace{-4mu}
}
\newcommand\rsmraise[1]{%
  \ifx#1\displaystyle .8\else
    \ifx#1\textstyle .8\else
      \ifx#1\scriptstyle .6\else
        .45%
      \fi
    \fi
  \fi}

\title{Algebraic analysis on scalar generalized Verma modules of Heisenberg parabolic type I.: $A_n$-series}

\author{Libor Křižka, Petr Somberg}

\AtEndDocument{\bigskip{\footnotesize%
  (L.\,Křižka) \textsc{Mathematical Institute of Charles University, Sokolovská 83, 180\,00 Praha 8, Czech Republic} \par
  \textit{E-mail address}: \texttt{krizka@karlin.mff.cuni.cz} \par
  \addvspace{\medskipamount}
  (P.\,Somberg) \textsc{Mathematical Institute of Charles University, Sokolovská 83, 180\,00 Praha 8, Czech Republic} \par
  \textit{E-mail address}: \texttt{somberg@karlin.mff.cuni.cz}
}}

\begin{document}
\date{}
\maketitle

\begin{abstract}

In the present article, we combine some techniques in the harmonic analysis together with
the geometric approach given by modules over sheaves of rings of twisted differential
operators ($\mcal{D}$-modules), and reformulate the composition series and
branching problems for objects
in the Bernstein-Gelfand-Gelfand parabolic category $\mcal{O}^\mfrak{p}$ geometrically
realized on certain orbits in the generalized flag manifolds.
The general framework is then applied to the scalar generalized Verma modules supported
on the closed
Schubert cell of the generalized flag manifold $G/P$ for $G=\SL(n+2,\C)$ and $P$ the Heisenberg
parabolic subgroup, and the algebraic analysis gives a complete classification of
$\mfrak{g}'_r$-singular vectors for all $\mfrak{g}'_r=\mfrak{sl}(n-r+2,\C)\,\subset\,
\mfrak{g}=\mfrak{sl}(n+2,\C)$, $n-r > 2$.
A consequence of our results is that we classify $\SL(n-r+2,\C)$-covariant differential
operators acting on homogeneous line bundles over the complexification of the odd dimensional
CR-sphere $S^{2n+1}$ and valued in homogeneous vector bundles over the complexification of
the CR-subspheres $S^{2(n-r)+1}$.


\end{abstract}

\thispagestyle{empty}

\tableofcontents


\section*{Introduction}
\addcontentsline{toc}{section}{Introduction}

Let $\mfrak{g}' \subset \mfrak{g}$ be a pair of reductive Lie algebras
of reductive Lie groups $G'\subset G$, $P\subset G$ and $P'\subset G'$
their parabolic subgroups, $P'=P\cap G'$, $\mfrak{p}$ and $\mfrak{p}'$
the parabolic Lie algebras of $P$ and $P'$, respectively.
The subject of the present paper is the study of the
$\mfrak{g}'$-composition series of objects
in the Bernstein-Gelfand-Gelfand (BGG for short) parabolic category
$\mcal{O}^\mfrak{p}$ associated to the pair $(\mfrak{g}, \mfrak{p})$.
The question of composition series for an object in $\mcal{O}^\mfrak{p}$
and a reductive Lie subalgebra $\mfrak{g}'$ can be rephrased
in terms of the existence of singular vectors or homomorphisms
between generalized Verma modules, and belongs to branching problems in
representation theory. A systematic approach to such problems was initiated
in \cite{Kobayashi1994, Kobayashi1998a, Kobayashi1998b, Kobayashi-Pevzner2013, koss}. Despite the fact that the branching
problem can have a continuous spectrum, \cite{Kobayashi2012}, there is a
wide class of discrete multiplicity free branching problems, \cite{Kobayashi-Pevzner2013}.
This property is fulfilled by our example discussed in the second half of the article.

There are two, mutually dual, motivations for writing the present article.
The first motivation is to develop a tool which allows a complete classification scheme
for a wide class of modules of interest in the BGG parabolic category
$\mcal{O}^\mfrak{p}$. Secondly, we would like to get an explicit description of
$G'$-covariant differential operators acting between homogeneous vector bundles
on generalized flag manifolds $G/P$ and $G'/P'$. The $G'$-covariant differential
operators differentiate sections of a homogeneous vector bundle
on $G/P$ and then restrict them on the submanifold $G'/P'\subset G/P$.
The existence of covariant differential operators is guaranteed by the discrete
decomposability for the branching problem in the BGG parabolic category $\mcal{O}^\mfrak{p}$, cf.\
\cite{Kobayashi1994, Kobayashi2007} for a general theory or \cite[Section 3]{koss} for
a review. The multiplicity-free condition of the branching problem then
assures uniqueness of covariant differential operators.

In \cite[Introduction]{koss}, two problems were formulated. For
$\mfrak{g}' \subset \mfrak{g}$ as above and
$M_\mfrak{p}^\mfrak{g}(\lambda)$ the generalized Verma
$\mfrak{g}$-module induced from an irreducible finite-dimensional
$\mfrak{p}$-module $V_\lambda$ with highest weight $\lambda$,
find the generators (or precise positions)
of irreducible $\mfrak{g}'$-submodules
in $M_\mfrak{p}^\mfrak{g}(\lambda)$ and find the associated Jordan-H\"older
series for $\mfrak{g}'$. Every irreducible $\mfrak{g}'$-submodule of an object in
$\mcal{O}^\mfrak{p}$ contains a singular vector and vice versa, every singular
vector generates a finite length $\mfrak{g}'$-submodule. In particular,
the singular vectors are responsible for the description of composition series
in the Grothendieck group $K(\mcal{O}^\mfrak{p})$ of the BGG parabolic category
$\mcal{O}^\mfrak{p}$.

A class of simple examples related to orthogonal Lie algebras $\mfrak{g}$, $\mfrak{g}'$
and parabolic subalgebras with commutative nilradicals is discussed in \cite{Kobayashi-Pevzner2013, koss}.
The technique to find the $\mfrak{g}'$-singular vectors is based on
the geometric realization of generalized Verma modules as twisted $\mcal{D}$-modules
supported on the closed (point) Schubert cell in $G/P$. Then the complicated
algebro-combinatorial problem
characterizing the singular vectors and thus the branching problem is converted,
by algebraic Fourier transform,
to a system of partial differential equations acting on the algebra of polynomials on
the opposite nilradical. This step allows to determine completely its space of solutions.

The aim of our article is twofold. We follow the seminal work of Beilinson and
Bernstein, \cite{Beilinson-Bernstein1981}, and generalize the approach
and results achieved in \cite{Kobayashi-Pevzner2013, koss}, towards the realization and
characterization of $(\mfrak{g}, N)$-modules in the BGG parabolic category
$\mcal{O}^\mfrak{p}$ as $N$-equivariant twisted $\mcal{D}$-modules supported
on $N$-orbits in $G/P$ ($N\subset P$ is the nilpotent subgroup of $P$).
The situation in \cite{Kobayashi-Pevzner2013, koss} is then recovered by considering
the generalized Verma $\mfrak{g}$-modules supported on the closed
point $N$-orbit $\{eP\} \subset G/P$.
In a forthcoming work, we employ the partial algebraic  Fourier transform and treat
the branching
problems analogous to ones formulated above for twisted generalized Verma $\mfrak{g}$-modules.
The language of $\mcal{D}$-modules is an indispensable tool to deal with
the parabolic subalgebras with non-commutative nilradials, and allows a natural implementation
of the symmetrization map between the universal
enveloping algebra and the symmetric algebra associated to the pair $(\mfrak{g}, \mfrak{p})$.
We would like to remark that even
the construction of non-standard homomorphisms between generalized Verma modules
in the case $\mfrak{g}'=\mfrak{g}$ is a difficult task, not known in
general. In the second part of our
article, we exploit the results of the first part and give a complete classification
of singular vectors (mentioned in the first paragraph of the Introduction) in a
particularly important CR (CR stands for Cauchy-Riemann or
Complex-Real) case of generalized flag manifolds and their CR flag submanifolds.
This means that we consider the case of complex Lie algebras
$\mfrak{g}=\mfrak{sl}(n+2,\C)$, $\mfrak{g}'_r=\mfrak{sl}(n-r+2,\C)$,
such that the nilradicals of their
parabolic subalgebras are isomorphic to Heisenberg Lie algebras and the
objects in $\mcal{O}^\mfrak{p}$ are the scalar
generalized Verma modules.

The content of our article is as follows. In Section \ref{sec:D-modules}, we briefly recollect
the notions of sheaves of rings of twisted differential operators $\mcal{D}_X^\lambda$
and equivariant
$\mcal{D}_X^\lambda$-modules. Focusing on the case of a generalized flag manifold
$X=G/P$, we describe the homomorphism from $\mfrak{g}$
into $\Gamma(X,\mcal{D}_X^\lambda)$ for a general pair $(\mfrak{g}, \mfrak{p})$.
In Section \ref{sec:Verma modules}, we discuss the geometric realization
of generalized Verma modules $M_\mfrak{p}^\mfrak{g}(\lambda)$ as
$\mcal{D}_X^\lambda$-modules supported on
the closed $N$-orbit $\{eP\}$ in $G/P$. An important device is the symmetrization map from the
polynomial algebra $S(\widebar{\mfrak{u}})$ to the universal enveloping algebra
$U(\widebar{\mfrak{u}})$ of the opposite nilradical $\widebar{\mfrak{u}}$, and the algebraic Fourier
transform of $M_\mfrak{p}^\mfrak{g}(\lambda)$. In the last part of this
section we introduce the key structure of $\mfrak{g}'$-singular vectors, and explain
our approach leading to the complete classification results in particular examples.
In Section \ref{sec:Heisenberg parabolic}, we apply all steps developed in Section \ref{sec:D-modules}
and Section \ref{sec:Verma modules} to the simplest
examples going beyond the Hermitean symmetric spaces (characterized by parabolic subalgebras
with commutative nilradicals).
Namely, we focus on the first in the row, and as for the applications the most important,
example of $A_n$-series
of Lie algebras and their parabolic subalgebras with Heisenberg nilradicals together
with $1$-dimensional
inducing representations. Using the tool of harmonic analysis, conveniently
organized by the Fischer decomposition of the symmetric algebra of $\widebar{\mfrak{u}}$
with respect to $\mfrak{sl}(n-r,\C) \oplus \mfrak{sl}(r,\C)$, $n-r > 2$,
we carry on the complete classification and explicit description of
$\mfrak{g}'_r=\mfrak{sl}(n-r+2,\C)$-singular vectors. In Section \ref{sec:CR operators}, we discuss the
construction of $G'$-covariant differential operators on the complexification of the
CR-sphere. Finally, in Appendix \ref{app:Fischer decompostion} we summarize for
reader's convenience the structure
of the Fischer decomposition for the action of $\mfrak{sl}(n,\C)$ on the polynomial algebra
$\C[(\C^n)^* \oplus \C^n]$.

Let us emphasize that due to the length of the article we decided to present
here just the results describing the composition and branching problems in the
Grothendieck group $K(\mcal{O}^\mfrak{p})$ of the BGG parabolic category
$\mcal{O}^\mfrak{p}$. The finer
property, determined by the Jordan-H\"older series, is intimately related to the
factorization properties of special polynomials responsible for the structure
of singular vectors (cf.\ \cite{koss} for examples related to the Gegenbauer
polynomials) and will appear in a forthcoming work. The cases of other
series of simple Lie algebras with Heisenberg parabolic subalgebras as well as the
cases of twisted Verma modules supported on the other $N$-orbits are also the
subject of our forthcoming work.

To summarize our achievements, the results in the present article contribute both
to the pure representation theoretical problems on the locus of reducibility (and their
generalizations for a class of reductive Lie subalgebras) of generalized Verma modules
(cf.\ \cite{Gyoja1994, Barchini-Kable-Zierau2008} for the formulation of Gyoja's conjectures), and to applications
related to the classification and explicit construction of covariant differential
operators for CR-geometries and their CR-subgeometries.


\section{${\mcal D}$-modules on generalized flag manifolds}
\label{sec:D-modules}

\subsection{Sheaves of rings of twisted differential operators}
\label{sec:Twisted differential}

As it follows from the seminal work of Beilinson and Bernstein (\cite{Beilinson-Bernstein1981}),
representation theory of Lie groups and algebras can be studied with the aid of the geometry
of $\mcal{D}$-modules on generalized flag manifolds. In this section we review several basic
notations and introduce some conventions useful to our further considerations.
More detailed information on $\mcal{D}$-modules can be found in
\cite{Kashiwara1989, Kashiwara-book, Hotta-book}.


Let $X$ be a complex manifold, $\mcal{O}_X$ the sheaf of holomorphic functions and
$\Theta_X$ the sheaf of holomorphic vector fields on $X$.
\medskip

\definition{A sheaf of rings of twisted differential operators on $X$
is a sheaf of rings
$\mcal{D}$ on $X$, equipped with an increasing exhausting filtration
$\{F_m\mcal{D}\}_{m \in \Z}$ of $\mcal{D}$ ($F_m\mcal{D}=0$ for $m<0$)
and a morphism of sheaves of rings $i \colon \mcal{O}_X \rarr \mcal{D}$
satisfying
\begin{enum}
  \item[(1)] the constant sheaf
	$\C_X$ is contained in the center of $\mcal{D}$, i.e.\ $\mcal{D}$ is a $\C_X$-algebra,
  \item[(2)] $i \colon \mcal{O}_X \rarr F_0\mcal{D}$ is an isomorphism of sheaves of rings,
  \item[(3)] $[a,i(f)] \in F_0\mcal{D}(U)$ for all $a \in F_1\mcal{D}(U)$ and $f \in \mcal{O}_X(U)$, and the morphism $\sigma_1 \colon F_1\mcal{D} \rarr \Theta_X$ of left $\mcal{O}_X$-modules defined by $i(\sigma_1(a)(f))=[a,i(f)]$
      induces the isomorphism $\sigma|_{\gr^F_1\!\mcal{D}} \colon \gr^F_1\!\mcal{D} \rarr \Theta_X$ of left $\mcal{O}_X$-modules,
  \item[(4)] $\smash{(\sigma|_{\gr^F_1\!\mcal{D}})^{-1}}$ defines the isomorphism $\sigma^{-1} \colon S_{\mcal{O}_X}\!(\Theta_X) \rarr \gr^F\!\mcal{D}$ of graded sheaves of rings, where $S_{\mcal{O}_X}\!(\Theta_X)$ is the symmetric algebra of $\Theta_X$ over $\mcal{O}_X$.
\end{enum}}

The simplest example of a sheaf of rings of twisted differential operators on
$X$ is the sheaf of rings of differential operators $\mcal{D}_X$.
If $\mcal{D}$ is a sheaf of rings of twisted differential operators on $X$ and
 $\mcal{L}$ an invertible $\mcal{O}_X$-module, then
$\mcal{L} \otimes_{\mcal{O}_X}\! \mcal{D} \otimes_{\mcal{O}_X}\! \mcal{L}^{-1}$
is a sheaf of rings of twisted differential operators on $X$ as well.


There are two basic operations on $\mcal{D}$-modules called the inverse image and
the direct image. We shall briefly introduce just the operation of the direct image
which will be important later on. Let $f \colon X \rarr Y$ be a morphism of complex
manifolds and let $\mcal{D}$ be a sheaf of rings of twisted differential operators on $Y$.
Then we denote by $f^\sharp\mcal{D}$ the sheaf of rings of twisted differential operators on $X$
given by the pull-back of $\mcal{D}$, in the case $\mcal{D}=\mcal{D}_Y$ we have
$f^\sharp\mcal{D}_Y \simeq \mcal{D}_X$. Denoting $\mcal{M}$ a left
$f^\sharp \mcal{D}$-module, its direct image is the left $\mcal{D}$-module
defined by $f_*(\mcal{D}_{Y \larr X} \otimes_{f^\sharp\mcal{D}} \mcal{M})$, where
the $(f^{-1}\mcal{D},f^\sharp\mcal{D})$-bimodule $\mcal{D}_{Y \larr X}$ is given by
\begin{align}
 \mcal{D}_{Y \larr X} = f^{-1}\mcal{D} \otimes_{f^{-1}\mcal{O}_Y} \omega_{X/Y}.
\end{align}
Here $\omega_{X/Y}= \omega_X \otimes_{f^{-1}\mcal{O}_X} f^{-1}\omega_Y$, and
$\omega_X, \omega_Y$ are the canonical sheaves of $X, Y$, respectively.


If $X$ carries an action of a complex Lie group $G$, then there is a notion of $G$-equivariant
sheaves of rings of twisted differential operators on $X$. We will not give a precise
definition of this structure here, because we will not need it in its full generality.


Let $X$ be a homogeneous space for a complex Lie group $G$. Then we get a principal
$H$-bundle $p \colon G \rarr X$, where $H$ is the stabilizer of a point in $X$. We
denote by $\mfrak{g}$ and $\mfrak{h}$ the Lie algebras of the Lie groups $G$ and $H$,
respectively, and by $\Hom_H(\mfrak{h},\C)$ the vector space of $H$-equivariant
homomorphisms from $\mfrak{h}$ to $\C$.
The left action of $G$ on $X$ induces the Lie algebra homomorphism
\begin{align}
  L_X \colon \mfrak{g} \rarr \Gamma(X,\Theta_X).
\end{align}
The constant sheaf of algebras $\mfrak{g}_X$ on $X$ naturally acts on $\mcal{O}_X$ and
by the left multiplication on $U(\mfrak{g})_X$, where $U(\mfrak{g})$ is the universal
enveloping algebra of $\mfrak{g}$ and $U(\mfrak{g})_X$ the corresponding constant
sheaf on $X$. Then $\mfrak{g}_X$ naturally acts on
$\mcal{U}_X(\mfrak{g})=\mcal{O}_X\! \otimes_{\C_X}\!U(\mfrak{g})_X$,
and extends to the action of $U(\mfrak{g})_X$ on $\mcal{U}_X(\mfrak{g})$. Because
$\mcal{O}_X$ acts on $\mcal{U}_X(\mfrak{g})$ as well, $\mcal{U}_X(\mfrak{g})$ is the
sheaf of rings on $X$.

Any element $\lambda \in \Hom_H(\mfrak{h},\C)$ gives rise to the $1$-dimensional
representation $\C_\lambda$ of $\mfrak{h}$, defined by
\begin{align}
  Av=\lambda(A)v, \quad A \in \mfrak{h},\, v \in \C.
\end{align}
On the other hand, any $H$-equivariant homomorphism $\lambda \colon \mfrak{h} \rarr \C$
yields the $G$-equivariant morphism $f_\lambda \colon \mcal{V}_X(\mfrak{h}) \rarr \mcal{V}_X(\C)$
of $\mcal{O}_X$-modules associated to the representations $\mfrak{h}$
and $\C$ of $H$. Since
${\textstyle \sum}_{s \in \mcal{V}_X\!(\mfrak{h})}\mcal{U}_X(\mfrak{g})(s-f_\lambda(s))$
is a two-sided ideal of $\mcal{U}_X(\mfrak{g})$, we obtain a $G$-equivariant sheaf of
rings of twisted differential operators
\begin{align}
  \mcal{D}_X(\lambda)=\mcal{U}_X(\mfrak{g})/{\textstyle \sum}_{s \in \mcal{V}_X\!(\mfrak{h})} \mcal{U}_X(\mfrak{g})(s-f_\lambda(s)) \label{eq:twisted sheaf realization}
\end{align}
on $X$, and any $G$-equivariant sheaf of rings of twisted differential operators on $X$
is isomorphic to $\mcal{D}_X(\lambda)$ for a uniquely determined $\lambda \in \Hom_H(\mfrak{h},\C)$.


There is an alternative description of $\mcal{D}_X(\lambda)$ given as follows.
We start with the sheaf of rings of differential operators $\mcal{D}_G$
on $G$ equipped with the natural structure of an $H$-equivariant
(with respect to the right action of $H$ on $G$) sheaf of rings of
twisted differential operators, and consider the left $\mcal{D}_G$-module
$\mcal{D}_G(\mfrak{h},\lambda)$ given by
\begin{align}
  \mcal{D}_G(\mfrak{h},\lambda)=\mcal{D}_G/{\textstyle \sum_{A \in \mfrak{h}}}\mcal{D}_G(R_G(A)+
	\lambda(A)),\quad \lambda \in \Hom_H(\mfrak{h},\C),
\end{align}
for $R_G \colon \mfrak{g} \rarr \Gamma(G,\Theta_G)$ the Lie algebra
homomorphism induced by the right action of $G$ on itself. Then we get
\begin{align}
  \mcal{D}_X(\lambda) \simeq (p_*\mcal{D}_G(\mfrak{h},\lambda))^H,
\end{align}
where $(p_*\mcal{D}_G(\mfrak{h},\lambda))^H$ means the subsheaf of $H$-invariants
of the direct image $p_*\mcal{D}_G(\mfrak{h},\lambda)$.

Since $\mcal{D}_X(\lambda)$ is a $G$-equivariant sheaf of rings of twisted differential operators on $X$, we obtain the Lie algebra homomorphism
\begin{align}
  \alpha_{\mcal{D}_X\!(\lambda)} \colon \mfrak{g} \rarr \Gamma(X,\mcal{D}_X(\lambda))
\end{align}
described in the following way. For any $Y \in \mfrak{g}$, we define the morphism
$\varphi_Y \colon \mcal{D}_G \rarr \mcal{D}_G$ of left $\mcal{D}_G$-modules by
\begin{align}
  P \mapsto PL_G(Y)|_U,\quad P \in \mcal{D}_G(U).
\end{align}
Because $[L_G(Y),R_G(A)]=0$ for all $A \in \mfrak{g}$, the morphism $\varphi_Y$ gives
rise to the $H$-equivariant morphism
$\psi_Y \colon \mcal{D}_G(\mfrak{h},\lambda) \rarr \mcal{D}_G(\mfrak{h},\lambda)$
of left $\mcal{D}_G$-modules. Then we set
\begin{align}
\alpha_{\mcal{D}_G\!(\lambda)}(Y)=\psi_Y(u)
\end{align}
with $u = 1\ {\rm mod}\ {\textstyle \sum}_{A \in \mfrak{h}}\mcal{D}_G(R_G(A)+\lambda(A))$.
Furthermore, if we denote by $\mcal{D}_X(\lambda)^{\rm op}$ the sheaf of rings opposite to
$\mcal{D}_X(\lambda)$, then $\mcal{D}_X(\lambda)^{\rm op}$ is also a $G$-equivariant sheaf
of rings of twisted differential operators on $X$ and we have
\begin{align}
  \mcal{D}_X(\lambda)^{\rm op} \simeq \mcal{D}_X(2\rho-\lambda),
\end{align}
where the element $\rho \in \Hom_H(\mfrak{h},\C)$ is given by
\begin{align}
  \rho(A) = - {1 \over 2} \tr_{\mfrak{g}/\mfrak{h}}\ad(A),\quad
	A \in \mfrak{h} . \label{eq:rho vector}
\end{align}


\subsection{Generalized flag manifolds}
\label{sec:Flag manifolds}

Let $G$ be a complex semisimple Lie group, $H$ be a maximal torus of $G$ and
$\mfrak{g}$ and $\mfrak{h}$ be the Lie algebras of the Lie groups $G$ and $H$,
respectively. Then $\mfrak{h}$ is a
Cartan subalgebra of $\mfrak{g}$. We denote by $\Delta$ the root system of $\mfrak{g}$
with respect to $\mfrak{h}$, $\Delta^+$ the positive root system in $\Delta$ and $\Pi\subset\Delta^+$
the set of simple roots. Furthermore, we associate to the positive root system
$\Delta^+$ the nilpotent Lie subalgebras
\begin{align}
  \mfrak{n}= \bigoplus_{\alpha \in \Delta^+} \mfrak{g}_\alpha \qquad \text{and} \qquad \widebar{\mfrak{n}}= \bigoplus_{\alpha \in \Delta^+} \mfrak{g}_{-\alpha},
\end{align}
and the solvable Lie subalgebras
\begin{align}
  \mfrak{b}= \mfrak{h} \oplus \mfrak{n} \qquad \text{and} \qquad \widebar{\mfrak{b}}= \mfrak{h} \oplus \widebar{\mfrak{n}}
\end{align}
of $\mfrak{g}$. The Lie algebras $\mfrak{b}$ and $\widebar{\mfrak{b}}$ are the (standard and opposite standard) Borel
subalgebras of $\mfrak{g}$.

Let us consider a subset $\Sigma$ of $\Pi$ and denote by $\Delta_\Sigma$ the root subsystem in $\mfrak{h}^*$ generated
by $\Sigma$. Then the standard parabolic subalgebra $\mfrak{p}$ of $\mfrak{g}$ associated to $\Sigma$ is defined by
\begin{align}
  \mfrak{p} = \mfrak{l} \oplus \mfrak{u},  \label{eq:levi decoposition}
\end{align}
where the reductive Levi factor $\mfrak{l}$ of $\mfrak{p}$ is defined through
\begin{align}
  \mfrak{l}= \mfrak{h} \oplus \bigoplus_{\alpha \in \Delta_\Sigma} \mfrak{g}_\alpha
\end{align}
and the nilradical $\mfrak{u}$ of $\mfrak{p}$ and the opposite nilradical $\widebar{\mfrak{u}}$ are given by
\begin{align}\label{parnilradical}
  \mfrak{u}= \bigoplus_{\alpha \in \Delta^+ \rsetminus \Delta_\Sigma^+}\mfrak{g}_\alpha \qquad \text{and} \qquad \widebar{\mfrak{u}}= \bigoplus_{\alpha \in \Delta^+ \rsetminus \Delta_\Sigma^+} \mfrak{g}_{-\alpha}.
\end{align}
We say that a weight $\lambda \in \mfrak{h}^*$ is $\mfrak{p}$-dominant integral provided
$\smash{{2(\lambda,\alpha) \over (\alpha,\alpha)}} \in \N_0$  for all $\alpha \in \Sigma$,
and we denote the set of all $\mfrak{p}$-dominant integral weight by $\Lambda^+(\mfrak{p})$.
Furthermore, let $P$ be a parabolic subgroup of $G$ with the Lie algebra $\mfrak{p}$. Then
$p \colon G \rarr G/P$ is a principal $P$-bundle and $X=G/P$ is called the generalized flag
manifold associated of $G$.

We define the subgroups $\widebar{U}$ and $N$ of $G$ to be the image of $\widebar{\mfrak{u}}$
and $\mfrak{n}$ under the exponential map $\exp \colon \mfrak{g} \rarr G$, respectively.
Moreover, the mapping
\begin{align}
  \exp \colon \widebar{\mfrak{u}} \rarr \widebar{U}
\end{align}
is a diffeomorphism and $\widebar{U}$ is a closed nilpotent subgroup of $G$. Therefore, for
$g \in G$ the mapping $f_g \colon \widebar{U} \rarr X$ given by
\begin{align}
  f_g(n)=p(gn)=gnP
\end{align}
is an open embedding, and the generalized flag manifold $X$ is covered by open subsets
$U_g=p(g\widebar{U}) \subset X$, i.e.\ we have
\begin{align}
  X = \bigcup_{g \in G} U_g.
\end{align}
Hence, we get the atlas $\{(U_g,u_g)\}_{g \in G}$ on $X$, where
$u_g \colon U_g \rarr \widebar{\mfrak{u}}$ is defined by
\begin{align}
  u_g^{-1}=f_g\circ \exp. \label{eq:atlas flag manifold}
\end{align}
Now, let $s_g \colon U_g \rarr G$ be the local section of the principal $P$-bundle
 $p \colon G \rarr X$ defined through
\begin{align}
  s_g(x)=gf_g^{-1}(x)=gs_e(g^{-1}.x),\quad x \in U_g,\, g \in G .\label{eq:trivialization principal bundle}
\end{align}
The local section $s_g \colon U_g \rarr G$ gives the trivialization of the principal $P$-bundle
$p \colon G \rarr X$ over $U_g$. Therefore, we get the isomorphism
\begin{align}
  \bfig
  \square<700,500>[p^{-1}(U_g)`U_g\times P`U_g`U_g;\varphi_g`p`p_1`{\rm id}_{U_g}]
  \efig
\end{align}
of principal $P$-bundles over $U_g$, where the mapping $\varphi_g \colon p^{-1}(U_g) \rarr U_g \times P$ is defined by
\begin{align}
  \varphi_g(h)=(p(h),s_g(p(h))^{-1}h)\quad \text{for all}\quad h \in p^{-1}(U_g) .\label{eq:principal P-bundle triv}
\end{align}
Denoting by $\boxtimes$ the outer tensor product, we obtain an isomorphism
\begin{align}
  \Phi_g \colon \mcal{D}_G|_{p^{-1}(U_g)} \riso \mcal{D}_{X \times P}|_{p_1^{-1}(U_g)} \riso \mcal{D}_X|_{U_g} \boxtimes \mcal{D}_P
\end{align}
of $P$-equivariant sheaves of rings of twisted differential operators and consequently the isomorphism
\begin{align}
  (p_*\mcal{D}_G(\mfrak{p},\lambda)|_{U_g})^P \riso \mcal{D}_X|_{U_g} \boxtimes \Gamma\big(P,\mcal{D}_P/{\textstyle \sum_{A \in \mfrak{p}}}\mcal{D}_P(R_P(A)+\lambda(A))\big)^P
\end{align}
of sheaves of rings of twisted differential operators on $U_g$. Since
$\mcal{D}_P/{\textstyle \sum}_{A \in \mfrak{p}}\mcal{D}_P(R_P(A)+\lambda(A))$
is isomorphic to $\mcal{O}_P$ as an $\mcal{O}_P$-module, we have the isomorphism
\begin{align}
  j_{s_g} \colon \mcal{D}_X(\lambda)|_{U_g} \riso (p_*\mcal{D}_G(\mfrak{p},\lambda)|_{U_g})^P \riso \mcal{D}_X|_{U_g} \label{eq:twisted sheaf trivialization}
\end{align}
of sheaves of rings of twisted differential operators on $U_g$.
\medskip

\proposition{Let $\lambda \in \Hom_P(\mfrak{p},\C)$. Then the Lie algebra homomorphism
\begin{align}
  \alpha_{\mcal{D}_X\!(\lambda)} \colon \mfrak{g} \rarr \Gamma(X,\mcal{D}_X(\lambda))
\end{align}
is given by
\begin{align}
  (j_{s_g} \circ \alpha_{\mcal{D}_X\!(\lambda)})(Y)= \pi^\lambda_g(Y),
\end{align}
where $\pi^\lambda_g(Y) \in \mcal{D}_X(U_g)$ is defined by
\begin{align}
  (\pi^\lambda_g(Y)f)(x)=(L_X(Y)f)(x)+\lambda((\Ad(s_g(x)^{-1})Y)_{\mfrak{p}})f(x)
\end{align}
for all $x \in V$, $f \in \mcal{O}_X(V)$ and $V \subset U_g$.}

\proof{For $Y\in \mfrak{g}$, we define the morphism $\varphi_Y \colon \mcal{D}_G \rarr \mcal{D}_G$ of left $\mcal{D}_G$-modules by
\begin{align*}
  P \mapsto PL_G(Y)|_U
\end{align*}
for all $P \in \mcal{D}_G(U)$. Since we have $[L_G(Y),R_G(A)]=0$ for all $A\in \mfrak{g}$, the morphism $\varphi_Y$ gives rise to a
$P$-equivariant morphism $\psi_Y \colon \mcal{D}_G(\mfrak{p},\lambda) \rarr \mcal{D}_G(\mfrak{p},\lambda)$ of left $\mcal{D}_G$-modules. Then we have
$\alpha_{\mcal{D}_X\!(\lambda)}(Y)=\psi_Y(u) \in \Gamma(X,\mcal{D}_X(\lambda)) \simeq \Gamma(X,p_*\mcal{D}_G(\mfrak{p},\lambda))^P$,
where $u = 1\ {\rm mod}\ {\textstyle \sum}_{A \in \mfrak{p}}\mcal{D}_G(R_G(A)+\lambda(A))$.

In the next step we compute $j_{s_g}(\psi_Y(u)|_{U_g}) \in \mcal{D}_X(U_g)$.
Let $\varphi_g \colon p^{-1}(U_g) \rarr U_g \times P$ be the isomorphism of
principal $P$-bundles over $U_g$ given by \eqref{eq:principal P-bundle triv}.
Then the corresponding element $\Phi_g(L_G(Y)|_{p^{-1}(U_g)}) \in \mcal{D}_X|_{U_g} \boxtimes \mcal{D}_P$ is given by the formula
\begin{align*}
\Phi_g(L_G(Y)|_{p^{-1}(U_g)})f= ((L_G(Y)|_{p^{-1}(U_g)})(f\circ \varphi_g))\circ \varphi_g^{-1}
\end{align*}
for all $f \in \mcal{O}_{X\times P}(U_g\times P)$. Hence we can write
\begin{align*}
 (\Phi_g(L_G(Y))f)(x,p) &={{\rm d}\over {\rm d}t}_{|t=0}
f(\exp(-tY).x,s_g(\exp(-tY).x)^{-1}\exp(-tY)s_g(x)p)\\
  &= {{\rm d}\over {\rm d}t}_{|t=0}f(\exp(-tY).x,p) + {{\rm d}\over {\rm d}t}_{|t=0}
	f(x,s_g(\exp(-tY).x)^{-1}\exp(-tY)s_g(x)p).
\end{align*}
The first term can be rewritten as
\begin{align*}
  {{\rm d}\over {\rm d}t}_{|t=0}f(\exp(-tY).x,p) = ((L_X(Y)|_{U_g} \boxtimes 1_P)f)(x,p).
\end{align*}
For the second term we can write
\begin{align*}
  {{\rm d}\over {\rm d}t}_{|t=0}f(x,s_g(\exp(-tY).x)^{-1}\exp(-tY)s_g(x)p)
	= {{\rm d}\over {\rm d}t}_{|t=0}f(x,pg(t)),
\end{align*}
where
\begin{align*}
  g(t)&=p^{-1}s_g(\exp(-tY).x)^{-1}\exp(-tY)s_g(x)p \\
&=  p^{-1}s_e(g^{-1}\exp(-tY).x)^{-1}g^{-1}\exp(-tY)gs_e(g^{-1}x)p.
\end{align*}
Further, if $g \in \widebar{U}P$, then there exist uniquely determined elements $\widebar{U}(g) \in \widebar{U}$ and  $P(g) \in P$ such that $g=\widebar{U}(g)P(g)$.
Since $g(t) \in P$, we can write
\begin{align*}
  g(t)&=P(p^{-1}s_e(g^{-1}\exp(-tY).x)^{-1}g^{-1}\exp(-tY)gs_e(g^{-1}.x)p) \\
  &= p^{-1}P(s_e(g^{-1}\exp(-tY).x)^{-1}g^{-1}\exp(-tY)gs_e(g^{-1}.x))p \\
  &= p^{-1}P(g^{-1}\exp(-tY)gs_e(g^{-1}.x))p \\
  &= p^{-1}P(s_e(g^{-1}.x)^{-1}g^{-1}\exp(-tY)gs_e(g^{-1}.x))p \\
  &= p^{-1}P(s_g(x)^{-1}\exp(-tY)s_g(x))p.
\end{align*}
This gives
\begin{align*}
  {{\rm d}\over {\rm d}t}_{|t=0}g(t)=-\Ad(p^{-1})(\Ad(s_g(x)^{-1})Y)_\mfrak{p},
\end{align*}
where $\Ad(s_g(x)^{-1})Y=(\Ad(s_g(x)^{-1})Y)_{\widebar{\mfrak{u}}}+(\Ad(s_g(x)^{-1})Y)_\mfrak{p}$ is the decomposition of the element $\Ad(s_g(x)^{-1})Y$ with respect to $\mfrak{g}=\widebar{\mfrak{u}} \oplus \mfrak{p}$. Therefore, we obtain
\begin{align*}
   {{\rm d}\over {\rm d}t}_{|t=0}f(x,s_g(\exp(-tY).x)^{-1}\exp(-tY)s_g(x)p)&
	= -(R_P(\Ad(p^{-1})(\Ad(s_g(x)^{-1})Y)_\mfrak{p})f)(x,p).
\end{align*}
But since the operator $Q \in \mcal{D}_G(p^{-1}(U_g))$ defined by
\begin{align*}
  (\Phi_g(Q)f)(x,p)&=(R_P(\Ad(p^{-1})(\Ad(s_g(x)^{-1})Y)_\mfrak{p})f)(x,p)+ \lambda(\Ad(p^{-1})(\Ad(s_g(x)^{-1})Y)_\mfrak{p})f(x,p)
\end{align*}
belongs to the ideal ${\textstyle \sum_{A \in \mfrak{p}}}\mcal{D}_G(R_G(A)+\lambda(A))$,
the $P$-equivariance of $\lambda \colon \mfrak{p} \rarr \C$ gives
\begin{align*}
  (\Phi_g(L_X(Y)|_{U_g})f)(x,p)&=((L_X(Y)|_{U_g} \boxtimes 1_P)f)(x,p)+ \lambda(\Ad(p^{-1})(\Ad(s_g(x)^{-1})Y)_\mfrak{p})f(x,p)\\
  & \quad + (\Phi_g(Q)f)(x,p) \\
  &=((L_X(Y)|_{U_g} \boxtimes 1_P)f)(x,p)+ \lambda((\Ad(s_g(x)^{-1})Y)_\mfrak{p})f(x,p)\\
  & \quad + (\Phi_g(Q)f)(x,p).
\end{align*}
If we put altogether, we obtain
\begin{align*}
  (j_{s_g}(\psi_Y(u)|_{U_g})f)(x)=(L_X(Y)f)(x)+\lambda((\Ad(s_g(x)^{-1})Y)_{\mfrak{p}})f(x)
\end{align*}
for all $f \in \mcal{O}_X(U_g)$, and the proof is complete.}
\medskip

Now, let us recall the atlas $\{(U_g,u_g)\}_{g\in G}$ on $X$ defined in
\eqref{eq:atlas flag manifold}, and let $(f_1,f_2,\dots,f_n)$ be a
basis of $\widebar{\mfrak{u}}$. Then for $x \in U_g$ we have
\begin{align}
u_g(x)=\sum_{i=1}^n u_g^i(x)f_i,
\end{align}
the functions $u_g^i \colon U_g \rarr \C$ are called the coordinate functions on $U_g$.
\medskip

\theorem{\label{prop:operator realization}Let $\lambda \in \Hom_P(\mfrak{p},\C)$. Then we have
\begin{align}
\pi^\lambda_g(Y)= -\sum_{i=1}^n\bigg[{\ad(u_g(x))e^{\ad(u_g(x))} \over e^{\ad(u_g(x))}-{\rm id}_{\widebar{\mfrak{u}}}}\,(e^{-\ad(u_g(x))}\!\Ad(g^{-1})Y)_{\widebar{\mfrak{u}}}\bigg]_i{\partial \over \partial u_g^i} + \lambda((e^{-\ad(u_g(x))}\!\Ad(g^{-1})Y)_\mfrak{p})
\end{align}
for all $Y\in \mfrak{g}$, where $[X]_i$ denotes the $i$-th coordinate of $X \in \widebar{\mfrak{u}}$ with respect to the basis $(f_1,f_2,\dots,f_n)$ of $\widebar{\mfrak{u}}$. In particular, we have
  \begin{align}
    \pi^\lambda_e(Y) = -\sum_{i=1}^n\bigg[{\ad(u_e(x)) \over e^{\ad(u_e(x))}-{\rm id}_{\widebar{\mfrak{u}}}}\,Y\bigg]_i{\partial \over \partial u_e^i}\quad \text{for}\quad Y \in \widebar{\mfrak{u}}
  \end{align}
   and
  \begin{align}
    \pi^\lambda_e(Y) = \sum_{i=1}^n\,[\ad(u_e(x))Y]_i{\partial \over \partial u_e^i} + \lambda(Y)
		\quad \text{for}\quad Y \in \mfrak{l} .
  \end{align}
}

\proof{Let $\{(U_g,u_g)\}_{g \in G}$ be the atlas on $X$ defined by the formula \eqref{eq:atlas flag manifold}.
Assuming $Y \in \mfrak{g}$, we can write by \eqref{eq:trivialization principal bundle}
\begin{align*}
  \Ad(s_g(x)^{-1})Y=\Ad(\exp(-u_g(x))g^{-1})Y=\Ad(\exp(-u_g(x)))\Ad(g^{-1})Y.
\end{align*}
The differential operator $L_X(Y)|_{U_g}$ can be rewritten in the local coordinates
$(u_g^1,u_g^2,\dots,u_g^n)$ on $U_g$ as follows. We have, for $V \subset U_g$ and
$f \in \mcal{O}_X(V)$,
\begin{align*}
  (L_X(Y)f)(x) &= {{\rm d}\over {\rm d}t}_{|t=0}f(\exp(-tY).x)=
	{{\rm d}\over {\rm d}t}_{|t=0}(f\circ u_g^{-1})(u_g(\exp(-tY).x)) \\
  & = \sum_{i=1}^n {{{\rm d} u_g^i(\exp(-tY).x) \over {\rm d}t}}_{|t=0} {\partial \over \partial u_g^i}_{|x}(f).
\end{align*}
If we denote $o=eP$, then we can write
\begin{align*}
  {{\rm d}\over {\rm d}t}_{|t=0}u_g(\exp(-tY).x) &=
	{{\rm d}\over {\rm d}t}_{|t=0}u_g(\exp(-tY)s_g(x).o) \\
  &= {{\rm d}\over {\rm d}t}_{|t=0}u_g(s_g(x)s_g(x)^{-1}\exp(-tY)s_g(x).o) \\
  &= {{\rm d}\over {\rm d}t}_{|t=0}u_g(gs_e(g^{-1}.x)s_g(x)^{-1}\exp(-tY)s_g(x).o) \\
  &= {{\rm d}\over {\rm d}t}_{|t=0}\exp^{-1}(\widebar{U}(s_e(g^{-1}.x)s_g(x)^{-1}\exp(-tY)s_g(x))) \\
  &= {{\rm d}\over {\rm d}t}_{|t=0}\exp^{-1}(s_e(g^{-1}.x) \widebar{U}(s_g(x)^{-1}\exp(-tY)s_g(x))) \\
  &= {{\rm d}\over {\rm d}t}_{|t=0}\exp^{-1}(\exp(u_g(x))\exp(-t(\Ad(s_g(x)^{-1})Y)_{\widebar{\mfrak{u}}})).
\end{align*}
The Baker-Campbell-Hausdorff formula for the nilpotent group $\widebar{U}$ gives
\begin{align*}
  \exp^{-1}(\exp(X)\exp(tZ))=X+{\ad X e^{\ad X} \over e^{\ad X}-{\rm id}_{\widebar{\mfrak{u}}}}\,tZ + t^2g(t),
\end{align*}
where $g(t)$ is a $\widebar{\mfrak{u}}$-valued polynomial in $t$, for all $X,Z \in \widebar{\mfrak{u}}$ and $t \in \R$, and so we obtain
\begin{align*}
 {{\rm d}\over {\rm d}t}_{|t=0}u_g(\exp(-tY).x) =
-{\ad(u_g(x)) e^{\ad(u_g(x))} \over e^{\ad(u_g(x))}-{\rm id}_{\widebar{\mfrak{u}}}}\,(\Ad(s_g(x)^{-1})Y)_{\widebar{\mfrak{u}}}.
\end{align*}
Hence we get
\begin{align*}
  L_X(Y)|_{U_g}= -\sum_{i=1}^n\bigg[{\ad(u_g(x))e^{\ad(u_g(x))} \over e^{\ad(u_g(x))}-{\rm id}_{\widebar{\mfrak{u}}}}\,(e^{-\ad(u_g(x))}\!\Ad(g^{-1})Y)_{\widebar{\mfrak{u}}}\bigg]_i{\partial \over \partial u_g^i}
\end{align*}
and thus we are done.}


\section{Generalized Verma modules}
\label{sec:Verma modules}

We remind the algebraic definition of a class of representations called
generalized Verma modules,
relying on the notation introduced in Section \ref{sec:Flag manifolds}.
The Bernstein-Gelfand-Gelfand  (BGG for short) category $\mcal{O}$ is
the full subcategory of the category of $U(\mfrak{g})$-modules whose objects are finitely
generated $U(\mfrak{g})$-modules, $\mfrak{h}$-semisimple and locally $\mfrak{n}$-finite.
Let $\mfrak{p}$ be a parabolic subalgebra of $\mfrak{g}$ containing
the Borel subalgebra $\mfrak{b}$, and $\mfrak{l}$ its reductive Levi factor, cf.\ \eqref{eq:levi decoposition}. The BGG parabolic category $\mcal{O}^\mfrak{p}$
is the full subcategory of $\mcal{O}$ whose objects are locally
$\mfrak{l}$-finite. The category $\mcal{O}^\mfrak{p}$ contains
the full subcategories $\mcal{O}_{\chi_\lambda}^\mfrak{p}$ whose objects
have generalized infinitesimal character $\chi_\lambda$ in the Harish-Chandra
parametrization, see \cite{Bernstein-Gelfand-Gelfand1971, Humphreys-book} for the detailed exposition.

We denote by $K(\mcal{O}^\mfrak{p})$ the Grothendieck group of the
abelian category $\mcal{O}^\mfrak{p}$ defined as the restricted product
$K(\mcal{O}^\mfrak{p})=\prod'_{[\lambda] \in \mfrak{h}^*\!/W} K(\mcal{O}_{\chi_\lambda}^\mfrak{p})$, where $W$ is the Weyl group of $\mfrak{g}$ with the affine action on $\mfrak{h}^*$ and $\prod'$ denotes the restricted product in which all up to
a countable number of components are zero.

Let $V_\lambda$ denote a finite-dimensional irreducible $\mfrak{l}$-module
with highest weight $\lambda$  regarded as ${\mfrak p}$-module with trivial action
of its nilradical $\mfrak{u}$, where $\lambda$ is a $\mfrak{p}$-dominant integral
weight $\lambda \in \Lambda^+(\mfrak{p})$.
The generalized Verma module $M^\mfrak{g}_\mfrak{p}(\lambda)$ is defined by
\begin{align}
M^\mfrak{g}_\mfrak{p}(\lambda)=
U(\mfrak{g})\otimes_{U(\mfrak{p})}\!V_\lambda,  \label{eq:generalized Verma module}
\end{align}
and $M^\mfrak{g}_\mfrak{p}(\lambda)$ is
called to be of scalar type provided $V_\lambda$ is a $1$-dimensional
$\mfrak{p}$-module. The highest weight modules $M^\mfrak{g}_\mfrak{p}(\lambda)$
are objects in $\mcal{O}^\mfrak{p}$, and any irreducible $U(\mfrak{g})$-module in
$\mcal{O}^\mfrak{p}$ is the quotient of some
$M^\mfrak{g}_\mfrak{p}(\lambda)$ by its maximal submodule.
Let us notice that in our conventions for the Harish-Chandra parametrization of maximal
ideals in the center $\mfrak{Z}(\mfrak{g})$ of the universal enveloping algebra $U(\mfrak{g})$, the trivial $U(\mfrak{g})$-module has
${\mfrak Z}({\mfrak g})$-infinitesimal character $\chi_{\rho_\mfrak{b}}$ and
$M^\mfrak{g}_\mfrak{p}(\lambda)$ has $\mfrak{Z}(\mfrak{g})$-infinitesimal
character $\chi_{\lambda+\rho_\mfrak{b}}$. Here $\rho_\mfrak{b} \in \mfrak{h}^*$ stands for the Weyl vector, i.e.\ $\rho_\mfrak{b}$ is the half-sum of positive roots of $\mfrak{g}$.


\subsection{Geometric realization of generalized Verma modules}

In the present section we rely on the notation of Section \ref{sec:Twisted differential}
and Section \ref{sec:Flag manifolds}.
Let us consider an $N$-orbit $Q$ in the generalized flag manifold $X$ and denote by
$i \colon Q \hookrightarrow X$ its embedding into $X$. Further, let $\tau$ be an irreducible
object in the category $\Con(i^\sharp\mcal{D}_X^\lambda,N)$ of $N$-equivariant
integrable $i^\sharp\mcal{D}_X^\lambda$-connections, i.e.\ the category of $N$-equivariant
$i^\sharp\mcal{D}_X^\lambda$-modules which are locally free $\mcal{O}_Q$-modules of
finite rank. Here $i^\sharp\mcal{D}_X^\lambda$
denotes the pull-back of $\mcal{D}_X^\lambda$ and
$\mcal{D}_X^\lambda= \mcal{D}_X(\lambda+\rho)$. Then to a given
pair $(Q,\tau)$ we attach an $N$-equivariant regular holonomic $\mcal{D}_X^\lambda$-module
$\mcal{I}(Q,\tau)$ defined as the maximal regular holonomic extension of
$i_*(\mcal{D}^\lambda_{X \larr Q} \otimes_{i^\sharp\mcal{D}^\lambda_X}\!\tau)|_{X\rsetminus \partial Q}$
for $\smash{\partial Q = \widebar{Q} \rsetminus Q}$, see \cite{Kashiwara1989}.

In the present article we are mostly interested in the generalized Verma modules, related
to the closed $N$-orbit $X_e=\{eP\}$. The case of other orbits realizing twisted generalized
Verma modules, will be discussed elsewhere. For the closed embedding
\begin{align}
i_e \colon X_e \hookrightarrow X
\end{align}
we have $i_e^\sharp\mcal{D}_X^\lambda \simeq \mcal{D}_{X_e}$ as $N$-equivariant sheaves of rings of twisted differential operators on $X_e$. Since $\mcal{O}_{X_e}$ is an irreducible object in the category $\Con(i_e^\sharp\mcal{D}_X^\lambda,N)$ and $X_e$ is a closed orbit, we get $\mcal{I}(X_e,\mcal{O}_{X_e})=i_{e*}(\mcal{D}^\lambda_{X \larr X_e} \otimes_{i_e^\sharp\mcal{D}^\lambda_X}\!\mcal{O}_{X_e})$.
\medskip

\proposition{Let $\lambda \in \Hom_P(\mfrak{p},\C)$. Then we have
\begin{align}
\Gamma(X,\mcal{I}(X_e,\mcal{O}_{X_e})) \simeq M^\mfrak{g}_\mfrak{p}(\lambda-\rho)
\end{align}
as $(\mfrak{g},N)$-modules.\label{lem:Verma module}}

\proof{Since $i_e \colon X_e \rarr X$ is a closed embedding, we obtain
\begin{align*}
  \mcal{I}(X_e,\mcal{O}_{X_e})=i_{e*}(\mcal{D}^\lambda_{X \larr X_e}\!\otimes_{i_e^\sharp\mcal{D}^\lambda_X}\!\mcal{O}_{X_e}).
\end{align*}
Further, from the definition of the $(i_e^{-1}\mcal{D}^\lambda_X,i_e^\sharp \mcal{D}^\lambda_X)$-bimodule $\mcal{D}^\lambda_{X \larr X_e}$ we have
\begin{align*}
 \mcal{D}^\lambda_{X \larr X_e}=i_e^{-1}\mcal{D}^\lambda_X \otimes_{i_e^{-1} \mcal{O}_X}\!\omega_{X_e/X},
\end{align*}
where $\omega_{X_e/X}=\omega_{X_e} \otimes_{i_e^{-1}\mcal{O}_X}\!i_e^{-1}\omega_X^{-1}$, $\omega_X$ and $\omega_{X_e}$ are the canonical sheaves of $X$ and $X_e$, respectively.
As a consequence of $\omega_{X_e/X} \simeq \mcal{O}_{X_e}$, we can write
\begin{align*}
  \Gamma(X,\mcal{I}(X_e,\mcal{O}_{X_e})) &= \Gamma(X,i_{e*}(\mcal{D}^\lambda_{X \larr X_e}\!\otimes_{i_e^\sharp\mcal{D}^\lambda_X}\! \mcal{O}_{X_e})) \\
  &= \Gamma(X_e,\mcal{D}^\lambda_{X \larr X_e}\!\otimes_{i_e^\sharp\mcal{D}^\lambda_X}\! \mcal{O}_{X_e}) \\
  &=\Gamma(X_e,(i_e^{-1}\mcal{D}^\lambda_X \otimes_{i_e^{-1} \mcal{O}_X}\!\omega_{X_e/X})\otimes_{i_e^\sharp\mcal{D}^\lambda_X}\! \mcal{O}_{X_e}) \\
  &\simeq \mcal{D}^\lambda_{X,o} \otimes_{\mcal{O}_{X,o}}\!\C ,
\end{align*}
where $o=eP$. Therefore, we need to know the geometric fiber of the right $\mcal{O}_X$-module $\mcal{D}^\lambda_X$ at the point $o$. The left action of $U(\mfrak{g})$ on $\mcal{D}^\lambda_{X,o} \otimes_{\mcal{O}_{X,o}}\!\C$ is given through the mapping
\begin{align*}
  U(\mfrak{g}) \rarr \Gamma(X,\mcal{D}^\lambda_X) \rarr \mcal{D}^\lambda_{X,o},
\end{align*}
which is a homomorphism of associative $\C$-algebras. From \eqref{eq:twisted sheaf realization} we have
\begin{align*}
  \mcal{D}^\lambda_X \simeq \mcal{U}_X(\mfrak{g})/{\textstyle \sum_{s \in \mcal{V}_X\!(\mfrak{p})}}(s-f_{\lambda+\rho}(s))\mcal{U}_X(\mfrak{g}),
\end{align*}
hence we get
\begin{align*}
  \C \otimes_{\mcal{O}_{X,o}}\! \mcal{D}^\lambda_{X,o} \simeq U(\mfrak{g}) / {\textstyle \sum_{A \in \mfrak{p}}} (A-(\lambda+\rho)(A))U(\mfrak{g}).
\end{align*}
Furthermore, if we use the isomorphism $\mcal{D}^\lambda_{X,o} \simeq (\mcal{D}^{-\lambda}_{X,o})^{\rm op}$ of associative $\C$-algebras, we obtain
\begin{align*}
  \mcal{D}^\lambda_{X,o} \otimes_{\mcal{O}_{X,o}}\!\C \simeq (\mcal{D}^{-\lambda}_{X,o})^{\rm op} \otimes_{\mcal{O}_{X,o}}\!\C \simeq \C \otimes_{\mcal{O}_{X,o}}\! \mcal{D}^{-\lambda}_{X,o} \simeq \C_{-\lambda+\rho} \otimes_{U(\mfrak{p})} U(\mfrak{g}).
\end{align*}
Hence $\C_{-\lambda+\rho} \otimes_{U(\mfrak{p})} U(\mfrak{g})$ has a natural right $U(\mfrak{g})$-module structure. The left $U(\mfrak{g})$-module structure on $\mcal{D}^\lambda_{X,o} \otimes_{\mcal{O}_{X,o}}\!\C$ is given as follows. Since the diagram
\begin{align*}
  \bfig
  \square<700,500>[U(\mfrak{g})`\mcal{D}^\lambda_{X,o}`U(\mfrak{g})^{\rm op}`(\mcal{D}^{-\lambda}_{X,o})^{\rm op};`\tau``]
  \efig
\end{align*}
is commutative, where $\tau \colon U(\mfrak{g}) \rarr U(\mfrak{g})^{\rm op}$ is the homomorphism of associative $\C$-algebras uniquely determined by $\tau(A)=-A$ for all $A \in \mfrak{g}$, we get
\begin{align*}
  a(1\otimes b)=1\otimes b\tau(a)
\end{align*}
for all $a,b \in U(\mfrak{g})$. It is easy to see that it is a highest weight module with highest weight $\lambda-\rho$. Moreover, it is also free $U(\widebar{\mfrak{u}})$-module of rank one with the free generator $1\otimes 1$, thus the generalized Verma module $M^\mfrak{g}_\mfrak{p}(\lambda-\rho)$.}

\lemma{\label{lem:quotient D module}The mapping
\begin{align}
  \tau_o \colon \mcal{D}_{X,o} \otimes_{\mcal{O}_{X,o}}\!\C \rarr \mcal{D}_{X,o}/\mcal{D}_{X,o}\mfrak{m}_o
\end{align}
defined by
\begin{align}
  \tau_o(Q\otimes v)=vQ\ {\rm mod}\ \mcal{D}_{X,o}\mfrak{m}_o
\end{align}
for $Q \in \mcal{D}_{X,o}$ and $v\in \C$, where $\mfrak{m}_o$ is the maximal ideal of $\mcal{O}_{X,o}$, is an isomorphism of left $\mcal{D}_{X,o}$-modules.}

\proof{For $Q \in \mcal{D}_{X,o}$ and $f\in \mcal{O}_{X,o}$, we have
\begin{align*}
  \tau_o(Qf\otimes 1) = Qf\ {\rm mod}\ \mcal{D}_{X,o}\mfrak{m}_o.
\end{align*}
On the other hand, since $Qf \otimes 1 = Q \otimes f(o)$, we obtain
\begin{align*}
  \tau_o(Q \otimes f(o)) =f(o)Q\ {\rm mod}\ \mcal{D}_{X,o}\mfrak{m}_o.
\end{align*}
But $f=f(o)+(f-f(o))$ implies $\tau_o(Qf\otimes 1)=\tau_o(Q \otimes f(o))$.
This means that the mapping $\tau_o$ is well-defined.

As for the injectivity of $\tau_o$, the condition $\tau_o(Q \otimes 1)=0$ implies $Q = \sum_{i=1}^m Q_i g_i$,
where $Q_i \in \mcal{D}_{X,o}$ and $g_i \in \mfrak{m}_o$. Therefore, we may write
\begin{align*}
Q \otimes 1 = \sum_{i=1}^m(Q_ig_i \otimes 1)= \sum_{i=1}^m(Q_i \otimes g_i(o))=0.
\end{align*}
The surjectivity of the mapping $\tau_o$ is obvious. The proof is complete.}

Let us consider the local chart $(U_e,u_e)$ on $X$ given by \eqref{eq:atlas flag manifold}. Since the mapping $u_e \colon U_e \rarr \widebar{\mfrak{u}}$ is a diffeomorphism, it induces the isomorphism
\begin{align}
  \Psi_{u_e} \colon \mcal{D}_X(U_e) \rarr \Gamma(\widebar{\mfrak{u}},\mcal{D}_{\widebar{\mfrak{u}}})
\end{align}
of associative $\C$-algebras. Let $(x_1,x_2,\dots,x_n)$ be linear coordinate functions on $\widebar{\mfrak{u}}$, then we get the coordinate functions $(u_e^1,u_e^2,\dots,u_e^n)$ on $U_e$ defined by
\begin{align}
  u_e^i = x_i \circ u_e
\end{align}
for $i=1,2,\dots,n$. Then we have
\begin{align}
  \Psi_{u_e}(u_e^i)=x_i \quad \text{and}\quad \Psi_{u_e}(\partial_{u_e^i})=\partial_{x_i}
\end{align}
for $i=1,2,\dots,n$.

We denote by $\eus{A}^\mfrak{g}_{\widebar{\mfrak{u}}}$ the Weyl algebra of the complex vector space $\widebar{\mfrak{u}}$ (see Section \ref{sec:Fourier transform} for the definition) and by $I_e$ the left ideal of $\eus{A}^\mfrak{g}_{\widebar{\mfrak{u}}}$ generated by polynomials on $\widebar{\mfrak{u}}$ vanishing at the point $0$. Let us note that the Weyl algebra $\eus{A}^\mfrak{g}_{\widebar{\mfrak{u}}}$ is contained in $\Gamma(\widebar{\mfrak{u}},\mcal{D}_{\widebar{\mfrak{u}}})$.
\medskip

\lemma{\label{lem:quotient A module}The mapping
\begin{align}
 \sigma_o \colon \eus{A}^\mfrak{g}_{\widebar{\mfrak{u}}}/I_e \rarr \mcal{D}_{X,o}/\mcal{D}_{X,o}\mfrak{m}_o
\end{align}
defined by
\begin{align}
  P\ {\rm mod}\ I_e \mapsto \Psi_{u_e}^{-1}(P)\ {\rm mod}\ \mcal{D}_{X,o}\mfrak{m}_o
\end{align}
is an isomorphism of left $\eus{A}^\mfrak{g}_{\widebar{\mfrak{u}}}$-modules.}

\proof{The canonical mapping $\eus{A}^\mfrak{g}_{\widebar{\mfrak{u}}} \rarr \mcal{D}_{X,o}$ induces the mapping $I_e \rarr \mcal{D}_{X,o}\mfrak{m}_o$, hence $\sigma_o$ is well-defined. Now, let us assume that $P \in \mcal{D}_{X,o}$, then we can write $P$ in the form
\begin{align*}
  P=\sum_{\gamma \in \N_0^n} (\partial_{u_e^1})^{\gamma_1}(\partial_{u_e^2})^{\gamma_2} \dots (\partial_{u_e^n})^{\gamma_n}a_\gamma,
\end{align*}
where $a_\gamma \in \mcal{O}_{X,o}$. Thus, if we define $Q \in \eus{A}^\mfrak{g}_{\widebar{\mfrak{u}}}$ by
\begin{align*}
  Q = \sum_{\gamma \in \N_0^n} (\partial_{x_1})^{\gamma_1}(\partial_{x_2})^{\gamma_2} \dots (\partial_{x_n})^{\gamma_n}a_\gamma(o),
\end{align*}
then we get $\sigma_o(Q\ {\rm mod}\ I_e) = P\ {\rm mod}\ \mcal{D}_{X,o}\mfrak{m}_o$. The injectivity of the mapping $\sigma_o$ is obvious.}

By Proposition \ref{lem:Verma module}, the generalized Verma module $M_\mfrak{p}^\mfrak{g}(\lambda-\rho)$ is isomorphic to $\smash{\mcal{D}^\lambda_{X,o} \otimes_{\mcal{O}_{X,o}}\!\C}$. If we use the isomorphism $j_{s_e} \colon \mcal{D}_X^\lambda|_{U_e} \rarr \mcal{D}_X|_{U_e}$ of sheaves of rings of twisted differential operators on $U_e$ given by \eqref{eq:twisted sheaf trivialization}, we obtain the $U(\mfrak{g})$-isomorphism
\begin{align}
  \varphi_\lambda \colon M_\mfrak{p}^\mfrak{g}(\lambda-\rho) \rarr \mcal{D}_{X,o} \otimes_{\mcal{O}_{X,o}}\!\C \label{eq:Verma module}
\end{align}
uniquely determined by
\begin{align}
  \varphi_\lambda(1\otimes 1)=1 \otimes 1. \label{eq:Verma module morphism identity}
\end{align}
The structure of a left $U(\mfrak{g})$-module on $\mcal{D}_{X,o} \otimes_{\mcal{O}_{X,o}}\!\C$
follows from the composition
\begin{align}
  U(\mfrak{g}) \rarr \Gamma(X,\mcal{D}^\lambda_X) \rarr \mcal{D}^\lambda_X(U_e) \rarr \mcal{D}_X(U_e) \rarr \mcal{D}_{X,o},
\end{align}
which is a homomorphism of associative $\C$-algebras. By Lemma \ref{lem:quotient D module} and  Lemma \ref{lem:quotient A module}, we obtain the isomorphism of $U(\mfrak{g})$-modules
\begin{align}
  \Phi_\lambda \colon M^\mfrak{g}_\mfrak{p}(\lambda-\rho) \rarr \eus{A}^\mfrak{g}_{\widebar{\mfrak{u}}}/I_e, \label{eq:Verma isomorphism}
\end{align}
uniquely determined by
\begin{align}
  \Phi_\lambda(1 \otimes 1)= 1\ {\rm mod}\ I_e.
\end{align}
The structure of a left $U(\mfrak{g})$-module on $\eus{A}^\mfrak{g}_{\widebar{\mfrak{u}}}$ is induced by
the homomorphism of associative $\C$-algebras
\begin{align}
  \pi_\lambda \colon U(\mfrak{g}) \rarr \Gamma(X,\mcal{D}^\lambda_X) \rarr \mcal{D}^\lambda_X(U_e) \rarr \mcal{D}_X(U_e) \rarr \Gamma(\widebar{\mfrak{u}},\mcal{D}_{\widebar{\mfrak{u}}}), \label{eq:algebra action}
\end{align}
where $\pi_\lambda(U(\mfrak{g})) \subset \eus{A}^\mfrak{g}_{\widebar{\mfrak{u}}}$
as it follows from Theorem \ref{prop:operator realization}.

\medskip

It the next Theorem we introduce coordinate-free description of a
linear surjective mapping from the symmetric algebra $S(\widebar{\mfrak{u}})$
to the universal enveloping algebra $U(\widebar{\mfrak{u}})$ of the Lie
algebra $\widebar{\mfrak{u}}$, completely characterized to be the
identity map on $\widebar{\mfrak{u}}$.
Hereby we realize $\Phi_\lambda$ in an explicit way. This result will
be used in the subsequent Section \ref{sec:Heisenberg parabolic} for
the construction of both singular vectors and equivariant differential
operators as elements of $U(\widebar{\mfrak{u}})$, where
$\widebar{\mfrak{u}}$ the opposite nilradical of the parabolic subalgebra,
see \eqref{parnilradical}.
\medskip

\theorem{\label{thm:symmetrization}Let $(f_1,f_2,\dots,f_n)$ be a basis
of $\widebar{\mfrak{u}}$, $(x_1,x_2,\dots,x_n)$ be the corresponding
linear coordinate functions on $\widebar{\mfrak{u}}$, and let
$\beta \colon S(\widebar{\mfrak{u}}) \rarr U(\widebar{\mfrak{u}})$ be the
symmetrization map defined by
\begin{align}
  \beta(f_1f_2\dots f_k) = {1 \over k!} \sum_{\sigma \in S_k} f_{\sigma(1)}f_{\sigma(2)} \dots f_{\sigma(k)}
\end{align}
for all $k \in \N$ and $f_1,f_2,\dots,f_k \in \widebar{\mfrak{u}}$.
Then
\begin{align}
  \Phi_\lambda(\beta(f_{i_1}f_{i_2}\dots f_{i_k}) \otimes 1) = (-1)^k \partial_{i_1}\partial_{i_2} \dots \partial_{i_k}\ {\rm mod}\ I_e, \label{eq:inverse mapping}
\end{align}
where $\partial_i=\partial_{x_i}$, for all $k \in \N$ and $i_1,i_2,\dots,i_k \in \{1,2,\dots,n\}$.}

\proof{For $f \in \widebar{\mfrak{u}}$, we get from Theorem \ref{prop:operator realization}
\begin{align*}
  \pi_\lambda(f) = - \sum_{i=1}^n \sum_{k=0}^\infty g^i_k(f) \partial_{x_i},
\end{align*}
where
\begin{align*}
  g^i_k(f) = \alpha_k \sum_{m_1,\dots,m_k} x_{m_1}\dots x_{m_k} [(\ad(f_{m_1}) \dots \ad(f_{m_k}))(f)]_i
\end{align*}
for the Bernoulli numbers $\alpha_k$, determined by the generating series
\begin{align*}
  {z \over e^z-1} = \sum_{k=0}^\infty \alpha_k z^k
\end{align*}
for $z \in \C$ satisfying $0< |z| < 1$. We shall prove the claim by induction. Because $g_0^i(f_{i_1})=\delta^i_{i_1}$, we can write
\begin{align*}
  \Phi_\lambda(\beta(f_{i_1})\otimes 1) &= \pi_\lambda(f_{i_1})\ {\rm mod}\ I_e = - \sum_{i=1}^n (g_0^i(f_{i_1})\partial_i+ g_1^i(f_{i_1})\partial_i)\ {\rm mod}\ I_e \\
  &= -\partial_{i_1} - \alpha_1 \sum_{i=1}^n [\ad(f_{i_1})(f_i)]_i\ {\rm mod}\ I_e \\
  & = -\partial_{i_1} - \alpha_1 \tr_{\widebar{\mfrak{u}}}(\ad(f_{i_1}))\ {\rm mod}\ I_e = -\partial_{i_1}\ {\rm mod}\ I_e.
\end{align*}
Let us now assume that \eqref{eq:inverse mapping} holds for some $k\in \N$. The formula
\begin{align*}
  \beta(f_{i_1}\dots f_{i_{k+1}}) = {1 \over k+1} \sum_{j=1}^{k+1} f_{i_j} \beta(f_{i_1} \dots \widehat{f}_{i_j}\dots f_{i_{k+1}})
\end{align*}
allow us to write
\begin{align*}
(k+1)\Phi_\lambda(\beta(f_{i_1}f_{i_2}\dots f_{i_{k+1}}) \otimes 1)&= {\textstyle \sum_{j=1}^{k+1}} \pi_\lambda(f_{i_j}) \Phi_\lambda(\beta(f_{i_1}\dots \widehat{f}_{i_j} \dots f_{i_{k+1}})\otimes 1)\ {\rm mod}\ I_e \\ &= (-1)^k {\textstyle \sum_{j=1}^{k+1}} \pi_\lambda(f_{i_j})\, \partial_{i_1} \dots \widehat{\partial}_{i_j} \dots \partial_{i_{k+1}}\ {\rm mod}\ I_e \\
&= (-1)^k {\textstyle \sum_{j=1}^{k+1}} \partial_{i_1} \dots \widehat{\partial}_{i_j} \dots \partial_{i_{k+1}}\pi_\lambda(f_{i_j}) \\
& \quad + (-1)^k {\textstyle \sum_{j=1}^{k+1}} [\pi_\lambda(f_{i_j}),\partial_{i_1} \dots \widehat{\partial}_{i_j} \dots \partial_{i_{k+1}}]\ {\rm mod}\ I_e \\
&= (-1)^{k+1}(k+1)\,\partial_{i_1} \dots \partial_{i_{k+1}} \\
& \quad + (-1)^k {\textstyle \sum_{j=1}^{k+1}} [\pi_\lambda(f_{i_j}),\partial_{i_1} \dots \widehat{\partial}_{i_j} \dots \partial_{i_{k+1}}]\ {\rm mod}\ I_e,
\end{align*}
so that we get
\begin{multline*}
  \Phi_\lambda(\beta(f_{i_1}f_{i_2}\dots f_{i_{k+1}}) \otimes 1) = (-1)^{k+1}\partial_{i_1} \dots \partial_{i_{k+1}} \\ + {(-1)^k \over (k+1)!} \sum_{\sigma \in S_{k+1}}  [\pi_\lambda(f_{i_{\sigma(1)}}),\partial_{i_{\sigma(2)}} \dots \partial_{i_{\sigma(k+1)}}]\ {\rm mod}\ I_e.
\end{multline*}
Therefore, it is enough to show that
\begin{align*}
   \sum_{\sigma \in S_r}[\pi_\lambda(f_{j_{\sigma(1)}}),\partial_{j_{\sigma(2)}}]\partial_{j_{\sigma(3)}} \dots \partial_{j_{\sigma(r)}}\ {\rm mod}\ I_e = 0
\end{align*}
for all $r \geq 2$. Because we have
\begin{align*}
  \sum_{\sigma \in S_r}[\pi_\lambda(f_{j_{\sigma(1)}}),\partial_{j_{\sigma(2)}}]\partial_{j_{\sigma(3)}} \dots \partial_{j_{\sigma(r)}}= \sum_{k=0}^\infty \sum_{i=1}^n \sum_{\sigma \in S_r}\,  [\partial_{j_{\sigma(2)}},g^i_k(f_{j_{\sigma(1)}})]\partial_{j_{\sigma(3)}} \dots \partial_{j_{\sigma(r)}}  \partial_i
  \end{align*}
for $r \geq 2$, it is clearly sufficient to prove
\begin{align}
  \sum_{k=0}^\infty \sum_{i=1}^n \sum_{\sigma \in S_r}  [\partial_{j_{\sigma(2)}},g^i_k(f_{j_{\sigma(1)}})]\,\partial_{j_{\sigma(3)}} \dots \partial_{j_{\sigma(r)}}  \partial_i\ {\rm mod}\ I_e = 0
	\label{eq:commutator}
\end{align}
for all $r \geq 2$. Again, we shall prove this fact by induction on $r$. If $r=2$, we can rewrite the left hand side of \eqref{eq:commutator} into the form
\begin{multline*}
  \sum_{k=0}^\infty \sum_{i=1}^n \sum_{\sigma \in S_2}\!\big(\partial_i  [\partial_{j_{\sigma(2)}},g^i_k(f_{j_{\sigma(1)}})] - [\partial_i,  [\partial_{j_{\sigma(2)}},g^i_k(f_{j_{\sigma(1)}})]]\big)\ {\rm mod}\ I_e \\ = \sum_{i=1}^n \sum_{\sigma \in S_2}\!\big(\partial_i  [\partial_{j_{\sigma(2)}},g^i_1(f_{j_{\sigma(1)}})] - [\partial_i,  [\partial_{j_{\sigma(2)}},g^i_2(f_{j_{\sigma(1)}})]]\big)\ {\rm mod}\ I_e,
\end{multline*}
where we used the fact that $g^i_k(f_j)$ is a polynomial of degree $k$ in the variables $\{x_1,x_2,\dots,x_n\}$ and $g^i_0(f_j)=\delta^i_j$. In addition, we have
\begin{align*}
  \sum_{i=1}^n \sum_{\sigma \in S_2} \partial_i[\partial_{j_{\sigma(2)}},g^i_1(f_{j_{\sigma(1)}})]= \alpha_1 \sum_{i=1}^n \sum_{\sigma \in S_2} \partial_i [\ad(f_{j_{\sigma(2)}})(f_{j_{\sigma(1)}})]_i=0
\end{align*}
and
\begin{align*}
  \sum_{i=1}^n \sum_{\sigma \in S_2}[\partial_i,  [\partial_{j_{\sigma(2)}},g^i_2(f_{j_{\sigma(1)}})]]&= \alpha_2\sum_{i=1}^n \sum_{\sigma \in S_2}\! \big([\ad(f_i)\ad(f_{j_{\sigma(2)}})(f_{j_{\sigma(1)}})]_i-[\ad(f_{j_{\sigma(2)}}) \ad(f_{j_{\sigma(1)}})(f_i)]_i\big) \\
  &= -\alpha_2\sum_{i=1}^n \sum_{\sigma \in S_2} [\ad(f_{j_{\sigma(2)}}) \ad(f_{j_{\sigma(1)}})(f_i)]_i \\ &= -\alpha_2 \sum_{\sigma \in S_2} \tr_{\widebar{\mfrak{u}}}(\ad(f_{j_{\sigma(2)}}) \ad(f_{j_{\sigma(1)}})) = 0,
\end{align*}
since $\ad(f_{j_{\sigma(2)}}) \ad(f_{j_{\sigma(1)}})$ is a nilpotent mapping on $\widebar{\mfrak{u}}$. Now, let us assume that \eqref{eq:commutator} holds for all $r \in \{2,3,\dots,r_0-1\}$, then the left hand side \eqref{eq:commutator} reduces to
\begin{multline*}
   (-1)^{r_0-1} \sum_{k=0}^\infty \sum_{i=1}^n \sum_{\sigma \in S_{r_0}} [\partial_i,[\partial_{j_{\sigma(r_0)}},[\dots , [\partial_{j_{\sigma(2)}},g^i_k(f_{j_{\sigma(1)}})]\dots]]]\ {\rm mod}\ I_e \\ = (-1)^{r_0-1} \sum_{i=1}^n \sum_{\sigma \in S_{r_0}} [\partial_i,[\partial_{j_{\sigma(r_0)}},[\dots , [\partial_{j_{\sigma(2)}},g^i_{r_0}(f_{j_{\sigma(1)}})]\dots]]]\ {\rm mod}\ I_e \\ = \alpha_{r_0}(-1)^{r_0}(r_0-1)! \sum_{i=1}^n\sum_{\sigma \in S_{r_0}}[\ad(f_{j_{\sigma(r_0)}})\dots \ad(f_{j_{\sigma(2)}})\ad(f_{j_{\sigma(1)}})(f_i)]_i\ {\rm mod}\ I_e \\ = \alpha_{r_0}(-1)^{r_0}(r_0-1)!\sum_{\sigma \in S_{r_0}}\! \tr_{\widebar{\mfrak{u}}}(\ad(f_{j_{\sigma(r_0)}})\dots \ad(f_{j_{\sigma(1)}}))\ {\rm mod}\ I_e =0.
\end{multline*}
The proof is complete.}


\subsection{Algebraic Fourier transform}
\label{sec:Fourier transform}

In this section we recall the Fourier transform of modules over Weyl algebra
and use it later to convert the algebraic characterization of certain vectors in
modules in the BGG parabolic category $\mcal{O}^\mfrak{p}$ into characterization by a system of
linear partial differential equations.


Let $V$ be a finite-dimensional complex vector space, regarded as a complex algebraic
variety $(V^{\rm alg},\mcal{O}_{V^{\rm alg}})$. The associative
$\C$-algebra $\eus{A}_V = \Gamma(V^{\rm alg},\mcal{D}_{V^{\rm alg}})$ is called the
Weyl algebra. Since the canonical morphism $i \colon V \rarr V^{\rm alg}$ of topological spaces
induces an injective morphism $\eus{A}_V \rarr \Gamma(V,\mcal{D}_V)$ of associative
$\C$-algebras, the Weyl algebra $\eus{A}_V$ can be regarded as a subalgebra of \
$\Gamma(V,\mcal{D}_V)$ (cf.\ Section \ref{sec:Twisted differential}).
We have
\begin{align}
  \eus{A}_V \simeq S(V^*)\otimes_\C S(V),
\end{align}
where $S(V^*)\simeq \C[V] = \Gamma(V^{\rm alg},\mcal{O}_{V^{\rm alg}})$ and we regard $S(V) \simeq \C[V^*]$ as the $\C$-algebra of constant coefficient differential operators.

Let $(x_1,x_2,\dots,x_n)$ be linear coordinate functions on $V$ and let $(y_1,y_2,\dots,y_n)$ be the dual linear coordinate functions on $V^*$. Then there is a canonical isomorphism
\begin{align}\label{eq:aft}
  \mcal{F} \colon \eus{A}_V \rarr \eus{A}_{V^*}
\end{align}
of associative $\C$-algebras given by
\begin{align}
  \mcal{F}(x_i) =-\partial_{y_i}, \qquad \mcal{F}(\partial_{x_i})= y_i \label{eq:Fourier transform}
\end{align}
for $i=1,2,\dots,n$. Let us note that the definition does not depend on the choice of linear coordinates on $V$.

Let $M$ be a left $\eus{A}_V$-module. The Fourier transform $\widehat{M}$ of $M$ has the same underlying vector space as $M$, and the left $\eus{A}_{V^*}$-module structure is given by
\begin{align}
  Pu=\mcal{F}^{-1}(P)u
\end{align}
for all $u \in \widehat{M}$ and $P\in \eus{A}_{V^*}$. The Fourier transform induces equivalences of categories
\begin{align}
  \widehat{\phantom{M}} \colon \Mod(\eus{A}_V) \rarr \Mod(\eus{A}_{V^*}\!), \qquad
  \widehat{\phantom{M}} \colon \Mod_f(\eus{A}_V) \rarr \Mod_f(\eus{A}_{V^*}\!),
\end{align}
where $\Mod_f(\eus{A}_V)$ is the category of finitely generated $\eus{A}_V$-modules.
\medskip

\lemma{Let $I$ be a left ideal of $\eus{A}_V$ and let $M$ be a left $\eus{A}_V$-module of the form
\begin{align}
  M=\eus{A}_V/I.
\end{align}
Then the Fourier transform $\smash{\widehat{M}}$ is isomorphic to the left $\eus{A}_{V^*}$-module $\smash{\widetilde{M}}$ of the form
\begin{align}
  \widetilde{M}=\eus{A}_{V^*}/\mcal{F}(I),
\end{align}
where the isomorphism
\begin{align}
  \varphi \colon \smash{\widetilde{M}} \rarr \smash{\widehat{M}}
\end{align}
is given by
\begin{align}
  Q\ {\rm mod}\ \mcal{F}(I) \mapsto \mcal{F}^{-1}(Q)\ {\rm mod}\ I
\end{align}
for all $Q \in \eus{A}_V$.}

\proof{The morphism $\psi \colon \eus{A}_{V^*} \rarr \smash{\widehat{M}}$
of left $\eus{A}_{V^*}$-modules defined by
\begin{align*}
  \psi(Q) = \mcal{F}^{-1}(Q)\ {\rm mod}\ I
\end{align*}
for all $Q \in \eus{A}_{V^*}$ is clearly surjective. It is also injective, because $\psi(Q)=0$ implies $Q \in \mcal{F}(I)$, since $\mcal{F}$ is an isomorphism of associative $\C$-algebras, and therefore $\ker \psi \subset \mcal{F}(I)$. The opposite inclusion is trivial, hence $\psi$ induces the isomorphism $\varphi \colon \smash{\widetilde{M}} \rarr \smash{\widehat{M}}$.}

For a left $\eus{A}_V$-module $M$ we consider a system of linear partial differential equations for a single unknown element $u \in M$ in the form
\begin{align}
  P_1u=P_2u=\dots =P_ku=0, \label{eq:system equation}
\end{align}
where $k\in \N$ and $P_1,P_2,\dots,P_k \in \eus{A}_V$. We denote by $\Sol(P_1,P_2,\dots,P_k;M)$ the complex vector space of solutions of the system \eqref{eq:system equation}.
\medskip

\lemma{Let $I$ be a left ideal of $\eus{A}_V$ and let $M$ be a left $\eus{A}_V$-module of the form $M=\eus{A}_V/I$. Let us consider elements $P_1,P_2,\dots,P_k \in \eus{A}_V$ for $k\in \N$.
Then the mapping
\begin{align}
 \tau \colon  \Sol(P_1,P_2,\dots,P_k;M)  \rarr  \Sol(\mcal{F}(P_1),\mcal{F}(P_2),\dots,\mcal{F}(P_k);\smash{\widetilde{M}}) \label{eq:Fourier transform module}
\end{align}
given by
\begin{align}
  Q\ {\rm mod}\ I \mapsto \mcal{F}(Q)\ {\rm mod}\ \mcal{F}(I)
\end{align}
is an isomorphism of complex vector spaces.}

\proof{First of all we show that $\tau$ maps into $\Sol(\mcal{F}(P_1),\mcal{F}(P_2),\dots,\mcal{F}(P_k);\smash{\widetilde{M}})$. For an element $Q\ {\rm mod}\ I \in \Sol(P_1,P_2,\dots,P_k;M)$, we get
\begin{align*}
  \mcal{F}(P_i)\tau(Q\ {\rm mod}\ I)=\mcal{F}(P_i)\mcal{F}(Q)\  {\rm mod}\ \mcal{F}(I) = \mcal{F}(P_iQ)\  {\rm mod}\ \mcal{F}(I) = 0\  {\rm mod}\ \mcal{F}(I)
\end{align*}
for $i=1,2,\dots,k$, where we used the fact that $P_iQ \in I$ and so $\mcal{F}(P_iQ) \in \mcal{F}(I)$ for $i=1,2,\dots,k$. Then it is easy to see that the mapping
\begin{align*}
  \tilde{\tau} \colon \Sol(\mcal{F}(P_1),\mcal{F}(P_2),\dots,\mcal{F}(P_k);\widetilde{M}) \rarr \Sol(P_1,P_2,\dots,P_k;M)
\end{align*}
defined by the formula
\begin{align*}
  Q\ {\rm mod}\ \mcal{F}(I) \mapsto \mcal{F}^{-1}(Q)\ {\rm mod}\ I.
\end{align*}
is the inverse mapping to $\tau$, and so $\tau$ is an isomorphism. The proof is complete.}


\subsection{Singular vectors in generalized Verma modules}

This section contains the definition of the key concept of singular vectors. The
complete collection of singular vectors in generalized Verma modules in the BGG parabolic category $\mcal{O}^\mfrak{p}$, or
in particular, in a given block $\mcal{O}_{\chi_\lambda}^\mfrak{p}$, fully encodes
the structure of its morphisms, \cite{Bernstein-Gelfand-Gelfand1971, Humphreys-book}.

Let us consider a complex reductive Lie subalgebra $\mfrak{g}'$ of a complex reductive Lie algebra
$\mfrak{g}$ together with its parabolic subalgebra $\mfrak{p}\subset\mfrak{g}$
such that $\mfrak{p}' = \mfrak{p} \cap \mfrak{g}'$ is a parabolic subalgebra of $\mfrak{g}'$.
We recall that the formula \eqref{eq:Verma isomorphism} associates to any element
$\lambda\in \Hom_P(\mfrak{p},\C)$ the mapping $\Phi_\lambda$, which
induces the isomorphism of $\mfrak{l}'$-modules
\begin{align}
  \smash{M^\mfrak{g}_\mfrak{p}(\lambda-\rho)^{\mfrak{u}'}} \riso
	\Sol(\mfrak{g},\mfrak{g}',\mfrak{p};\eus{A}^\mfrak{g}_{\widebar{\mfrak{u}}}/I_e). \label{eq:sv}
\end{align}
Here $\smash{M^\mfrak{g}_\mfrak{p}(\lambda-\rho)^{\mfrak{u}'}}$ denotes the subspace in the
scalar generalized Verma module $M^\mfrak{g}_\mfrak{p}(\lambda-\rho)$ (cf. (\ref{eq:generalized Verma module}))
annihilated by the nilradical $\mfrak{u}'$ of $\mfrak{p}'$, and the elements in
$\smash{M^\mfrak{g}_\mfrak{p}(\lambda-\rho)^{\mfrak{u}'}}$ are called
$\mfrak{g}'$-singular vectors. The right hand side of \eqref{eq:sv}
is defined by
\begin{align}
  \Sol(\mfrak{g},\mfrak{g}',\mfrak{p};\eus{A}^\mfrak{g}_{\widebar{\mfrak{u}}}/I_e) = \Sol(\{\pi_\lambda(X);\, X \in \mfrak{u}'\}; \eus{A}^\mfrak{g}_{\widebar{\mfrak{u}}}/I_e).
\end{align}
Since the Levi factor $\mfrak{l}'$ of $\mfrak{p}'$ is a reductive Lie algebra, we have
$\mfrak{l}'=\mfrak{z}(\mfrak{l}') \oplus [\mfrak{l}',\mfrak{l}']$ with $\mfrak{z}(\mfrak{l}')$
the center of $\mfrak{l}'$, and $\mfrak{z}(\mfrak{l}')$ acts on any finite-dimensional
irreducible $\mfrak{l}'$-module $V$ by a character
$\mu \colon \mfrak{z}(\mfrak{l}') \rarr \C$ fulfilling $Xv=\mu(X)v$ for all
$X \in \mfrak{z}(\mfrak{l}')$, $v \in V$.

Let $V\subset\smash{M_\mfrak{p}^\mfrak{g}(\lambda-\rho)^{\mfrak{u}'}}$ be a finite-dimensional
irreducible $\mfrak{l}'$-module, where the action of $\mfrak{z}(\mfrak{l}')$ on $V$ is given
by a character $\mu \colon \mfrak{z}(\mfrak{l}') \rarr \C$. Then we have
\begin{align}
  \Phi_\lambda(V) \subset \Sol(\{\pi_\lambda(X), \pi_\lambda(Y)-\mu(Y);\, X \in \mfrak{u}',Y \in \mfrak{z}(\mfrak{l}')\}; \eus{A}^\mfrak{g}_{\widebar{\mfrak{u}}}/I_e).
\end{align}
The application of Fourier transform, see
\eqref{eq:Fourier transform module}, implies
\begin{align}
  \tau(\Phi_\lambda(V)) \subset \Sol(\{\hat{\pi}_\lambda(X),\hat{\pi}_\lambda(Y)-\mu(Y);\, X \in \mfrak{u}', Y \in \mfrak{z}(\mfrak{l}')\}; \eus{A}^\mfrak{g}_{\widebar{\mfrak{u}}^*}\!/\mcal{F}(I_e)),
\end{align}
where
\begin{align}
\Sol(\mfrak{g},\mfrak{g}',\mfrak{p};\eus{A}^\mfrak{g}_{\widebar{\mfrak{u}}^*}\!/\mcal{F}(I_e))_\mu^\mcal{F}=
\Sol(\{\hat{\pi}_\lambda(X), \hat{\pi}_\lambda(Y)-\mu(Y);\, X \in \mfrak{u}',Y \in \mfrak{z}(\mfrak{l}')\}; \eus{A}^\mfrak{g}_{\widebar{\mfrak{u}}^*}\!/\mcal{F}(I_e))
\end{align}
is an ${\mathfrak l}'$-submodule in the polynomial ring
$\C[\widebar{\mfrak{u}}^*] \simeq \eus{A}^\mfrak{g}_{\widebar{\mfrak{u}}^*}\!/\mcal{F}(I_e)$
given by the solution of a system of partial differential equations.

A ${\mfrak g}'$-singular vector is a generator of a $\mfrak{g}'$-submodule in
$M^\mfrak{g}_\mfrak{p}(\lambda-\rho)$, and its existence is equivalent to a discrete
component in the branching problem $M^\mfrak{g}_\mfrak{p}(\lambda)|_{\mfrak{g}'}$.
The results in our article classify discrete components in the branching problem
for a class of generalized Verma modules regarded as objects in the BGG
parabolic category
$\smash{\mcal{O}^{\mfrak{p}}}$ and a class of $\mfrak{g}'$-compatible
parabolic subalgebras $\mfrak{p}$. The condition of $\mfrak{g}'$-compatibility for
parabolic subalgebras guarantees that the decompositions are discrete and thus
realized by singular vectors, \cite{Kobayashi1994, Kobayashi2012}.
Let us recall that a parabolic subalgebra
$\mfrak{p}\subset\mfrak{g}$ is $\mfrak{g}'$-compatible provided there exists a hyperbolic element
$E'\in\mfrak{g}'$ such that $\mfrak{p}$ is the direct sum of eigenspaces of $\ad(E')$ with non-negative integral
eigenvalues. This then implies that $G'\!P$ is a closed submanifold of $G$
and a geometric argument inspired by $\mcal{D}$-module theory gives the discrete decomposability
of all modules in $\mcal{O}^\mfrak{p}$ with respect to $\mfrak{g}'$.

Assuming that $\mfrak{p}=\mfrak{l} \oplus \mfrak{u}$ is a $\mfrak{g}'$-compatible
parabolic subalgebra of
$\mfrak{g}$ and $V_\lambda$ a finite-dimensional irreducible $\mfrak{l}$-module
with highest weight $\lambda$, we consider the induced $\mfrak{l}'$-module structure on the
symmetric algebra
$S(\widebar{\mfrak{u}}/(\widebar{\mfrak{u}}\cap\mfrak{g}'))$ and define for any
finite-dimensional irreducible $\mfrak{l}'$-module $V_\mu'$ with highest weight
$\mu$ the multiplicity function
\begin{align}
m(\lambda,\mu)=\dim_{\mathbb C} \Hom_{\mfrak{l}'}(V'_\mu, V_\lambda|_{\mfrak{l}'}
\otimes S(\widebar{\mfrak{u}}/(\widebar{\mfrak{u}}\cap\mfrak{g}'))).  \label{eq:multiplicity formula}
\end{align}
Then $m(\lambda,\mu)<\infty$ for all $\mu\in \Lambda^+(\mfrak{p}')$, and in the Grothendieck group
$\smash{K(\mcal{O}^{{\mfrak p}})}$ of the BGG parabolic category $\smash{\mcal{O}^{{\mfrak p}}}$ holds
\begin{align}
M^\mfrak{g}_\mfrak{p}(\lambda)|_{\mfrak{g}'}\simeq \bigoplus_{\mu \in \Lambda^+(\mfrak{p}')} m(\lambda,\mu)\,
M^{\mfrak{g}'}_\mfrak{{p}'}\!(\mu) \label{eq:reduction}
\end{align}
for a generalized Verma module $M^\mfrak{g}_\mfrak{p}(\lambda)$. We shall
use \eqref{eq:reduction} in Section \ref{sec:qbp} below to produce in a particular case
of our interest a characterization of the generators of
${\mfrak g}'$-submodules.


\section{$A_n$-series of Lie algebras with Heisenberg parabolic subalgebras}
\label{sec:Heisenberg parabolic}

In the present section we implement and apply the previous general exposition
to the case of scalar generalized Verma modules for the pair given by the classical
complex Lie algebra $A_{n+1}$ and its parabolic subalgebra with the nilradical
isomorphic to the Heisenberg Lie algebra, and its complex Lie subalgebra $A_{n-r+1}$
($n-r \geq 1$) together with compatible parabolic subalgebra whose nilradical is
isomorphic to the Heisenberg Lie algebra of less dimension. In the Dynkin diagrammatic
notation, this type of parabolic subalgebra is determined by omitting the first and the
last simple nodes in the diagram.

Let us briefly indicate the road map for the present section. Following
Section \ref{sec:Flag manifolds}, we set our representation theoretical conventions
and use Theorem \ref{prop:operator realization} to describe the embedding of
the Lie algebra $\mfrak{sl}(n+2,\C)$ into the Weyl algebra
$\eus{A}^\mfrak{g}_{\widebar{\mfrak{u}}}$. As a $\mfrak{g}$-module we take
the left $\eus{A}^\mfrak{g}_{\widebar{\mfrak{u}}}$-module given as the quotient
of $\eus{A}^\mfrak{g}_{\widebar{\mfrak{u}}}$ by the left ideal of
$\eus{A}^\mfrak{g}_{\widebar{\mfrak{u}}}$ generated by polynomials
vanishing at the origin of $\widebar{\mfrak{u}}$, and apply to it the algebraic Fourier transform
\eqref{eq:aft}. This converts the former task of finding
$\mfrak{sl}(n-r+2,\C)$-singular vectors into the problem of solving a system
of partial differential equations, realized in
$\eus{A}^\mfrak{g}_{\widebar{\mfrak{u}}^*}$ and acting on the polynomial algebra
$\C[\widebar{\mfrak{u}}^*]$. The main technical device
allowing the complete classification of the solution space of this system of
partial differential equations is the Fischer decomposition, describing the
$(\mfrak{sl}(n-r,\C) \oplus \mfrak{sl}(r,\C))$-module structure on
$\C[(\mfrak{g}_{-1})^*]$, see Appendix \ref{app:Fischer decompostion}.

\subsection{Representation theoretical conventions}
\label{sub:rtc}

In the rest of the article we consider the complex semisimple Lie group $G=\SL(n+2,\C)$, $n\in \N$, and its Lie algebra
$\mfrak{g}=\mfrak{sl}(n+2,\C)$. The Cartan subalgebra $\mfrak{h}$ of $\mfrak{g}$ is given by diagonal
matrices
\begin{align}
  \mfrak{h}=\{\diag(a_1,a_2,\dots,a_{n+2});\, a_1,a_2,\dots,a_{n+2} \in \C,\ {\textstyle \sum_{i=1}^{n+2}}a_i=0\}.
\end{align}
For $i=1,2,\dots,n+2$ we define $\veps_i \in \mfrak{h}^*$
by $\veps_i(\diag(a_1,a_2,\dots,a_{n+2}))=a_i$. Then the root system of $\mfrak{g}$ with respect to $\mfrak{h}$ is
$\Delta = \{\veps_i-\veps_j;\, 1 \leq i \neq j \leq n+2\}$. The root space
 $\mfrak{g}_{\veps_i-\veps_j}$ is the complex linear span of $e_{ij}$,
the $(n+2 \times n+2)$-matrix such that $(e_{ij})_{k\ell}=\delta_{ik}\delta_{j\ell}$.
The positive root system is $\Delta^+=\{\veps_i-\veps_j;\, 1\leq i<j \leq n+2\}$, in which the set of simple roots is
$\Pi=\{\alpha_1,\alpha_2,\dots,\alpha_{n+1}\}$, $\alpha_i=\veps_i-\veps_{i+1},\, i=1,2,\dots,n+1$.
Then the fundamental weights are
$\omega_i= \smash{\sum_{j=1}^i \veps_j}$, $i=1,2,\dots,n+1$.
The Lie subalgebras $\mfrak{b}$ and $\widebar{\mfrak{b}}$ defined by the linear span of positive and negative root
spaces together with the Cartan subalgebra are called the standard Borel subalgebra and the opposite standard Borel subalgebra of $\mfrak{g}$, respectively.
The subset
$\Sigma=\{\alpha_2,\alpha_3,\dots,\alpha_n\}$
of $\Pi$ generates the root subsystem $\Delta_\Sigma$ in $\mfrak{h}^*$, and we associate to $\Sigma$ the standard parabolic
subalgebra $\mfrak{p}$ of $\mfrak{g}$ by
  $\mfrak{p} = \mfrak{l} \oplus \mfrak{u}$.
The reductive Levi factor $\mfrak{l}$ of $\mfrak{p}$ is defined through
\begin{align}
  \mfrak{l}= \mfrak{h} \oplus \bigoplus_{\alpha \in \Delta_\Sigma} \mfrak{g}_\alpha,
\end{align}
and the nilradical $\mfrak{u}$ of $\mfrak{p}$ and the opposite nilradical $\widebar{\mfrak{u}}$ are
\begin{align}
  \mfrak{u}= \bigoplus_{\alpha \in \Delta^+ \rsetminus \Delta_\Sigma^+}\mfrak{g}_\alpha \qquad \text{and}
  \qquad \widebar{\mfrak{u}}= \bigoplus_{\alpha \in \Delta^+ \rsetminus \Delta_\Sigma^+} \mfrak{g}_{-\alpha},
\end{align}
respectively. We define the $\Sigma$-height $\htt_\Sigma(\alpha)$ of $\alpha \in \Delta$ by
\begin{align}
  \htt_\Sigma\!\big({\textstyle \sum_{i=1}^{n+1}} a_i \alpha_i\big) = a_1 + a_{n+1},
\end{align}
so $\mfrak{g}$ is a $|2|$-graded Lie algebra with respect to the grading given by
$\mfrak{g}_i = \bigoplus_{\alpha \in \Delta,\, \htt_\Sigma(\alpha)=i} \mfrak{g}_\alpha$
for $0 \neq i \in \Z$, and
$\mfrak{g}_0= \mfrak{h} \oplus \smash{\bigoplus}_{\alpha \in \Delta,\, \htt_\Sigma(\alpha)=0} \mfrak{g}_\alpha$. Moreover, we have
$\mfrak{u}=\mfrak{g}_1 \oplus \mfrak{g}_2$,
$\widebar{\mfrak{u}}=\mfrak{g}_{-2} \oplus \mfrak{g}_{-1}$ and $\mfrak{l}= \mfrak{g}_0$.

The basis $(f_1,\dots,f_n,g_1,\dots,g_n,c)$ of the root spaces in the opposite nilradical $\widebar{\mfrak{u}}$ is given by
\begin{align}
  f_i =\begin{pmatrix}
    0 & 0 & 0 \\
    1_i & 0 & 0 \\
    0 & 0 & 0
  \end{pmatrix}\!,\quad
  g_i =\begin{pmatrix}
    0 & 0 & 0 \\
    0 & 0 & 0 \\
    0 & 1_i^{\rm T} & 0
  \end{pmatrix}\!,\quad
  c =\begin{pmatrix}
    0 & 0 & 0 \\
    0 & 0 & 0 \\
    1 & 0 & 0
  \end{pmatrix}\!,
\end{align}
where the only non-trivial Lie brackets are $[f_i,g_i]=-c$
for all $i=1,2,\dots,n$. Analogously, the basis $(d_1,\dots,d_n,e_1,\dots,e_n,a)$ of the root spaces in $\mfrak{u}$ is given by
\begin{align}
  d_i =\begin{pmatrix}
    0 & 1_i^{\rm T} & 0 \\
    0 & 0 & 0 \\
    0 & 0 & 0
  \end{pmatrix}\!,\quad
  e_i =\begin{pmatrix}
    0 & 0 & 0 \\
    0 & 0 & 1_i \\
    0 & 0 & 0
  \end{pmatrix}\!,\quad
  a =\begin{pmatrix}
    0 & 0 & 1 \\
    0 & 0 & 0 \\
    0 & 0 & 0
  \end{pmatrix}\!,
\end{align}
where $[d_i,e_i]=a \label{eq:commutation rel h}$ for all $i=1,2,\dots,n$.
The Levi subalgebra $\mfrak{l}$ of $\mfrak{p}$ is the linear span of
\begin{align}
  h_1 =\begin{pmatrix}
    1 & 0 & 0 \\
    0 & 0 & 0 \\
    0 & 0 & -1
  \end{pmatrix}\!,\quad
  h_2 =\begin{pmatrix}
    1 & 0 & 0 \\
    0 & -{2\over n} I_n & 0 \\
    0 & 0 & 1
  \end{pmatrix}\!,\quad
  h_A =\begin{pmatrix}
    0 & 0 & 0 \\
    0 & A & 0 \\
    0 & 0 & 0
  \end{pmatrix}\!,
\end{align}
where $A \in M_{n\times n}(\C)$ satisfies $\tr A = 0$. Moreover, the elements $h_1$ and $h_2$ form a basis of the center $\mfrak{z}(\mfrak{l})$ of $\mfrak{l}$.

Finally, the parabolic subgroup $P$ of $G$ with the Lie algebra $\mfrak{p}$ is defined by
\begin{align}
  P=\left\{\begin{pmatrix}
    a & x^{\rm T} & c \\
    0 & A & y \\
    0 & 0 & b
  \end{pmatrix}\!;\,a,b,c\in \C,\, x,y \in \C^n,\, A \in M_{n\times n}(\C),\, ab\det(A)=1\right\}\!. \label{eq:parabolic subgroup}
\end{align}
Any character $\lambda \in \Hom_P(\mfrak{p},\C)$ is given by
\begin{align}
  \lambda=\lambda_1 \widetilde{\omega}_1 + \lambda_2 \widetilde{\omega}_{n+1} \label{eq:scalarchar}
\end{align}
for some $\lambda_1,\lambda_2 \in \C$, where $\widetilde{\omega}_1, \widetilde{\omega}_{n+1}\in \Hom_P(\mfrak{p},\C)$
are equal to $\omega_1, \omega_{n+1} \in \mfrak{h}^*$ and regarded as trivially extended to
$\mfrak{p} =\mfrak{h} \oplus (\bigoplus_{\alpha \in \Delta_\Sigma} \mfrak{g}_\alpha) \oplus \mfrak{u}$.
Then the vector $\rho \in \Hom_P(\mfrak{p},\C)$ defined by the formula \eqref{eq:rho vector} is given by
\begin{align}
  \rho={\textstyle {n+1 \over 2}}\,\widetilde{\omega}_1 +{\textstyle {n+1 \over 2}}\,\widetilde{\omega}_{n+1}.
\end{align}
For $r \in \{0,1,\dots,n-1\}$ we define an injective Lie algebra homomorphism
\begin{align}
i_r \colon \mfrak{sl}(n-r+2,\C) \rarr \mfrak{sl}(n+2,\C)
\end{align}
by
\begin{align}
  \begin{pmatrix}
    a & x^{\rm T} & b \\
    y & A & u \\
    c & v^{\rm T} & d
  \end{pmatrix} \mapsto
  \begin{pmatrix}
    a & x^{\rm T} & 0 & b \\
    y & A & 0 & u \\
    0 & 0 & 0 & 0 \\
    c & v^{\rm T} & 0 & d
  \end{pmatrix}\!,
\end{align}
where $a,b,c,d \in \C$, $x,y,u,v \in \C^{n-r}$, $A \in M_{n-r \times n-r}(\C)$ and $a+d+\tr A =0$.
Then we regard the semisimple Lie subalgebra $\mfrak{g}'_r$ of $\mfrak{g}$ to be the image
of the mapping $i_r$. We set $\mfrak{p}'_r = \mfrak{g}'_r \cap \mfrak{p}$, so that $\mfrak{p}'_r$
is a parabolic subalgebra of $\mfrak{g}'_r$ and
$\mfrak{p}'_r = \mfrak{l}'_r \oplus \mfrak{u}'_r$ with $\mfrak{l}'_r = \mfrak{g}'_r \cap \mfrak{l}$,
$\mfrak{u}'_r = \mfrak{g}'_r \cap \mfrak{u}$. Finally, we obtain
$\mfrak{l}'_r=\mfrak{z}(\mfrak{l}'_r) \oplus [\mfrak{l}'_r,\mfrak{l}'_r]$ for
$[\mfrak{l}'_r,\mfrak{l}'_r]= \mfrak{g}'_r \cap [\mfrak{l},\mfrak{l}]$, and $h_1'$, $h_2'$, where
\begin{align}
  h_1' = h_1, \qquad h_2'= h_2 + h_{A_r}
\end{align}
for $A_r$ given by
\begin{align}
  A_r = \begin{pmatrix}
    -{2r \over n(n-r)}I_{n-r} & 0 \\
    0 & {2\over n} I_r
  \end{pmatrix}\!,
\end{align}
form a basis of the center $\mfrak{z}(\mfrak{l}_r')$ of $\mfrak{l}_r'$.


\subsection{The branching problem for the pair $(\mfrak{g},\mfrak{g}'_r)$}
\label{sec:qbp}

This section is a short digression in which we explain
and determine, as a consequence of character formulas, qualitative
properties of the branching
problem for the pair $(\mfrak{g}, \mfrak{g}'_r)$ and the scalar generalized Verma module $M^\mfrak{g}_\mfrak{p}(\lambda)$
(see \eqref{eq:generalized Verma module}) induced from a character $\lambda \in \Hom_P(\mfrak{p},\C)$ (see \eqref{eq:scalarchar}).

Let us denote by $\C_\lambda$ the $1$-dimensional $\mfrak{l}$-module
given by a highest weight $\lambda=\lambda_1\omega_1+\lambda_2\omega_{n+1}$.
The induced $\mfrak{l}'_r$-module on the symmetric algebra $S(\widebar{\mfrak{u}}/\widebar{\mfrak{u}}'_r)$, which is a free commutative
$\C$-algebra generated by the $2r$-dimensional complex vector space
$\widebar{\mfrak{u}}/\widebar{\mfrak{u}}'_r$, is isomorphic to the direct sum of $1$-dimensional $\mfrak{l}_r'$-modules. Then it follows that in the branching problem for the pair $(\mfrak{g},\mfrak{g}'_r)$ and the
scalar generalized Verma $\mfrak{g}$-module $M^\mfrak{g}_\mfrak{p}(\lambda)$
appear only scalar generalized Verma $\mfrak{g}'_r$-modules
$\smash{M^{\mfrak{g}'_r}_{\mfrak{p}'_r}(\mu)}$, their multiplicities (cf.\ \eqref{eq:multiplicity formula}) are given by
\begin{align}
  m(\lambda, \mu)=\dim\Hom_{\mfrak{l}'_r}\!(\C_\mu,
  \C_\lambda|_{\mfrak{l}'_r}\otimes_\C S(\widebar{\mfrak{u}}/\widebar{\mfrak{u}}'_r)), \label{eq:multiplicity}
\end{align}
and are equal to $\binom{a+r-1}{a}\binom{b+r-1}{b}$ provided the highest weight $\mu$ of $\mfrak{l}'_r$ is equal to $(\lambda_1-a)\omega_1+(\lambda_2-b)\omega_{n-r+1}$ for $a,b \in \N_0$, and zero
otherwise.

In the Grothendieck group $K(\mcal{O}^\mfrak{p})$ of the Bernstein-Gelfand-Gelfand parabolic category
$\mcal{O}^\mfrak{p}$ (cf.\ \eqref{eq:reduction}) holds
\begin{align}
M^\mfrak{g}_\mfrak{p}(\lambda_1\omega_1+\lambda_2\omega_{n+1})|_{\mfrak{g}'_r}\simeq \bigoplus_{a,b \in \N_0}
{\textstyle \binom{a+r-1}{a}\binom{b+r-1}{b}} M^{\mfrak{g}'_r}_{\mfrak{p}'_r}((\lambda_1-a)\omega_1+(\lambda_2-b)\omega_{n-r+1}). \label{eq:diag-branch-tk}
\end{align}
In the rest of the article we construct, among others, the generators of
$\mfrak{g}'_r$-submodules on the right hand side of \eqref{eq:diag-branch-tk}.

\subsection{Embedding of $\mfrak{g}$ into the Weyl algebras
$\eus{A}^\mfrak{g}_{\bar{\mfrak{u}}}$ and $\eus{A}^\mfrak{g}_{\bar{\mfrak{u}}^*}$}

Let us denote by $(\hat{x}_1,\dots,\hat{x}_n,\hat{y}_1,\dots,\hat{y}_n,\hat{z})$ the linear coordinate functions on
$\widebar{\mfrak{u}}$ with respect to the basis $(f_1,\dots,f_n,g_1,\dots,g_n,c)$ of the opposite nilradical $\widebar{\mfrak{u}}$, and by $(x_1,\dots,x_n,y_1,\dots,y_n,z)$ the dual linear coordinate functions on $\widebar{\mfrak{u}}^*$. Then the Weyl algebra
$\eus{A}^\mfrak{g}_{\widebar{\mfrak{u}}}$ is generated by
\begin{align}
\{\hat{x}_1,\dots,\hat{x}_n,\hat{y}_1,\dots,\hat{y}_n,\hat{z},\partial_{\hat{x}_1},\dots,\partial_{\hat{x}_n}, \partial_{\hat{y}_1},\dots, \partial_{\hat{y}_n},\partial_{\hat{z}}\}
\label{eq:Weyl algebra generators}
\end{align}
and the Weyl algebra
$\eus{A}^\mfrak{g}_{\widebar{\mfrak{u}}^*}$ is generated by
\begin{align}
\{x_1,\dots,x_n,y_1,\dots,y_n,z,\partial_{x_1},\dots,\partial_{x_n}, \partial_{y_1},\dots, \partial_{y_n},\partial_{z}\}.
\label{eq:dual Weyl algebra generators}
\end{align}
The coordinate functions $(u_e^1,u_e^2,\dots,u_e^{2n+1})$ on $U_e$ are defined by
\begin{align}
  u_e(x)=\sum_{i=1}^n u^i_e(x)f_i + \sum_{i=1}^n u_e^{n+i}(x)g_i + u_e^{2n+1}(x)c \label{eq:coordinate functions}
\end{align}
for all $x \in U_e$.

Now we apply Theorem \ref{prop:operator realization} to
$\lambda \in \Hom_P(\mfrak{p},\C)$ and find an explicit realization of $\mfrak{g}$ as a Lie subalgebra of
$\eus{A}^\mfrak{g}_{\widebar{\mfrak{u}}}$. Thus, we get a homomorphism
\begin{align}
  \pi_\lambda=\Psi_{u_e} \circ \pi_e^{\lambda+\rho} \colon U(\mfrak{g}) \rarr \eus{A}^\mfrak{g}_{\widebar{\mfrak{u}}}
\end{align}
of associative $\C$-algebras. Further, we define the Fourier transform
\begin{align}
\mcal{F}\colon \eus{A}^\mfrak{g}_{\widebar{\mfrak{u}}} \rarr \eus{A}^\mfrak{g}_{\widebar{\mfrak{u}}^*} \label{eq:Fourier transform ex}
\end{align}
by \eqref{eq:Fourier transform} with respect to the generators \eqref{eq:Weyl algebra generators} and \eqref{eq:dual Weyl algebra generators}.
Since $\mcal{F}$ is an isomorphism of associative $\C$-algebras, the composition
\begin{align}
  \hat{\pi}_\lambda = \mcal{F} \circ \pi_\lambda,
\end{align}
gives the homomorphism
\begin{align}
  \hat{\pi}_\lambda \colon U(\mfrak{g}) \rarr \eus{A}^\mfrak{g}_{\widebar{\mfrak{u}}^*}.
\end{align}
of associative $\C$-algebras.
\smallskip

Let us introduce the notation
\begin{align}
  E_x= \sum_{j=1}^n x_j\partial_{x_j},\quad  E_z= z\partial_z,\quad  E_y= \sum_{j=1}^n y_j\partial_{y_j}
\end{align}
and
\begin{align}
  E_{\hat{x}}= \sum_{j=1}^n \hat{x}_j\partial_{\hat{x}_j},\quad  E_{\hat{z}}= \hat{z}\partial_{\hat{z}},\quad  E_{\hat{y}}= \sum_{j=1}^n \hat{y}_j\partial_{\hat{y}_j}
\end{align}
for the Euler homogeneity operators, and
\begin{align}
   \squares = \sum_{j=1}^n \partial_{x_j} \partial_{y_j},\quad q=\sum_{j=1}^n x_jy_j, \quad \hat{q}=\sum_{j=1}^n \hat{x}_j\hat{y}_j
\end{align}
for the Laplace operator and quadratic polynomials. In the following Theorem
we retain the notation of previous sections.
\medskip

\theorem{\label{reproper}Let $\lambda \in \Hom_P(\mfrak{p},\C)$. Then the embedding of
$\mfrak{g}$ into
$\eus{A}^\mfrak{g}_{\widebar{\mfrak{u}}}$ and $\eus{A}^\mfrak{g}_{\widebar{\mfrak{u}}^*}$
is given by
\begin{enumerate}
\item[1)]
\begin{align}
\begin{aligned}
  \pi_\lambda(f_i)&=-\partial_{\hat{x}_i}+{\textstyle {1\over 2}}\hat{y}_i\partial_{\hat{z}}, \\
  \pi_\lambda(g_i)&=-\partial_{\hat{y}_i}-{\textstyle {1\over 2}}\hat{x}_i\partial_{\hat{z}}, \\
  \pi_\lambda(c)&=-\partial_{\hat{z}}
\end{aligned}
\end{align}
for $i=1,2,\dots,n$;
\item[2)]
\begin{align}
\begin{aligned}
  \hat{\pi}_\lambda(f_i)&=-x_i-{\textstyle {1\over 2}}z\partial_{y_i}, \\
  \hat{\pi}_\lambda(g_i)&=-y_i+{\textstyle {1\over 2}}z\partial_{x_i}, \\
  \hat{\pi}_\lambda(c)&=-z
\end{aligned}
\end{align}
for $i=1,2,\dots,n$;
\item[3)]
\begin{align}
\begin{aligned}
  \pi_\lambda(h_1)&= E_{\hat{x}}  + E_{\hat{y}}+ 2E_{\hat{z}} + \lambda_1 + \lambda_2 + n+1, \\
  \pi_\lambda(h_2)&= \big(1+{\textstyle {2\over n}}\big)E_{\hat{x}}- \big(1+{\textstyle {2\over n}}\big)E_{\hat{y}} + \lambda_1 - \lambda_2, \\
  \pi_\lambda(h_A)&=- {\textstyle \sum_{i,j=1}^n} a_{ij}(\hat{x}_j\partial_{\hat{x}_i}-\hat{y}_i\partial_{\hat{y}_j})
\end{aligned}
\end{align}
for all $A \in M_{n\times n}(\C)$ satisfying $\tr A = 0$;
\item[4)]
\begin{align}
\begin{aligned}
  \hat{\pi}_\lambda(h_1)&= -E_x - E_y- 2E_z + \lambda_1 + \lambda_2 -(n+1), \\
  \hat{\pi}_\lambda(h_2)&= -\big(1+{\textstyle {2\over n}}\big)E_x+ \big(1+{\textstyle {2\over n}}\big)E_y + \lambda_1 - \lambda_2, \\
  \hat{\pi}_\lambda(h_A)&= {\textstyle \sum_{i,j=1}^n} a_{ij}(x_i\partial_{x_j}-y_j\partial_{y_i})
\end{aligned}
\end{align}
for all $A \in M_{n\times n}(\C)$ satisfying $\tr A = 0$;
\item[5)]
\begin{align}
\begin{aligned}
  \pi_\lambda(d_i)&=\hat{z}\partial_{\hat{y}_i}+ \hat{x}_i\big(E_{\hat{x}}+{\textstyle {1 \over 2}}E_{\hat{z}} + \lambda_1 + {\textstyle {1 \over 2}}(n+1)\!\big) - {\textstyle {1\over 2}} \hat{q}\big(\partial_{\hat{y}_i} - {\textstyle {1\over 2}}\hat{x}_i\partial_{\hat{z}}\big), \\
  \pi_\lambda(e_i)&=-\hat{z}\partial_{\hat{x}_i}+ \hat{y}_i\big(E_{\hat{y}}+{\textstyle {1 \over 2}}E_{\hat{z}} + \lambda_2 + {\textstyle {1 \over 2}}(n+1)\!\big) - {\textstyle {1\over 2}}\hat{q}\big(\partial_{\hat{x}_i} + {\textstyle {1\over 2}} \hat{y}_i\partial_{\hat{z}}\big),  \\
  \pi_\lambda(a)&=\hat{z}(E_{\hat{x}}+ E_{\hat{y}}+E_{\hat{z}}+\lambda_1 + \lambda_2 + n+1)+{\textstyle {1\over 2}} \hat{q}\big(E_{\hat{x}} - E_{\hat{y}}+ \lambda_1-\lambda_2 + {\textstyle {1 \over 2}}\hat{q}\partial_{\hat{z}} \big)
\end{aligned}
\end{align}
for $i=1,2,\dots,n$;
\item[6)]
\begin{align}
\begin{aligned}
  \hat{\pi}_\lambda(d_i)&=-y_i\partial_z+ \partial_{x_i}\!\big(E_x+{\textstyle {1\over 2}}E_z -\lambda_1 + {\textstyle {1\over 2}}(n-1)\!\big)- {\textstyle {1 \over 2}}\big(y_i+ {\textstyle {1 \over 2}}z\partial_{x_i}\big)\squares, \\
  \hat{\pi}_\lambda(e_i)&=x_i\partial_z+ \partial_{y_i}\!\big(E_y+{\textstyle{1\over 2}}E_z -\lambda_2 + {\textstyle {1\over 2}}(n-1)\!\big) - {\textstyle {1 \over 2}}\big(x_i - {\textstyle {1 \over 2}}z\partial_{y_i}\big)\squares, \\
  \hat{\pi}_\lambda(a)&=\partial_z(E_x+E_y+E_z - \lambda_1 - \lambda_2 + n)-{\textstyle {1\over 2}} \big(E_x - E_y -\lambda_1 + \lambda_ 2 - {\textstyle {1 \over 2}} z\squares\big)\squares
\end{aligned}
\end{align}
for $i=1,2,\dots,n$.
\end{enumerate}
}

\proof{We have from Theorem \ref{prop:operator realization}
\begin{align*}
  \pi^{\lambda+\rho}_e(Y) = -\sum_{i=1}^{2n+1}\bigg[{\ad_{u_e(x)} \over e^{\ad_{u_e(x)}}-{\rm id}_{\widebar{\mfrak{u}}}}\,Y\bigg]_i \partial_{u_e^i}
\end{align*}
for all $Y \in \widebar{\mfrak{u}}$. Since $\ad_{u_e(x)}^2\!Y=0$, the expansion of the exponential implies
\begin{align*}
  \pi^{\lambda+\rho}_e(Y) = -\sum_{i=1}^{2n+1}\Big[Y-{\textstyle {1\over 2}}\ad_{u_e(x)}\!Y\Big]_i \partial_{u_e^i}
\end{align*}
and the first item follows by \eqref{eq:coordinate functions}.
Similarly, we get from Theorem \ref{prop:operator realization}
\begin{align*}
  \pi^{\lambda+\rho}_e(Y) = \sum_{i=1}^{2n+1}\,[\ad_{u_e(x)}\!Y]_i\,\partial_{u_e^i} + (\lambda+\rho)(Y)
\end{align*}
for all $Y \in \mfrak{l}$. The statement of the third item then follows from \eqref{eq:coordinate functions}.
Finally, we have
\begin{align*}
  \pi^{\lambda+\rho}_e(Y)= -\sum_{i=1}^{2n+1}\bigg[{\ad_{u_e(x)}e^{\ad_{u_e(x)}} \over e^{\ad_{u_e(x)}}-{\rm id}_{\widebar{\mfrak{u}}}}\,(e^{-\ad_{u_e(x)}}Y)_{\widebar{\mfrak{u}}}\bigg]_i\,\partial_{u_e^i} + (\lambda+\rho)((e^{-\ad_{u_e(x)}}Y)_\mfrak{p})
\end{align*}
for all $Y \in \mfrak{g}$. Since $e^{-\ad_{u_e(x)}}Y=(e^{-\ad_{u_e(x)}}Y)_{\widebar{\mfrak{u}}} + (e^{-\ad_{u_e(x)}}Y)_\mfrak{p}$,
we can write
\begin{align*}
 (e^{-\ad_{u_e(x)}}Y)_{\widebar{\mfrak{u}}}= e^{-\ad_{u_e(x)}}Y- (e^{-\ad_{u_e(x)}}Y)_\mfrak{p}= e^{-\ad_{u_e(x)}}Y-Y- (e^{-\ad_{u_e(x)}}Y-Y)_\mfrak{p}
\end{align*}
for $Y \in \mfrak{p}$. Therefore, we have
\begin{align*}
  {\ad_{u_e(x)}e^{\ad_{u_e(x)}} \over e^{\ad_{u_e(x)}}-{\rm id}_{\widebar{\mfrak{u}}}}\,(e^{-\ad_{u_e(x)}}Y)_{\widebar{\mfrak{u}}}&= -\ad_{u_e(x)}\!Y-\ad_{u_e(x)}(e^{-\ad_{u_e(x)}}Y-Y)_\mfrak{p}  \\
  &\quad - {\ad_{u_e(x)} \over e^{\ad_{u_e(x)}}-{\rm id}_{\widebar{\mfrak{u}}}}\,(e^{-\ad_{u_e(x)}}Y-Y)_\mfrak{p}
\end{align*}
for $Y \in \mfrak{p}$. Furthermore, since $\mfrak{g}$ is $|2|$-graded, $\ad_{u_e(x)}^5\!Y=0$ for all $Y\in\mfrak{p}$, we obtain
\begin{align*}
   {\ad_{u_e(x)}e^{\ad_{u_e(x)}} \over e^{\ad_{u_e(x)}}-{\rm id}_{\widebar{\mfrak{u}}}}\,(e^{-\ad_{u_e(x)}}Y)_{\widebar{\mfrak{u}}}&= -\ad_{u_e(x)}\!Y-\ad_{u_e(x)}(e^{-\ad_{u_e(x)}}Y-Y)_\mfrak{p}- (e^{-\ad_{u_e(x)}}Y-Y)_\mfrak{p} \\
   &\quad + {1\over 2}\,\ad_{u_e(x)}(e^{-\ad_{u_e(x)}}Y-Y)_\mfrak{p}  -{1\over 12}\,\ad_{u_e(x)}^2(e^{-\ad_{u_e(x)}}Y-Y)_\mfrak{p}.
\end{align*}
Putting all ingredients together, we conclude
\begin{align*}
  {\ad_{u_e(x)}e^{\ad_{u_e(x)}} \over e^{\ad_{u_e(x)}}-{\rm id}_{\widebar{\mfrak{u}}}}\,(e^{-\ad_{u_e(x)}}Y)_{\widebar{\mfrak{u}}}&= -(\ad_{u_e(x)}\!Y)_{\widebar{\mfrak{u}}} + {1\over 2}\,(\ad_{u_e(x)}(\ad_{u_e(x)}\!Y)_\mfrak{p})_{\widebar{\mfrak{u}}} -{1\over 4}\,\ad_{u_e(x)}(\ad_{u_e(x)}^2\!Y)_\mfrak{p} \\
  &\quad +{1\over 12}\,\ad_{u_e(x)}^2(\ad_{u_e(x)}\!Y)_\mfrak{p} -{1\over 24}\,\ad_{u_e(x)}^2(\ad_{u_e(x)}^2\!Y)_\mfrak{p}
\end{align*}
for $Y \in \mfrak{p}$, and the fifth item follows by \eqref{eq:coordinate functions}.

The computation of the Fourier transform of all operators is straightforward.
}


\subsection{Algebraic analysis on generalized Verma modules and singular vectors}

In what follows the generators $x_1,\dots,x_n,y_1,\dots,y_n,z$ of the symmetric algebra
$S(\widebar{\mfrak{u}})\simeq \C[\widebar{\mfrak{u}}^*]$ have the grading
$\deg(x_i)=\deg(y_i)=1$ for $i=1,2,\dots ,n$
and $\deg(z)=2$. Thus we regard $\C[\widebar{\mfrak{u}}^*]$ as a graded commutative $\C$-algebra.

Let us recall the existence of a canonical isomorphism of left
$\eus{A}^\mfrak{g}_{\widebar{\mfrak{u}}^*}$-modules
\begin{align}
  \C[\widebar{\mfrak{u}}^*]  \riso \eus{A}^\mfrak{g}_{\widebar{\mfrak{u}}^*}\!/\mcal{F}(I_e).
\end{align}
Let $\mu \colon \mfrak{z}(\mfrak{l}'_r) \rarr \C$ be a character of $\mfrak{z}(\mfrak{l}'_r)$. Then we have
\begin{align}
  \hat{\pi}_\lambda(h_1')R=\mu(h_1')R, \qquad
  \hat{\pi}_\lambda(h_2')R=\mu(h_2')R
\end{align}
for $R \in \Sol(\mfrak{g},\mfrak{g}'_r,\mfrak{p};\C[\widebar{\mfrak{u}}^*])_\mu^\mcal{F}$,
which is by Theorem \ref{reproper} equivalent to
\begin{gather}
(E_x'+E_y'+E_x''+E_y''+2E_z)R = (\lambda_1+\lambda_2-(n+1)-\mu(h_1'))R, \\
((n-r+2)(E_x' - E_y')+(n-r)(E_x''-E_y''))R = (n-r)(\lambda_1-\lambda_2-\mu(h_2'))R.
\end{gather}
Here we introduced the Euler homogeneity operators
\begin{align}
E_x'=\sum_{j=1}^{n-r} x_j\partial_{x_j},\quad E_x''=\sum_{j=n-r+1}^{n} x_j\partial_{x_j}, \quad E_y'=\sum_{j=1}^{n-r} y_j\partial_{y_j},\quad E_y''=\sum_{j=n-r+1}^{n} y_j\partial_{y_j}.
\end{align}
This gives a restriction on $\mu$ for which
$\Sol(\mfrak{g},\mfrak{g}'_r,\mfrak{p}; \C[\widebar{\mfrak{u}}^*])_\mu^\mcal{F}$
is a non-zero vector space, because the Euler homogenity operators
$E_x'$, $E_y'$, $E_x''$, $E_y''$ and $E_z$ acting on $\C[\widebar{\mfrak{u}}^*]$
have eigenvalues in $\N_0$. Thus, we may assume
\begin{gather}
(E_x'+E_y'+E_x''+E_y''+2E_z)R=mR, \label{eq:homogenity m} \\
((n-r+2)(E_x' - E_y')+(n-r)(E_x''-E_y''))R=tR \label{eq:homogenity t}
\end{gather}
for some $m \in \N_0$ and $t \in \Z$. Consequently, we shall apply the solution operators
\begin{align}
  \hat{\pi}_\lambda(d_i)R=0\qquad \text{and} \qquad \hat{\pi}_\lambda(e_i)R=0
\end{align}
for $i=1,2,\dots,n-r$ to polynomials $R$ of the form
\begin{align}
  R=\sum_{k=0}^{\lfloor {m \over 2} \rfloor} R_{m-2k}z^k,
\end{align}
where $R_{m-2k} \in \C[(\mfrak{g}_{-1})^*]$ satisfy
\begin{gather}
(E_x' + E_y'+E_x''+E_y'')R_{m-2k}=(m-2k)R_{m-2k}, \label{eq:Rm homogenity I} \\
((n-r+2)(E_x' - E_y')+(n-r)(E_x''-E_y''))R_{m-2k}=tR_{m-2k}, \label{eq:Rm homogenity II}
\end{gather}
and obtain the recurrence relations
\begin{multline}
  -(k+1)y_i R_{m-2k-2}+\partial_{x_i}\!\big(E_x +{\textstyle {1 \over 2}}k - \lambda_1 +{\textstyle {1 \over 2}}(n-1)\!\big) R_{m-2k} \\- {\textstyle {1 \over 2}}y_i \squares R_{m-2k} - {\textstyle {1 \over 4}} \partial_{x_i} \squares R_{m-2k+2}=0 \label{eq:recurrence P rel}
\end{multline}
and
\begin{multline}
   \phantom{-k}(k+1)x_i R_{m-2k-2}+\partial_{y_i}\!\big(E_y+{\textstyle {1 \over 2}}k - \lambda_2 +{\textstyle {1 \over 2}}(n-1)\!\big)R_{m-2k} \\ - {\textstyle {1 \over 2}}x_i \squares R_{m-2k} + {\textstyle {1 \over 4}} \partial_{y_i} \squares R_{m-2k+2}=0 \label{eq:recurrence Q rel}
\end{multline}
for $k=0,1,\dots, \lfloor {m\over 2} \rfloor$, where $R_{m-2k}=0$ for $k < 0$ and for $k> \lfloor {m\over 2} \rfloor$.
In particular, for $k=0$ we get
\begin{align}
  -y_iR_{m-2}+\partial_{x_i}\!\big(E_x - \lambda_1 +{\textstyle {1 \over 2}}(n-1)\!\big)R_m - {\textstyle {1 \over 2}}y_i \squares R_m &=0, \label{eq:recurrence P0 rel} \\
   x_iR_{m-2}+\partial_{y_i}\!\big(E_y - \lambda_2 +{\textstyle {1 \over 2}}(n-1)\!\big) R_m - {\textstyle {1 \over 2}}x_i \squares R_m &=0, \label{eq:recurrence Q0 rel}
\end{align}
which implies
\begin{align}
  \partial_{x_i}\!\big(E_x - \lambda_1 +{\textstyle {1 \over 2}}(n-1)\!\big)R_m &\in (y_i), \label{eq:ideal y} \\
  \partial_{y_i}\!\big(E_y - \lambda_2 +{\textstyle {1 \over 2}}(n-1)\!\big)R_m &\in (x_i) \label{eq:ideal x}
\end{align}
for $i=1,2,\dots,n-r$, where $(x_i)$ and $(y_i)$ are the ideals of $\C[(\mfrak{g}_{-1})^*]$ generated by $x_i$ and $y_i$, respectively.

Beside the structure of recurrence relations \eqref{eq:recurrence P rel}
and \eqref{eq:recurrence Q rel} we see that $R$ is uniquely
determined by $R_m$. Therefore, we can define a linear mapping
\begin{align}
  \Sol(\mfrak{g},\mfrak{g}'_r,\mfrak{p};\C[\widebar{\mfrak{u}}^*])^\mcal{F}_\mu
& \rarr  \C[(\mfrak{g}_{-1})^*],
\nonumber \\
R & \mapsto R_m,
\end{align}
which is injective and $\mfrak{l}'_r$-equivariant.
Since $\smash{\Sol(\mfrak{g},\mfrak{g}'_r,\mfrak{p}; \C[\widebar{\mfrak{u}}^*])^\mcal{F}_\mu}$
is a completely reducible $\mfrak{l}'_r$-module, any irreducible $\mfrak{l}'_r$-submodule
is contained in an isotypical component. Therefore, we can restrict to the case when $R_m$ is in an
isotypical component of the representation $\mfrak{l}'_r$ on $\C[(\mfrak{g}_{-1})^*]$.

Furthermore, the Fischer decomposition (cf.\ Appendix \ref{app:Fischer decompostion}) implies the isomorphism of vector spaces
\begin{align}
  \varphi \colon \C[\widebar{\mfrak{u}}^*] \riso  \bigoplus_{(a,b,c,d) \in \N_0^4}\C[q',q'',z]\otimes_\C
  \mcal{H}'_{a,b} \otimes_\C \mcal{H}''_{c,d}, \label{eq:isomorphism}
\end{align}
where
\begin{align}
q'=\sum_{j=1}^{n-r} x_jy_j,\qquad q''=\sum_{j=n-r+1}^{n} x_jy_j
\end{align}
are quadratic polynomials, $\mcal{H}'_{a,b}$ and $\mcal{H}''_{c,d}$ are the spaces of $(a,b)$-homogeneous
polynomials in the variables $(x_1,\dots ,x_{n-r},{y}_1,\dots ,{y}_{n-r})$ and $(c,d)$-homogeneous
polynomials in the variables $(x_{n-r+1},\dots ,x_{n},{y}_{n-r+1},\dots ,{y}_{n})$ harmonic for
\begin{align}
\squares'=\sum_{j=1}^{n-r}\partial_{x_j}\partial_{y_j}\qquad\text{and}\qquad \squares''=\sum_{j=n-r+1}^{n}\partial_{x_j}\partial_{y_j},
\end{align}
respectively. Consequently, we have
\begin{align}
  \varphi \circ \squares \circ \varphi^{-1} |_{\C[q',q'',z] \otimes_\C \mcal{H}'_{a,b} \otimes_\C \mcal{H}''_{c,d}} = Q_{a+b,c+d} \otimes 1 \otimes 1, \label{eq:Laplace tensor}
\end{align}
where
\begin{align}
  Q_{\alpha,\beta}=q'\partial_{q'}^2 + (n-r + \alpha)\partial_{q'} + q''\partial_{q''}^2 + (r+\beta) \partial_{q''}
\end{align}
for $\alpha,\beta\in \N_0$, and
\begin{align}
  \partial_{x_i}|_{\varphi^{-1}(\C[q',q'',z] \otimes_\C \mcal{H}'_{0,b} \otimes_\C \mcal{H}''_{c,d})} &= y_i(\varphi^{-1} \circ (\partial_{q'} \otimes 1 \otimes 1) \circ \varphi), \label{eq:partial x}\\
  \partial_{y_i}|_{\varphi^{-1}(\C[q',q'',z] \otimes_\C \mcal{H}'_{a,0} \otimes_\C \mcal{H}''_{c,d})} &= x_i(\varphi^{-1} \circ (\partial_{q'} \otimes 1 \otimes 1) \circ \varphi) \label{eq:partial y}
\end{align}
for $i=1,2,\dots,n-r$.
\medskip

\lemma{\label{lem:isotypical decomposition} Let us assume $n-r>2$. Then the isotypical components of the $\mfrak{l}'_r$-module $\C[\widebar{\mfrak{u}}^*]$ are of the form
\begin{align}\label{isotcomp}
  \bigoplus_{\ell=0}^{\min\{c_0,d_0\}} \varphi^{-1}(\mathbb{C}[q',q'',z]_\ell \otimes_\C \mcal{H}'_{a_0,b_0} \otimes_\C \mcal{H}''_{c_0-\ell,d_0-\ell})
\end{align}
for $a_0,b_0,c_0,d_0 \in \N_0$, where $\C[q',q'',z]_\ell \subset \C[q',q'',z]$ is the subspace of polynomials of degree $\ell$. These $\mfrak{l}_r'$-isotypical components are of the highest weight
$\mu_1\omega_1+a_0\omega_2+b_0\omega_{n-r} + \mu_2 \omega_{n-r+1}$, where
\begin{align}
  \mu_1&= \lambda_1 - {\textstyle {1 \over 2}}(n+1) - 2a_0 - c_0, \\
  \mu_2&= \lambda_2 - {\textstyle {1 \over 2}}(n+1) - 2b_0 - d_0.
\end{align}}

\proof{Since the mapping
\begin{align*}
   \varphi \colon \C[\widebar{\mfrak{u}}^*] \riso {\textstyle \bigoplus_{(a,b,c,d) \in \N_0^4}}\,\C[q',q'',z]\otimes_\C \mcal{H}'_{a,b} \otimes_\C \mcal{H}''_{c,d}
\end{align*}
is an isomorphism of vector spaces, we can define the structure of an $\mfrak{l}_r'$-module on the vector space ${\textstyle \bigoplus_{(a,b,c,d) \in \N_0^4}}\,\C[q',q'',z]\otimes_\C \mcal{H}'_{a,b} \otimes_\C \mcal{H}''_{c,d}$, so that $\varphi$ is an $\mfrak{l}'_r$-equivariant mapping. For the action of $[\mfrak{l}_r',\mfrak{l}_r']$ we get
\begin{align*}
  \varphi \circ \hat{\pi}_\lambda(h_{A'}) \circ \varphi^{-1} = 1 \otimes \hat{\pi}_\lambda(h_{A'})\otimes 1,
\end{align*}
where $h_{A'} \in [\mfrak{l}_r',\mfrak{l}_r']$. Because $\mcal{H}_{a_0,b_0}'$ is an irreducible $[\mfrak{l}_r',\mfrak{l}_r']$-module with the highest weight $a_0\omega_2 + b_0 \omega_{n-r}$ and the highest weight vector \smash{$x_1^{a_0}y_{n-r}^{b_0}$} (cf.\ Appendix \ref{app:Fischer decompostion}), we obtain that the isotypical component of $[\mfrak{l}_r',\mfrak{l}_r']$ with highest weight $a_0\omega_2+b_0\omega_{n-r}$ is
\begin{align*}
{\textstyle \bigoplus}_{(c,d) \in \N_0^2}\, \C[q',q'',z] \otimes_\C \mcal{H}_{a_0,b_0}'\otimes_\C \mcal{H}_{c,d}''.
\end{align*}
The generators of the center $\mfrak{z}(\mfrak{l}_r')$ act by
\begin{align*}
  \varphi \circ \hat{\pi}_\lambda(h_1') \circ \varphi^{-1} &= - 1\otimes (E_x'+E_y') \otimes 1 - 1\otimes 1 \otimes (E_x''+E_y'') \\ & \hspace{4cm} - 2 (q'\partial_{q'}+q''\partial_{q''}+E_z) \otimes 1 \otimes 1 + \lambda_1 + \lambda_2 -(n+1)
\intertext{and}
  \varphi \circ \hat{\pi}_\lambda(h_2') \circ \varphi^{-1} &=-\big(1+{\textstyle {2\over n-r}}\big) (1 \otimes (E_x'-E_y') \otimes 1) - 1 \otimes 1 \otimes (E_x''-E_y'') + \lambda_1-\lambda_2.
\end{align*}
Therefore, we have
\begin{align*}
  (\varphi \circ \hat{\pi}_\lambda(h_1') \circ \varphi^{-1})v&= -(a_0+b_0+c+d+2\ell-\lambda_1-\lambda_2+n+1)v, \\
  (\varphi \circ \hat{\pi}_\lambda(h_2') \circ \varphi^{-1})v&=-\big(\!\big(1+{\textstyle {2\over n-r}}\big)(a_0-b_0)+c-d -\lambda_1+\lambda_2 \big)v
\end{align*}
for $v \in \C[q',q'',z]_\ell \otimes_\C \mcal{H}_{a_0,b_0}'\otimes_\C \mcal{H}_{c,d}''$,
which implies that the center $\mfrak{z}(\mfrak{l}_r')$ acts by a character on the subspace
\begin{align*}
  \bigoplus_{\ell=0}^{\min\{c_0,d_0\}} \C[q',q'',z]_\ell \otimes_\C \mcal{H}'_{a_0,b_0} \otimes_\C \mcal{H}''_{c_0-\ell,d_0-\ell}
\end{align*}
of the vector space ${\textstyle \bigoplus}_{(c,d) \in \N_0^2}\, \C[q',q'',z] \otimes_\C \mcal{H}_{a_0,b_0}'\otimes_\C \mcal{H}_{c,d}''$.}

\medskip

Lemma \ref{lem:isotypical decomposition} yields the restrictions on
the rank of the Lie subalgebra $\mfrak{g}_r'=\mfrak{sl}(n-r+2,\C)$, because
for $n-r=2$ or $n-r=1$ the decomposition of $\C[\widebar{\mfrak{u}}^*]$
on $\mfrak{l}_r'$-isotypical components differs from \eqref{isotcomp}.
So we shall focus on $n-r>2$ and omit these two special cases from our
discussion.
Futhermore, by Lemma \ref{lem:isotypical decomposition} we see that the
isotypical components of the $\mfrak{l}'_r$-module $\C[(\mfrak{g}_{-1})^*]$
are uniquely determined by $a_0,b_0,c_0,d_0 \in \N_0$ and the highest weight
of the corresponding isotypical component is
\begin{align}
  \big(\lambda_1 - {\textstyle {1 \over 2}}(n+1) - 2a_0 - c_0\big) \omega_1 + a_0 \omega_2 + b_0 \omega_{n-r} + \big(\lambda_2 - {\textstyle {1 \over 2}}(n+1) - 2b_0 - d_0\big) \omega_{n-r+1}.
\end{align}
Let us consider
\begin{align}
  \varphi(R_m) \in  \bigoplus_{\ell=0}^{\min\{c_0,d_0\}} \mcal{P}_\ell \otimes_\C \mcal{H}'_{a_0,b_0} \otimes_\C \mcal{H}''_{c_0-\ell,d_0-\ell},
\end{align}
where $\mcal{P}_\ell\subset \C[q',q'']$ is the subspace of polynomials of degree $\ell$ and
$a_0+b_0+c_0+d_0=m$, $(n-r+2)(a_0-b_0)+(n-r)(c_0-d_0)=t$ as a consequence of \eqref{eq:homogenity m}, \eqref{eq:homogenity t}.
The notation
\begin{align}
  a_1(a,b,c,d)&=a+c-\lambda_1 + {\textstyle {1 \over 2}}(n-1), \\
  a_2(a,b,c,d)&=b+d-\lambda_2 + {\textstyle {1 \over 2}}(n-1)
\end{align}
for $a,b,c,d \in \N_0$, turns the conditions \eqref{eq:ideal y} and \eqref{eq:ideal x} into
\begin{align}
   \partial_{x_i}\!\big(E_x - \lambda_1 +{\textstyle {1 \over 2}}(n-1)\!\big)R_m &= a_1(a_0,b_0,c_0,d_0)\partial_{x_i}\!R_m,   \\
  \partial_{y_i}\!\big(E_y - \lambda_2 +{\textstyle {1 \over 2}}(n-1)\!\big)R_m & = a_2(a_0,b_0,c_0,d_0)\partial_{y_i}\!R_m
\end{align}
for $i=1,2,\dots,n-r$. Therefore, we get four mutually exclusive cases:
\begin{enum}
  \item[1)] $a_0 \neq 0$,\quad $b_0 \neq 0$,\quad $a_1(a_0,b_0,c_0,d_0)=0$,\quad $a_2(a_0,b_0,c_0,d_0)=0$;
  \item[2)] $a_0 \neq 0$,\quad $b_0=0$,\quad $a_1(a_0,b_0,c_0,d_0)=0$;
  \item[3)] $a_0 = 0$,\quad $b_0 \neq 0$,\quad $a_2(a_0,b_0,c_0,d_0)=0$;
  \item[4)] $a_0 = 0$,\quad $b_0 = 0$.
\end{enum}
\medskip

In what follows, we present a particular discussion of the construction
of $R$ for all these possibilities. First of all, we shall start with
two preparatory technical Lemmas.
\medskip

\lemma{\label{lem:condition lambda2 zero rel}
Let $\partial_{y_i}\!\big(E_y- \lambda_2 +{1 \over 2}(n-1)\!\big)R_m =0$ for $i=1,2,\dots,n-r$, then
\begin{align}
  R=\sum_{k=0}^{\lfloor {m\over 2} \rfloor} {z^k \over 2^k k!}\, \squares^k\! R_m, \label{eq:R_m lambda2 rel}
\end{align}
where $R_m$ satisfies
\begin{align}
  y_i\squares R_m=\partial_{x_i}\!\big(E_x - \lambda_1 +{\textstyle {1 \over 2}}(n-1)\!\big)R_m \label{eq:Laplace con2 I rel}
\end{align}
for $i=1,2,\dots,n-r$.}

\proof{First of all, we shall prove by induction that
$R_{m-2k}= {1 \over 2^k k!} \squares^k\!R_m$ and
$\partial_{y_i}\!\big(E_y - \lambda_2 + {\textstyle {1 \over 2}}(n-1)\!\big)R_{m-2k}=-kR_{m-2k}$ for $i=1,2,\dots,n-r$ and all $k \in \N_0$. This statement holds for $k=0$.
Let us assume that $R_{m-2k}= {1 \over 2^k k!} \squares^k\!R_m$ and $\partial_{y_i}\!\big(E_y - \lambda_2 + {\textstyle {1 \over 2}}(n-1)\!\big)R_{m-2k}=-kR_{m-2k}$ for $k=0,1,\dots,k_0$.
Then \eqref{eq:recurrence Q rel} for $k=k_0$ gives
\begin{align*}
   (k_0+1)x_i R_{m-2k_0-2}-{\textstyle {1 \over 2}}k_0 \partial_{y_i}R_{m-2k_0} - {\textstyle {1 \over 2}}x_i \squares R_{m-2k_0} + {\textstyle {1 \over 4}} \partial_{y_i} \squares R_{m-2k_0+2}=0,
\end{align*}
which is equivalent to
\begin{align*}
  x_i\big(\!(k_0+1)R_{m-2k_0-2}-{\textstyle {1 \over 2}}\,\squares R_{m-2k_0}\big)-{\textstyle {1 \over 2}}\,\partial_{y_i}\!\big(k_0R_{m-2k_0} - {\textstyle {1 \over 2}}\,\squares R_{m-2k_0+2}\big)=0.
\end{align*}
Since $k_0R_{m-2k_0} - {\textstyle {1 \over 2}}\,\squares R_{m-2k_0+2}=0$, we obtain
\begin{align*}
  R_{m-2k_0-2} = {\textstyle {1 \over 2(k_0+1)}}\,\squares R_{m-2k_0} = {\textstyle {1 \over 2^{k_0+1}(k_0+1)!}}\,\squares^{k_0+1}\!R_m
\end{align*}
and
\begin{align*}
  \partial_{y_i}\!\big(E_y - \lambda_2 + {\textstyle {1 \over 2}}(n-1)\!\big)R_{m-2k_0-2} &=  {\textstyle {1 \over 2(k_0+1)}}\, \partial_{y_i}\!\big(E_y - \lambda_2 + {\textstyle {1 \over 2}}(n-1)\!\big) \squares R_{m-2k_0} \\
  &= {\textstyle {1 \over 2(k_0+1)}}\, \squares \partial_{y_i}\!\big(E_y -1 - \lambda_2 + {\textstyle {1 \over 2}}(n-1)\!\big)  R_{m-2k_0} \\
  &={\textstyle {1 \over 2(k_0+1)}}\,(-k_0-1) \squares  R_{m-2k_0} \\
  &=-(k_0+1)  R_{m-2k_0-2},
\end{align*}
where we used $[\squares,E_y]= \squares$ in $\eus{A}^\mfrak{g}_{\widebar{\mfrak{u}}^*}$.

Now the formulas $2k R_{m-2k}=\squares R_{m-2k+2}$ and $2(k+1)R_{m-2k-2}=\squares R_{m-2k}$ reduce
\eqref{eq:recurrence P rel} into
\begin{align}
  y_i\squares R_{m-2k}=\partial_{x_i}\!\big(E_x - \lambda_1 +{\textstyle {1 \over 2}}(n-1)\!\big)R_{m-2k} \label{eq:compatibility Rm I rel}
\end{align}
for $k=0,1,\dots,\lfloor {m \over 2} \rfloor$. Let us assume \eqref{eq:compatibility Rm I rel}
holds for some $k \in \{0,1,\dots,\lfloor {m\over 2} \rfloor\}$. Then the application of the
Laplace operator $\squares$ to \eqref{eq:compatibility Rm I rel}, together with the relations
$[\squares,\partial_{x_i}]=0$, $[\squares,y_i]=\partial_{x_i}$ and $[\squares,E_x]=\squares$
in $\eus{A}^\mfrak{g}_{\widebar{\mfrak{u}}^*}$, gives
\begin{align*}
  y_i\squares^2\!R_{m-2k}=\partial_{x_i}\!\big(E_x - \lambda_1 +{\textstyle {1 \over 2}}(n-1)\!\big)\squares R_{m-2k},
\end{align*}
and the equality $2(k+1)R_{m-2k-2} = \squares R_{m-2k}$ implies
\begin{align*}
 y_i\squares R_{m-2k-2}=\partial_{x_i}\!\big(E_x - \lambda_1 +{\textstyle {1 \over 2}}(n-1)\!\big)R_{m-2k-2}.
\end{align*}
Therefore, the system of equations \eqref{eq:compatibility Rm I rel} for $k=0,1,\dots, \lfloor {m \over 2} \rfloor$ reduces to the single equation
\begin{align*}
  y_i\squares R_m=\partial_{x_i}\!\big(E_x - \lambda_1 +{\textstyle {1 \over 2}}(n-1)\!\big)R_m
\end{align*}
for $R_m$ and we are done.}

\lemma{\label{lem:condition lambda1 zero rel}Let $\partial_{x_i}\!\big(E_x - \lambda_1 +{\textstyle {1 \over 2}}(n-1)\!\big)R_m =0$ for $i=1,2,\dots,n-r$. Then
\begin{align}
  R=\sum_{k=0}^{\lfloor {m\over 2} \rfloor} {(-1)^kz^k \over 2^k k!}\, \squares^k\!R_m, \label{eq:R_m lambda1}
\end{align}
where $R_m$ satisfies
\begin{align}
  x_i\squares R_m=\partial_{y_i}\!\big(E_y - \lambda_2 +{\textstyle {1 \over 2}}(n-1)\!\big)R_m, \label{eq:Laplace con1 I rel}
\end{align}
for $i=1,2,\dots,n-r$.}

\proof{The proof can be given either by similar argument as in the previous
Lemma \ref{lem:condition lambda2 zero rel} or we can use the fact that
the system of equations \eqref{eq:recurrence Q rel} can be obtained from the system of equations
\eqref{eq:recurrence P rel} by the substitution $x_i \rightarrow y_i$, $y_i \rightarrow x_i$,
$z \rightarrow -z$ and $\lambda_1 \rarr \lambda_2$.}
\bigskip

\noindent
{\bf Case 1)} Let us assume
\begin{align}
 \varphi(R_m) \in \bigoplus_{\ell=0}^{\min\{c_0,d_0\}} \mcal{P}_\ell \otimes_\C \mcal{H}'_{a_0,b_0} \otimes_\C \mcal{H}''_{c_0-\ell,d_0-\ell},
\end{align}
where $a_1(a_0,b_0,c_0,d_0)=0$, $a_2(a_0,b_0,c_0,d_0)=0$ and $a_0 \neq 0$, $b_0 \neq 0$.
Then we have $\partial_{x_i}\!\big(E_x - \lambda_1 +{\textstyle {1 \over 2}}(n-1)\!\big)R_m = 0$ and  $\partial_{y_i}\!\big(E_y - \lambda_2 +{\textstyle {1 \over 2}}(n-1)\!\big)R_m = 0$ for $i=1,2,\dots,n-r$. Hence, from Lemma \ref{lem:condition lambda1 zero rel} we get
\begin{align}
  R=R_m,
\end{align}
where $R_m$ satisfies
\begin{align}
  \squares R_m=0. \label{eq:Laplace zero}
\end{align}
By the isomorphism \eqref{eq:isomorphism} and the formula \eqref{eq:Laplace tensor}, we can rewrite the equation \eqref{eq:Laplace zero}, into the form
\begin{align}
  (Q_{a_0+b_0,c_0+d_0-2\ell} \otimes 1 \otimes 1) \varphi(R_m)_\ell = 0
\end{align}
for $\ell=0,1,\dots,\min\{c_0,d_0\}$, where $\varphi(R_m)_\ell \in \mcal{P}_\ell \otimes_\C \mcal{H}'_{a_0,b_0} \otimes_\C \mcal{H}''_{c_0-\ell,d_0-\ell}$ and
\begin{align}
\varphi(R_m)=\sum_{\ell=0}^{\min\{c_0,d_0\}} \varphi(R_m)_\ell.
\end{align}
As we have
\begin{align}
  Q_{a_0+b_0,c_0+d_0-2\ell} =Q_{-c_0-d_0+\lambda_1+\lambda_2 -n+1,c_0+d_0-2\ell},
\end{align}
where we used $a_1(a_0,b_0,c_0,d_0)=0$ and $a_2(a_0,b_0,c_0,d_0) =0$, we get
\begin{align}
  \varphi(R) \in \bigoplus_{\ell=0}^{\min\{c_0,d_0\}} \mcal{S}_{\ell,c_0+d_0}^{\lambda-\rho} \otimes_\C \mcal{H}'_{a_0,b_0} \otimes_\C, \mcal{H}''_{c_0-\ell,d_0-\ell}
\end{align}
where the subspace $\mcal{S}_{\ell,c_0+d_0}^{\lambda-\rho}$ of $\mcal{P}_\ell$ is defined by
\begin{align}
  \mcal{S}_{\ell,s}^{\lambda-\rho}=\{u \in \mcal{P}_\ell;\, Q_{-s+\lambda_1+\lambda_2-n+1,s-2\ell} u=0\} \label{eq:solution Q}
\end{align}
for $s,\ell \in \N_0$.
\medskip

\noindent
{\bf Case 2)} Let us now suppose that
 \begin{align}
    \varphi(R_m) \in \bigoplus_{\ell=0}^{\min\{c_0,d_0\}} \mcal{P}_\ell \otimes_\C \mcal{H}'_{a_0,0} \otimes_\C \mcal{H}''_{c_0-\ell,d_0-\ell},
  \end{align}
where $a_1(a_0,0,c_0,d_0)=0$ and $a_0 \neq 0$. Then we have $\partial_{x_i}\!\big(E_x - \lambda_1 +{\textstyle {1 \over 2}}(n-1)\!\big)R_m = 0$ for $i=1,2,\dots,n-r$, and so by Lemma \ref{lem:condition lambda1 zero rel}
\begin{align}
    R=\sum_{k=0}^{\lfloor {m\over 2} \rfloor} {(-1)^kz^k \over 2^k k!}\, \squares^k\! R_m,
\end{align}
where $R_m$ satisfies
\begin{align}
   x_i\squares R_m=a_2(a_0,0,c_0,d_0)\partial_{y_i}R_m \label{eq:solution x}
\end{align}
for $i=1,2,\dots,n-r$. Using \eqref{eq:partial y}, we can write \eqref{eq:solution x} as
\begin{align}
  \squares R_m = a_2(a_0,0,c_0,d_0)(\varphi^{-1}\circ (\partial_{q'}\otimes 1 \otimes 1) \circ \varphi)R_m. \label{eq:Laplace a2}
\end{align}
The combination of \eqref{eq:isomorphism} and \eqref{eq:Laplace tensor} rewrites the equation \eqref{eq:Laplace a2} into the form
\begin{align}
  (Q_{a_0,c_0+d_0-2\ell} \otimes 1 \otimes 1)\varphi(R_m)_\ell= a_2(a_0,0,c_0,d_0)(\partial_{q'}\otimes 1 \otimes 1)\varphi(R_m)_\ell
\end{align}
for $\ell=0,1,\dots,\min\{c_0,d_0\}$. As we have
\begin{align}
  Q_{a_0,c_0+d_0-2\ell}-a_2(a_0,0,c_0,d_0)\partial_{q'} &= Q_{a_0-a_2(a_0,0,c_0,d_0),c_0+d_0-2\ell} \nonumber \\ &=Q_{-c_0-d_0+\lambda_1+\lambda_2 -n+1,c_0+d_0-2\ell},
\end{align}
due to $a_1(a_0,0,c_0,d_0)=0$, we get
\begin{align}
  \varphi(R_m) \in \bigoplus_{\ell=0}^{\min\{c_0,d_0\}} \mcal{S}_{\ell,c_0+d_0}^{\lambda-\rho} \otimes_\C \mcal{H}'_{a_0,0} \otimes_\C \mcal{H}''_{c_0-\ell,d_0-\ell}
\end{align}
for the subspace $\mcal{S}_{\ell,c_0+d_0}^{\lambda-\rho}$ of $\mcal{P}_\ell$ defined by \eqref{eq:solution Q}.
We shall now prove by induction
\begin{align}
  \squares^k\!R_m = k!\binom{a_2(a_0,0,c_0,d_0)}{k} (\varphi^{-1} \circ  (\partial_{q'}^k \otimes 1 \otimes 1) \circ \varphi)R_m
\end{align}
for all $k\in \N_0$, which evidently holds for $k=0$ and $k=1$. We may write
\begin{align}
  [\varphi \circ \squares \circ \varphi^{-1},\partial_{q'}^k \otimes 1 \otimes 1] = - k (\partial_{q'}^{k+1} \otimes 1 \otimes 1) \label{eq:laplace commutator}
\end{align}
for $k \in \N_0$, where we used \eqref{eq:Laplace tensor}. The previous formula and the induction hypothesis
lead to
\begin{align*}
  \squares^{k+1}\!R_m&=k!\binom{a_2(a_0,0,c_0,d_0)}{k} \squares (\varphi^{-1} \circ  (\partial_{q'}^k \otimes 1 \otimes 1)\circ \varphi)R_m \\
  &=k!\binom{a_2(a_0,0,c_0,d_0)}{k} (- k (\varphi^{-1} \circ (\partial_{q'}^{k+1} \otimes 1 \otimes 1) \circ \varphi)+(\varphi^{-1} \circ (\partial_{q'}^k \otimes 1 \otimes 1)\circ \varphi)\squares)R_m \\
  &= k!\binom{a_2(a_0,0,c_0,d_0)}{k} (a_2(a_0,0,c_0,d_0)-k)(\varphi^{-1} \circ (\partial_{q'}^{k+1}\otimes 1 \otimes 1)\circ \varphi))R_m \\
  &=(k+1)!\binom{a_2(a_0,0,c_0,d_0)}{k+1}(\varphi^{-1} \circ (\partial_{q'}^{k+1}\otimes 1 \otimes 1) \circ \varphi)R_m,
\end{align*}
where we used \eqref{eq:Laplace a2}. Denoting
\begin{align}
  T^{\lambda-\rho}_{a_0+c_0,d_0}= \sum_{k=0}^{\lfloor {m \over 2} \rfloor} {(-1)^k \over 2^k} \binom{a_2(a_0,0,c_0,d_0)}{k} (z\partial_{q'})^k,
\end{align}
we get
\begin{align}
  \varphi(R) \in \bigoplus_{\ell=0}^{\min\{c_0,d_0\}} T^{\lambda-\rho}_{a_0+c_0,d_0}\mcal{S}_{\ell,c_0+d_0}^{\lambda-\rho} \otimes_\C \mcal{H}'_{a_0,0} \otimes_\C \mcal{H}''_{c_0-\ell,d_0-\ell}.
\end{align}
\medskip

\noindent
{\bf Case 3)} We now assume
 \begin{align}
    \varphi(R_m) \in \bigoplus_{\ell=0}^{\min\{c_0,d_0\}} \mcal{P}_\ell \otimes_\C \mcal{H}'_{0,b_0} \otimes_\C \mcal{H}''_{c_0-\ell,d_0-\ell},
  \end{align}
where $a_2(0,b_0,c_0,d_0)=0$ and $b_0 \neq 0$. Then we have $\partial_{y_i}\!\big(E_y - \lambda_2 +{\textstyle {1 \over 2}}(n-1)\!\big)R_m = 0$ for $i=1,2,\dots,n-r$, and so by Lemma \ref{lem:condition lambda2 zero rel} we get
\begin{align}
    R=\sum_{k=0}^{\lfloor {m\over 2} \rfloor} {z^k \over 2^k k!}\, \squares^k\! R_m,
\end{align}
where $R_m$ satisfies
\begin{align}
   y_i\squares R_m=a_1(0,b_0,c_0,d_0)\partial_{x_i}R_m \label{eq:solution y}
\end{align}
for $i=1,2,\dots,n-r$. Using \eqref{eq:partial x}, we can rewrite \eqref{eq:solution y} as \begin{align}
  \squares R_m = a_1(0,b_0,c_0,d_0)(\varphi^{-1}\circ (\partial_{q'} \otimes 1 \otimes 1)\circ \varphi)R_m. \label{eq:Laplace a1}
\end{align}
The isomorphism \eqref{eq:isomorphism} and the formula
\eqref{eq:Laplace tensor} reduce the equation \eqref{eq:Laplace a1} into
\begin{align}
  (Q_{b_0,c_0+d_0-2\ell} \otimes 1 \otimes 1)\varphi(R_m)_\ell = a_1(0,b_0,c_0,d_0)(\partial_{q'}\otimes 1 \otimes 1)\varphi(R_m)_\ell,
\end{align}
for $\ell=0,1,\dots,\min\{c_0,d_0\}$. As we have
\begin{align}
  Q_{b_0,c_0+d_0-2\ell}-a_1(0,b_0,c_0,d_0)\partial_{q'} &= Q_{b_0-a_1(0,b_0,c_0,d_0),c_0+d_0-2\ell} \nonumber \\
  & = Q_{-c_0-d_0+\lambda_1+\lambda_2 -n+1,c_0+d_0-2\ell},
\end{align}
where we used $a_2(0,b_0,c_0,d_0)=0$, we get
\begin{align}
  \varphi(R_m) \in \bigoplus_{\ell=0}^{\min\{c_0,d_0\}} \mcal{S}_{\ell,c_0+d_0}^{\lambda-\rho} \otimes_\C \mcal{H}'_{0,b_0} \otimes_\C \mcal{H}''_{c_0-\ell,d_0-\ell}
\end{align}
for the subspace $\mcal{S}_{\ell,c_0+d_0}^{\lambda-\rho}$ of $\mcal{P}_\ell$ defined by \eqref{eq:solution Q}.
We shall prove by induction that
\begin{align}
  \squares^k\!R_m = k!\binom{a_1(0,b_0,c_0,d_0)}{k} (\varphi^{-1} \circ (\partial_{q'}^k \otimes 1 \otimes 1) \circ \varphi)R_m
\end{align}
for all $k\in \N_0$, which evidently holds for $k=0$ and $k=1$. Using \eqref{eq:laplace commutator} and the induction hypothesis, we obtain
\begin{align*}
  \squares^{k+1}\!R_m&=k!\binom{a_1(0,b_0,c_0,d_0)}{k} \squares (\varphi^{-1} \circ  (\partial_{q'}^k\otimes 1 \otimes 1) \circ \varphi)R_m \\
  &= k!\binom{a_1(0,b_0,c_0,d_0)}{k} (- k (\varphi^{-1} \circ (\partial_{q'}^{k+1} \otimes 1 \otimes 1) \circ \varphi)+(\varphi^{-1} \circ (\partial_{q'}^k \otimes 1 \otimes 1)\circ \varphi)\squares)R_m  \\
  &=k!\binom{a_1(0,b_0,c_0,d_0)}{k} (a_1(0,b_0,c_0,d_0)-k)(\varphi^{-1} \circ (\partial_{q'}^{k+1}\otimes 1 \otimes 1)\circ \varphi))R_m \\
  &=(k+1)!\binom{a_1(0,b_0,c_0,d_0)}{k+1}(\varphi^{-1} \circ (\partial_{q'}^{k+1} \otimes 1 \otimes 1) \circ \varphi)R_m,
\end{align*}
where we used \eqref{eq:Laplace a2}. Denoting
\begin{align}
  T^{\lambda-\rho}_{c_0,b_0+d_0} = \sum_{k=0}^{\lfloor {m \over 2} \rfloor} {1\over 2^k} \binom{a_1(0,b_0,c_0,d_0)}{k} (z\partial_{q'})^k,
\end{align}
we have
\begin{align}
  \varphi(R) \in \bigoplus_{\ell=0}^{\min\{c_0,d_0\}} T^{\lambda-\rho}_{c_0,b_0+d_0}\mcal{S}_{\ell,c_0+d_0}^{\lambda-\rho} \otimes_\C \mcal{H}'_{0,b_0} \otimes_\C \mcal{H}''_{c_0-\ell,d_0-\ell}.
\end{align}
\medskip

\noindent
{\bf Case 4)} Finally, let us suppose that
 \begin{align}
    \varphi(R_m) \in \bigoplus_{\ell=0}^{\min\{c_0,d_0\}} \mcal{P}_\ell \otimes_\C \mcal{H}'_{0,0} \otimes_\C  \mcal{H}''_{c_0-\ell,d_0-\ell}.
 \end{align}
Using \eqref{eq:partial x} and \eqref{eq:partial y}, we can rewrite \eqref{eq:recurrence P0 rel}
and \eqref{eq:recurrence Q0 rel} into the form
\begin{align}
  -R_{m-2}+a_1(0,0,c_0,d_0)(\varphi^{-1} \circ (\partial_{q'}\otimes 1 \otimes 1) \circ \varphi)R_m - {\textstyle {1 \over 2}} \squares R_m &=0, \\
   R_{m-2}+a_2(0,0,c_0,d_0)(\varphi^{-1} \circ (\partial_{q'}\otimes 1 \otimes 1) \circ \varphi)R_m - {\textstyle {1 \over 2}} \squares R_m &=0,
\end{align}
 and their sum yields
\begin{align}
  \squares R_m = (a_1(0,0,c_0,d_0)+a_2(0,0,c_0,d_0))(\varphi^{-1}\circ (\partial_{q'} \otimes 1 \otimes 1) \circ \varphi)R_m. \label{eq:Laplace a3}
\end{align}
The isomorphism \eqref{eq:isomorphism} and the formula \eqref{eq:Laplace tensor} allow
to rewrite \eqref{eq:Laplace a3} into
\begin{align}
  (Q_{0,c_0+d_0-2\ell} \otimes 1 \otimes 1)\varphi(R_m)_\ell = (a_1(0,0,c_0,d_0)+a_2(0,0,c_0,d_0))(\partial_{q'}\otimes 1 \otimes 1)\varphi(R_m)_\ell
\end{align}
for $\ell=0,1,\dots,\min\{c_0,d_0\}$. As we have
\begin{align}
  Q_{0,c_0+d_0-2\ell}-(a_1(0,0,c_0,d_0)+a_2(0,0,c_0,d_0))\partial_{q'}  = Q_{-c_0-d_0+\lambda_1+\lambda_2 -n+1,c_0+d_0-2\ell},
\end{align}
we get
\begin{align}
  \varphi(R_m) \in \bigoplus_{\ell=0}^{\min\{c_0,d_0\}} \mcal{S}_{\ell,c_0+d_0}^{\lambda-\rho} \otimes_\C \mcal{H}'_{0,0} \otimes_\C \mcal{H}''_{c_0-\ell,d_0-\ell},
\end{align}
where the subspace $\mcal{S}_{\ell,c_0+d_0}^{\lambda-\rho}$ of $\mcal{P}_\ell$ is defined by \eqref{eq:solution Q}.
We shall prove by induction that
\begin{align}
   \squares R_{m-2k} = (a_1(0,0,c_0,d_0)+a_2(0,0,c_0,d_0)-k)(\varphi^{-1}\circ (\partial_{q'}\otimes 1 \otimes 1) \circ \varphi) R_{m-2k} \label{eq:Laplace a12}
\end{align}
and
\begin{align}
  \varphi(R_{m-2k}) \in \bigoplus_{\ell=k}^{\min\{c_0,d_0\}} \mcal{P}_{\ell-k} \otimes_\C \mcal{H}'_{0,0} \otimes_\C  \mcal{H}''_{c_0-\ell,d_0-\ell} \label{eq:structure}
\end{align}
for all $k \in \N_0$. The claim holds for $k = 0$.
The recurrence relation \eqref{eq:recurrence P rel} for $k=k_0$ and the induction hypothesis imply
\begin{multline}
  (k_0+1)R_{m-2k_0-2}={\textstyle {1 \over 2}}(a_1(0,0,c_0,d_0) -a_2(0,0,c_0,d_0))(\varphi^{-1} \circ (\partial_{q'}\otimes 1 \otimes 1) \circ \varphi)R_{m-2k_0} \\- {\textstyle {1 \over 4}}(a_1(0,0,c_0,d_0)+a_2(0,0,c_0,d_0)-k_0+1)(\varphi^{-1} \circ (\partial_{q'}^2 \otimes 1 \otimes 1) \circ \varphi) R_{m-2k_0+2},
\end{multline}
which means that
\begin{align}
  \varphi(R_{m-2k_0-2}) \in \bigoplus_{\ell=k_0+1}^{\min\{c_0,d_0\}} \mcal{P}_{\ell-k_0-1} \otimes_\C \mcal{H}'_{0,0} \otimes_\C  \mcal{H}''_{c_0-\ell,d_0-\ell}.
\end{align}
By \eqref{eq:recurrence P rel} and \eqref{eq:recurrence Q rel} for $k=k_0+1$, we have
\begin{multline}
  -(k_0+2)R_{m-2k_0-4}+(a_1(0,0,c_0,d_0)- {\textstyle {1\over 2}}(k_0+1))(\varphi^{-1} \circ (\partial_{q'}\otimes 1 \otimes 1) \circ \varphi)R_{m-2k_0-2} \\ - {\textstyle {1 \over 2}} \squares R_{m-2k_0-2} - {\textstyle {1 \over 4}} (a_1(0,0,c_0,d_0)+a_2(0,0,c_0,d_0)-k_0)(\varphi^{-1}\circ (\partial_{q'}^2\otimes 1 \otimes 1) \circ \varphi)R_{m-2k_0}=0
\end{multline}
and
\begin{multline}
   \phantom{-k}(k_0+2)R_{m-2k_0-4}+(a_2(0,0,c_0,d_0)- {\textstyle {1\over 2}}(k_0+1))(\varphi^{-1} \circ (\partial_{q'}\otimes 1 \otimes 1) \circ \varphi)R_{m-2k_0-2} \\ - {\textstyle {1 \over 2}} \squares R_{m-2k_0-2} + {\textstyle {1 \over 4}} (a_1(0,0,c_0,d_0)+a_2(0,0,c_0,d_0)-k_0)(\varphi^{-1}\circ (\partial_{q'}^2\otimes 1 \otimes 1) \circ \varphi)R_{m-2k_0}=0,
\end{multline}
which implies
\begin{align}
   \squares R_{m-2k_0-2} = (a_1(0,0,c_0,d_0)+a_2(0,0,c_0,d_0)-k_0-1)(\varphi^{-1}\circ (\partial_{q'} \otimes 1 \otimes 1) \circ \varphi) R_{m-2k_0-2}.
\end{align}
As it follows from \eqref{eq:Laplace a12} and \eqref{eq:structure}, the recurrence relations \eqref{eq:recurrence P rel} and \eqref{eq:recurrence Q rel} are equivalent to
\begin{multline}
  (k+1)R_{m-2k-2}={\textstyle {1 \over 2}}(a_1(0,0,c_0,d_0) -a_2(0,0,c_0,d_0))(\varphi^{-1} \circ (\partial_{q'}\otimes 1\otimes 1) \circ \varphi)R_{m-2k} \\- {\textstyle {1 \over 4}}(a_1(0,0,c_0,d_0)+a_2(0,0,c_0,d_0)-k+1)(\varphi^{-1} \circ (\partial_{q'}^2 \otimes 1 \otimes 1) \circ \varphi) R_{m-2k+2} \label{eq:recurrence Q 12}
\end{multline}
for $k=0,1,\dots, \lfloor {m\over 2} \rfloor$, so we get by induction
\begin{align}
  R_{m-2k}={\alpha_k \over 2^k} (\varphi^{-1} \circ (\partial_{q'}^k \otimes 1 \otimes 1) \circ \varphi) R_m,
\end{align}
where $\alpha_k \in \C$ satisfy the following recurrence relation
\begin{align}
  (k+2)\alpha_{k+2}=(c_0-d_0-\lambda_1+\lambda_2)\alpha_{k+1}- (c_0+d_0-\lambda_1-\lambda_2+n-1-k)\alpha_k
\end{align}
with $\alpha_0=1$ and $\alpha_{-1}=0$. The notation
\begin{align}
  T^{\lambda-\rho}_{c_0,d_0} = \sum_{k=0}^{\lfloor {m \over 2} \rfloor} {\alpha_k \over 2^k}\,(z\partial_{q'})^k
\end{align}
allows to write
\begin{align}
  \varphi(R) \in \bigoplus_{\ell=0}^{\min\{c_0,d_0\}} T^{\lambda-\rho}_{c_0,d_0}\mcal{S}_{\ell,c_0+d_0}^{\lambda-\rho} \otimes_\C \mcal{H}'_{0,0} \otimes_\C \mcal{H}''_{c_0-\ell,d_0-\ell}.
\end{align}

\subsection{General structure of singular vectors}

In this subsection we summarize and organize the results achieved in
preceding sections
into particular statements describing the structure of singular vectors.
Let us recall the isomorphism
\begin{align}
  \tau \circ \Phi_{\lambda+\rho} \colon M_\mfrak{p}^\mfrak{g}(\lambda)^{\mfrak{u}'_r} \riso \Sol(\mfrak{g},\mfrak{g}'_r,\mfrak{p};\eus{A}^\mfrak{g}_{\widebar{\mfrak{u}}}/I_e) \riso \Sol(\mfrak{g},\mfrak{g}'_r,\mfrak{p};\C[\widebar{\mfrak{u}}^*])^\mcal{F}
\end{align}
of $\mfrak{l}_r'$-modules. Since $M_\mfrak{p}^\mfrak{g}(\lambda)^{\mfrak{u}'_r}$ is a completely reducible $\mfrak{l}'_r$-module, we can find its decomposition into the isotypical components. As a consequence of the previous considerations we see that the isotypical components can be uniformly written as
\begin{align}
  \bigoplus_{\ell=0}^{\min\{c_0,d_0\}} \varphi^{-1}(T^\lambda_{a_0+c_0,b_0+d_0}\mcal{S}_{\ell,c_0+d_0}^\lambda \otimes_\C \mcal{H}'_{a_0,b_0} \otimes_\C \mcal{H}''_{c_0-\ell,d_0-\ell}) \subset \Sol(\mfrak{g},\mfrak{g}_r',\mfrak{p};\C[\widebar{\mfrak{u}}^*])^\mcal{F},
\end{align}
where $a_0,b_0,c_0,d_0 \in \N_0$ and the subspace $\mcal{S}_{\ell,s}^\lambda \subset \mcal{P}_\ell$ is defined by
\begin{align}
  \mcal{S}_{\ell,s}^\lambda = \{u \in \mcal{P}_\ell;\, Q_{-s+\lambda_1+\lambda_2+2,s-2\ell} u=0\}.
\end{align}
The differential operator $T_{r_1,r_2}^\lambda \colon \C[q',q'',z] \rarr \C[q',q'',z]$ has the form
\begin{align}
  T_{r_1,r_2}^\lambda= \sum_{k=0}^{\lfloor {r_1+r_2 \over 2} \rfloor} {\alpha_k \over 2^k}\, (z\partial_{q'})^k,
\end{align}
where $\alpha_k \in \C$ satisfy the following recurrence relation
\begin{align}
  (k+2)\alpha_{k+2}=(r_1-r_2-\lambda_1+ \lambda_2)\alpha_{k+1}- (r_1+r_2-\lambda_1-\lambda_2-2-k)\alpha_k \label{eq:recurrence T}
\end{align}
with $\alpha_0=1$ and $\alpha_{-1}=0$. We denote this isotypical component by
$\mcal{V}_{a_0,b_0,c_0,d_0}^\lambda$.

Let us introduce the generating function for the collection $\{\alpha_k\}_{k\in{\mathbb N}_0}$,
\begin{align}
  g(w)= \sum_{k=0}^\infty \alpha_k w^k.
\end{align}
Then the recurrence relation \eqref{eq:recurrence T} is equivalent to the first order linear differential equation
\begin{align}
  \big(\!(1-w^2)\partial_w+(r_1+r_2-\lambda_1-\lambda_2-2)w-(r_1-r_2-\lambda_1+\lambda_2)\!\big)g(w)=0
\end{align}
for $g(w)$ with $g(0)=1$, whose unique solution is given by
\begin{align}
  g(w)=(1+w)^{r_1-\lambda_1-1}(1-w)^{r_2-\lambda_2-1}
\end{align}
with $g(0)=1$.
\bigskip

For the reader's convenience, we split our discussion according to the integrality
of the inducing weight $\lambda \in \Hom_P(\mfrak{p},\C)$.
\medskip

\theorem{\label{thm:decomposition I}Let us suppose that $\lambda_1,\lambda_2 \in \N_0$, and let $n-r>2$.
\begin{enum}
\item[i)] If $r>0$, then we have
\begin{equation}
\begin{split}
  \tau \circ \Phi_{\lambda+\rho} \colon M^\mfrak{g}_\mfrak{p}(\lambda)^{\mfrak{u}_r'} \riso & \bigoplus_{c_0=0}^{\lambda_1} \bigoplus_{d_0=0}^{\lambda_2} \mcal{V}^\lambda_{\lambda_1+1-c_0,\lambda_2+1-d_0,c_0,d_0} \oplus \bigoplus_{c_0=0}^{\lambda_1} \bigoplus_{d_0=0}^\infty \mcal{V}^\lambda_{\lambda_1+1-c_0,0,c_0,d_0}  \\ &
  \oplus \bigoplus_{c_0=0}^\infty \bigoplus_{d_0=0}^{\lambda_2} \mcal{V}^\lambda_{0,\lambda_2+1-d_0,c_0,d_0} \oplus \bigoplus_{c_0=0}^\infty \bigoplus_{d_0=0}^\infty \mcal{V}^\lambda_{0,0,c_0,d_0}.
\end{split}
\end{equation}
\item[ii)] If $r=0$, then we have
\begin{align}
  \tau \circ \Phi_{\lambda+\rho} \colon M^\mfrak{g}_\mfrak{p}(\lambda)^{\mfrak{u}} \riso
    \mcal{H}_{0,0} \oplus \mcal{H}_{\lambda_1+1,0} \oplus \mcal{H}_{0,\lambda_2+1} \oplus \mcal{H}_{\lambda_1+1,\lambda_2+1} \oplus P_\lambda \mcal{H}_{0,0},
\end{align}
where $P_\lambda=T^\lambda_{\lambda_1+\lambda_2+n+1,\lambda_1+\lambda_2+n+1} (q^{\lambda_1+\lambda_2+n+1})$.
\end{enum}}

\proof{The isotypical component $\mcal{V}^\lambda_{a,b,c,d}$ appears in the decomposition if and only if $a,b,c,d \in \N_0$ satisfy the constraints in Case 1 up to Case 4, i.e.\ $a \neq 0$, $b \neq 0$, $a+c-\lambda_1-1=0$, $b+d-\lambda_2-1=0$ in Case 1, $a \neq 0$, $b = 0$, $a+c-\lambda_1-1=0$ in Case 2, $a = 0$, $b \neq 0$, $b+d-\lambda_2-1=0$ in Case 3, and $a =0$, $b = 0$ in Case 4.}
\medskip

\theorem{\label{thm:decomposition II}Let us suppose that $\lambda_1 \in \N_0$ and $\lambda_2 \notin \N_0$, and let $n-r>2$.
\begin{enum}
\item[i)] If $r>0$, then we have
\begin{align}
  \tau \circ \Phi_{\lambda+\rho} \colon M^\mfrak{g}_\mfrak{p}(\lambda)^{\mfrak{u}_r'} \riso \bigoplus_{c_0=0}^{\lambda_1} \bigoplus_{d_0=0}^\infty \mcal{V}^\lambda_{\lambda_1+1-c_0,0,c_0,d_0}  \oplus \bigoplus_{c_0=0}^\infty \bigoplus_{d_0=0}^\infty \mcal{V}^\lambda_{0,0,c_0,d_0}.
\end{align}
\item[ii)] If $r=0$, then we have
\begin{align}
  \tau \circ \Phi_{\lambda+\rho} \colon M^\mfrak{g}_\mfrak{p}(\lambda)^{\mfrak{u}} \riso \begin{cases}
  \mcal{H}_{0,0} \oplus \mcal{H}_{\lambda_1+1,0} & \text{for $\lambda_1+\lambda_2+n \notin \N_0$},\\
  \mcal{H}_{0,0} \oplus \mcal{H}_{\lambda_1+1,0} \oplus P_\lambda\mcal{H}_{0,0} &
      \!\!\!\begin{array}{l}
      \text{for } \lambda_1+\lambda_2+n \in \N_0  \\[-1mm]
      \text{and } {-\lambda_2-n} \notin \N,
     \end{array}\\
  \mcal{H}_{0,0} \oplus \mcal{H}_{\lambda_1+1,0} \oplus P_\lambda\mcal{H}_{0,0} \oplus P_\lambda^1\mcal{H}_{-\lambda_2-n,0} &
      \!\!\!\begin{array}{l}
      \text{for } \lambda_1+\lambda_2+n \in \N_0  \\[-1mm]
      \text{and } {-\lambda_2-n} \in \N,
     \end{array}
  \end{cases}
\end{align}
where $P_\lambda=T^\lambda_{\lambda_1+\lambda_2+n+1,\lambda_1+\lambda_2+n+1} (q^{\lambda_1+\lambda_2+n+1})$ and $P_\lambda^1=T^\lambda_{\lambda_1+1,\lambda_1+\lambda_2+n+1} (q^{\lambda_1+\lambda_2+n+1})$.
\end{enum}}

\proof{The structure of the proof is identical as in Theorem \ref{thm:decomposition I}.}
\medskip

\theorem{\label{thm:decomposition III}Let us suppose that $\lambda_1 \notin \N_0$ and $\lambda_2 \in \N_0$, and let $n-r>2$.
\begin{enum}
\item[i)] If $r>0$, then we have
\begin{align}
  \tau \circ \Phi_{\lambda+\rho} \colon M^\mfrak{g}_\mfrak{p}(\lambda)^{\mfrak{u}'_r} \riso \bigoplus_{c_0=0}^{\infty} \bigoplus_{d_0=0}^{\lambda_2} \mcal{V}^\lambda_{0,\lambda_2+1-d_0,c_0,d_0}  \oplus \bigoplus_{c_0=0}^\infty \bigoplus_{d_0=0}^\infty \mcal{V}^\lambda_{0,0,c_0,d_0}.
\end{align}
\item[ii)] If $r=0$, then we have
\begin{align}
  \tau \circ \Phi_{\lambda+\rho} \colon M^\mfrak{g}_\mfrak{p}(\lambda)^{\mfrak{u}} \riso \begin{cases}
  \mcal{H}_{0,0} \oplus \mcal{H}_{0,\lambda_2+1} & \text{for $\lambda_1+\lambda_2+n \notin \N_0$},\\
  \mcal{H}_{0,0} \oplus \mcal{H}_{0,\lambda_2+1} \oplus P_\lambda \mcal{H}_{0,0} &
      \!\!\!\begin{array}{l}
      \text{for } \lambda_1+\lambda_2+n \in \N_0  \\[-1mm]
      \text{and } {-\lambda_2-n} \notin \N,
     \end{array}\\
  \mcal{H}_{0,0} \oplus \mcal{H}_{0,\lambda_2+1} \oplus P_\lambda \mcal{H}_{0,0} \oplus P_\lambda^2 \mcal{H}_{0,-\lambda_1-n} &
      \!\!\!\begin{array}{l}
      \text{for } \lambda_1+\lambda_2+n \in \N_0  \\[-1mm]
      \text{and } {-\lambda_2-n} \in \N,
     \end{array}
  \end{cases}
\end{align}
where $P_\lambda=T^\lambda_{\lambda_1+\lambda_2+n+1,\lambda_1+\lambda_2+n+1} (q^{\lambda_1+\lambda_2+n+1})$ and $P_\lambda^2=T^\lambda_{\lambda_1+\lambda_2+n+1,\lambda_2+1} (q^{\lambda_1+\lambda_2+n+1})$.
\end{enum}}

\proof{The structure of the proof is identical as in Theorem \ref{thm:decomposition I}.}
\medskip

\theorem{\label{thm:decomposition IV}Let us suppose that $\lambda_1,\lambda_2 \notin \N_0$, and let $n-r>2$.
\begin{enum}
\item[i)] If $r>0$, then we have
\begin{align}
   \tau \circ \Phi_{\lambda+\rho} \colon M^\mfrak{g}_\mfrak{p}(\lambda)^{\mfrak{u}'_r} \riso \bigoplus_{c_0=0}^\infty \bigoplus_{d_0=0}^\infty \mcal{V}^\lambda_{0,0,c_0,d_0}.
\end{align}
\item[ii)] If $r=0$, then we have
\begin{align}
  \tau \circ \Phi_{\lambda+\rho} \colon M^\mfrak{g}_\mfrak{p}(\lambda)^{\mfrak{u}} \riso \begin{cases}
    \mcal{H}_{0,0} & \text{for $\lambda_1+\lambda_2+n \notin \N_0$}, \\
    \mcal{H}_{0,0} \oplus P_\lambda\mcal{H}_{0,0}  & \text{for $\lambda_1+\lambda_2+n \in \N_0$},
  \end{cases}
\end{align}
where $P_\lambda=T^\lambda_{\lambda_1+\lambda_2+n+1,\lambda_1+\lambda_2+n+1} (q^{\lambda_1+\lambda_2+n+1})$.
\end{enum}}

\proof{The structure of the proof is identical as in Theorem \ref{thm:decomposition I}.}
\medskip

For $r=0$, it is straightforward to identify in previous theorems
the singular vectors corresponding to standard and
non-standard homomorphisms of generalized Verma modules,
respectively, cf.\ \cite{Lepowsky1977}.

\subsection{Examples}

In this subsection we illustrate the general results given in Theorem \ref{thm:decomposition I}, Theorem \ref{thm:decomposition II},
Theorem \ref{thm:decomposition III} and Theorem \ref{thm:decomposition IV} and write down explicit formulas for homomorphisms
between generalized Verma modules in several examples for $r=0$ and $r=1$.
\medskip

We shall start with $r=0$. Let $v_\lambda$ and $v_\mu$
be the highest weight vectors of $M^\mfrak{g}_\mfrak{p}(\lambda)$ and
$M^\mfrak{g}_\mfrak{p}(\mu)$, respectively. A homomorphism
\begin{align}
  \varphi \colon M^\mfrak{g}_\mfrak{p}(\mu) \rarr M^\mfrak{g}_\mfrak{p}(\lambda)
\end{align}
of generalized Verma modules is uniquely determined by $\varphi(v_\mu) \in M^\mfrak{g}_\mfrak{p}(\lambda)$.
For $n>2$, we have in the Poincaré-Birkhoff-Witt basis of $U(\widebar{\mfrak{u}})$:
\begin{enumerate}
  \item[(1)] If $\lambda=\lambda_1 \omega_1 + \lambda_2 \omega_{n+1}$,
	$\mu=(\lambda_1-1)\omega_1+(\lambda_2-1)\omega_{n+1}$, $\lambda_1 + \lambda_2 + n=0$,
  then the singular vector
  \begin{align}
    \varphi(v_\mu)=\big({\textstyle \sum_{i=1}^n} f_ig_i + {\textstyle {1 \over 2}}(\lambda_1-\lambda_2+n) c\big)v_\lambda,
	\end{align}
	in $M^\mfrak{g}_\mfrak{p}(\lambda)$ of graded homogeneity two induces a homomorphism of scalar generalized Verma modules,
	cf.\ Section \ref{sub:rtc} for the notation of root spaces.
  \item[(2)] If $\lambda=\lambda_1 \omega_1 + \lambda_2 \omega_{n+1}$, $\mu=(\lambda_1-2)\omega_1+(\lambda_2-2)\omega_{n+1}$, $\lambda_1 + \lambda_2 + n=1$, then the singular vector
  \begin{align}
    \varphi(v_\mu)=\big({\textstyle \sum_{i=1}^n} f_ig_i +
		{\textstyle {1 \over 2}}(\lambda_1-\lambda_2+n-1) c\big)\big({\textstyle \sum_{i=1}^n} f_ig_i
		+ {\textstyle {1 \over 2}}(\lambda_1-\lambda_2+n+1) c\big)v_\lambda
  \end{align}
 in $M^\mfrak{g}_\mfrak{p}(\lambda)$ of graded homogeneity four induces a homomorphism of scalar generalized Verma modules.
 \item[(3)] If $\lambda=0$, $\mu=-2\omega_1 + \omega_2 + \omega_n - 2 \omega_{n+1}$, then Theorem \ref{thm:decomposition I} implies that
  \begin{align}
    \varphi(v_\mu)=(f_1g_n) v_\lambda
  \end{align}
	induces a homomorphism from the vector valued generalized Verma module
	$M^\mfrak{g}_\mfrak{p}(\mu)$ to $M^\mfrak{g}_\mfrak{p}(\lambda)$.
\end{enumerate}
The examples (1) and (2) suggest the following remarkable conjectural factorization property of singular vectors
$\varphi(v_\mu) \in M^\mfrak{g}_\mfrak{p}(\lambda)$ inducing homomorphisms of scalar generalized Verma modules, see Theorem \ref{thm:symmetrization}, regarded by \eqref{eq:inverse mapping} as elements of the (non-commutative)
universal enveloping algebra $U(\widebar{\mfrak{u}})$. For all
$\lambda=\lambda_1\omega_1+\lambda_2 \omega_{n+1}$, $\mu=(\lambda_1-a)\omega_1 + (\lambda_2-a)\omega_{n+1}$,
 $\lambda_1+\lambda_2+n=a-1$ and $a \in \N$, we conjecture to hold
\begin{align}
  \varphi(v_\mu)= \prod_{j=0}^{a-1}\big({\textstyle \sum_{i=1}^n} f_ig_i
	+ {\textstyle {1 \over 2}}(\lambda_1-\lambda_2 + n - a + 2j+1)c \big)v_\lambda.
\end{align}
We do not know a direct proof of this observation.
\medskip

Now we pass to $r=1$.
Denoting $v'_\mu$ the highest weight vector of $\smash{M^{\mfrak{g}_1'}_{\mfrak{p}_1'}(\mu)}$,
a $\mfrak{g}_1'$-homomorphism
\begin{align}
  \varphi \colon \smash{M^{\mfrak{g}_1'}_{\mfrak{p}_1'}(\mu)} \rarr M^\mfrak{g}_\mfrak{p}(\lambda)
\end{align}
of generalized Verma modules is again uniquely determined by
$\varphi(v_\mu') \in M^\mfrak{g}_\mfrak{p}(\lambda)$. For $n>3$, we
have in the Poincaré-Birkhoff-Witt basis of $U(\widebar{\mfrak{u}})$:
\begin{enumerate}
  \item[(1)] If $\lambda=\lambda_1 \omega_1 + \lambda_2 \omega_{n+1}$,
	$\mu=(\lambda_1-a)\omega_1+\lambda_2\omega_n$, $a \in \N_0$,
  then the singular vector of graded homogeneity $a$ in $M^\mfrak{g}_\mfrak{p}(\lambda)$,
  \begin{align}
    \varphi(v'_\mu)=f_n^a v_\lambda,
	\end{align}
	induces a $\mfrak{g}_1'$-homomorphism of scalar generalized Verma modules.
 \item[(2)] If $\lambda=\lambda_1 \omega_1 + \lambda_2 \omega_{n+1}$,
	$\mu=\lambda_1\omega_1+(\lambda_2-a)\omega_n$, $a \in \N_0$,
  then the singular vector of graded homogeneity $a$ in $M^\mfrak{g}_\mfrak{p}(\lambda)$,
  \begin{align}
    \varphi(v'_\mu)=g_n^a v_\lambda,
	\end{align}
	induces a $\mfrak{g}_1'$-homomorphism of scalar generalized Verma modules.
\item[(3)] If $\lambda=\lambda_1 \omega_1 + \lambda_2 \omega_{n+1}$,
	$\mu=(\lambda_1-a-1)\omega_1+(\lambda_2-1)\omega_n$, $a \in \N_0$,
  then the singular vector of graded homogeneity $a+2$ in
	$M^\mfrak{g}_\mfrak{p}(\lambda)$,
 \begin{align}
    \varphi(v'_\mu)=\big({\textstyle \sum_{i=1}^{n-1}} g_if_i + (\lambda_1-a)c
		- {\textstyle {\lambda_1 + \lambda_2 + n -1-a \over a+1}}\,g_nf_n \big)f_n^a v_\lambda,
	\end{align}
	induces a $\mfrak{g}_1'$-homomorphism of scalar generalized Verma modules.
\item[(4)] If $\lambda=\lambda_1 \omega_1 + \lambda_2 \omega_{n+1}$,
	$\mu=(\lambda_1-1)\omega_1+(\lambda_2-a-1)\omega_n$, $a \in \N_0$,
  then the singular vector of graded homogeneity $a+2$
	in $M^\mfrak{g}_\mfrak{p}(\lambda)$,
\begin{align}
    \varphi(v'_\mu)=\big({\textstyle \sum_{i=1}^{n-1}} f_ig_i -(\lambda_2-a)c
		- {\textstyle {\lambda_1 + \lambda_2 + n -1-a \over a+1}}\,f_ng_n \big)g_n^av_\lambda,
	\end{align}
	induces a $\mfrak{g}_1'$-homomorphism of scalar generalized Verma modules.
\end{enumerate}

We notice that for $r=0$, the analogues of scalar valued singular vectors
for the pair of orthogonal Lie algebras and their compatible conformal parabolic
subalgebras with commutative nilradicals were constructed in
\cite{Juhl-book, koss}.


\section{CR-equivariant differential operators}
\label{sec:CR operators}

A CR  (a shorthand notation for Cauchy-Riemann or Complex-Real) structure is given by a real
differentiable manifold $M$ together with a complex distribution $H$, i.e.\ a complex
vector subbundle of the complexified tangent bundle $TM_\C$ such that
\begin{enumerate}
\item[(1)] $[H,H]\subset H$, i.e.\ $H$ is integrable,
\item[(2)] $H\cap \widebar{H}=\{0\}$, i.e.\ $H$ is almost Lagrangian.
\end{enumerate}
A CR-density on $M$ is a section of the complex line bundle
$(\Lambda^{n,0}H^*)^\lambda \otimes (\Lambda^{0,n}H^*)^{\lambda'}$,
defined for all $\lambda,\lambda'\in \C$ such that
$\lambda-\lambda'\in \Z$.
Manifolds endowed with CR-structure emerge on the boundaries of
strictly pseudo-convex domains
in $\C^{n+1}$, and invariants of integrable strictly pseudo-convex
CR-structures are used to obtain an expression for the asymptotic
expansion of the Bergman kernel in $\C^{n+1}$.
For the introduction and various aspects of CR-structure, we refer to
\cite{Tanaka-book, Harvey-Lawson1975, Fefferman1976, Branson-Fontana-Morpurgo2013, Collingwood-Shelton1990, Cap-Slovak-book} and references therein.

Our results are directly related to the homogeneous model of CR-structure, which
is the real $(2n+1)$-dimensional sphere $S^{2n+1}\subset\C^{n+1}$ realized
as a generalized flag manifold $S^{2n+1}=G_\R/P_\R$ for the real form
$G_\R={\rm SU}(n+1,1)\subset \SL(n+2, {\mathbb C})$
and $P_\R=G_\R\cap P$ the real form of the complex parabolic subgroup $P$ in \eqref{eq:parabolic subgroup}.
The vector subbundle $H$ is given by $H= (TS^{2n+1})_\C \cap T^{1,0}\C^{n+1}$,
where $T^{1,0}\C^{n+1}$ is the bundle of holomorphic vector fields on $\C^{n+1}$. The real
form of $H$ is the bundle ${\rm Re}(H\oplus \widebar{H})$, whose fiber at a point
$p\in S^{2n+1}$ is in terms of the complex structure $J$ on $\mathbb{C}^{n+1}$ given by
${\rm Re}(H\oplus \widebar{H})_p = \{ X\in T_pS^{2n+1};\, JX \in T_pS^{2n+1}\subset T_p\C^{n+1}\}$
and the almost complex structure on ${\rm Re}(H\oplus \widebar{H})$ is just the restriction of $J$.
The flat model of CR-manifold is the Heisenberg algebra.
Moreover, the notion of a CR submanifold in CR manifold $M$ is given by the real differentiable submanifold
$N\subset M$ and a complex vector subbundle $H'\subset H|_N$ on $N$ fulfilling analogous conditions
as $H$ does. In the case of homogeneous model of CR manifold $M=S^{2n+1}$, the real
codimension $2r$ homogeneous CR submanifold is the odd dimensional sphere $S^{2(n-r)+1}$ defined by
$x_{n}=x_{n-1}=\dots =x_{n-r+1}=0,\, y_{n}=y_{n-1}=\dots =y_{n-r+1}=0$.

The folklore results in \cite{Collingwood-Shelton1990} for $G'=G$, and in
\cite{koss, Kobayashi-Pevzner2013}
for the pair $G' \subset G$
and the pair of compatible parabolic subgroups $P' \subset P$, relate the existence of
homomorphisms between scalar generalized Verma modules to $G'$-covariant
differential operators acting on principal series representations
for $G$ and $G'$ supported on $G/P$ and $G'\!/P'$ ($G'\!/P'\subset G/P$), respectively.
There is a contravariant bijection
\begin{align}
\Hom_{(\mfrak{g}',N')}(M^{\mfrak{g}'}_{\mfrak{p}'}\!(W), M^\mfrak{g}_\mfrak{p}(V))
\simeq \Hom_{\Diff}(\Ind^G_P(V^*), \Ind^{G'}_{P'}(W^*)), \label{gvmdo}
\end{align}
where the subscript $\mathrm{Diff}$ denotes the homomorphisms given
by differential operators.
The generalized flag manifolds $G/P$ and $G'\!/P'$
are complexifications of the real generalized flag manifolds $G_\R/P_\R\simeq S^{2n+1}$ and
$G'_\R/P'_\R\simeq S^{2(n-r)+1}$
($G'_\R/P'_\R\subset G_\R/P_\R$), respectively. Thus
the $G'$-covariant differential operators in \eqref{gvmdo} induce, by restriction to the real
locus of the complex generalized flag manifolds, $G'_\R$-covariant differential
operators acting between sections of homogeneous vector bundles on $S^{2n+1}$ and $S^{2(n-r)+1}$,
respectively.

It is rather straightforward to transfer the singular vectors in Theorem \ref{thm:decomposition I},
Theorem \ref{thm:decomposition II}, Theorem \ref{thm:decomposition III} and Theorem \ref{thm:decomposition IV}
to covariant differential operators acting on the non-compact model of induced representations, cf.\
\cite{koss} for similar questions related to parabolic subalgebras of Hermitean symmetric type.


\begin{appendices}

\section{The Fischer decomposition for $\mfrak{sl}(n,\C)$}
\label{app:Fischer decompostion}

In this section we recall the Fischer decomposition proved by E.\,Fischer (\cite{Fischer1918}) in 1917 and apply it to one specific example.

Let $\eus{A}_V$ be the Weyl algebra of a finite-dimensional complex vector space $V$ (see Section \ref{sec:Fourier transform} for the definition). Then $\C[V]$ is a left $\eus{A}_V$-module and
\begin{align}
  \eus{A}_V \simeq \C[V] \otimes_\C \C[V^*],
\end{align}
where we regard $\C[V^*]$ as the $\C$-algebra of constant coefficient differential operators. Let us consider a polynomial $q \in \C[V]$ and a constant coefficient differential operator $P \in \C[V^*]$. We say that the pair $(q,P)$ forms a Fischer pair on $V$, if we have the decomposition
\begin{align}
  \C[V] \simeq \C[q] \otimes_\C\! \mcal{H},
\end{align}
where $\mcal{H} = \ker P$ and $\C[q]$ is the $\C$-subalgebra of $\C[V]$ generated by $q$.

To construct some Fischer pair on $V$, let us consider linear coordinate functions $(x_1,x_2,\dots,x_n)$ on $V$, which give us a canonical isomorphism $\C[V] \simeq \C[x]$, $x=(x_1,x_2,\dots,x_n)$. We introduce the scalar product $\langle \cdot\,,\cdot \rangle$ on $\C[x]$ by
\begin{align}
  \langle p(x),q(x) \rangle = p^*\!(\partial_x)q(x)|_{x=0}, \label{eq:scalar product}
\end{align}
where $\partial_x=(\partial_{x_1},\partial_{x_2},\dots,\partial_{x_n})$ and $p^*\!(x)=\smash{\overline{p(\widebar{x})}}$. Now, let $q \in \C[x]$ be a homogenous polynomial of degree $n \in \N$, and let us consider the mapping $T_q \colon \C[x] \rarr \C[x]$ given by $T_q(p)=qp$. Then we have
\begin{align}
  \langle T_q(p_1),p_2 \rangle = \langle p_1,q^*\!(\partial_x)p_2\rangle
\end{align}
for all $p_1,p_2 \in \C[x]$, hence $q^*\!(\partial_x) \colon \C[x] \rarr \C[x]$ is the adjoint mapping to $T_q$. Therefore, $\C[x]$ decomposes into two orthogonal subspaces $\im T_q$ and $\ker q^*\!(\partial_x)$, and hence $(q,q^*\!(\partial_x))$ is a Fischer pair.

We use this decomposition to the following specific example. Let $\C^n$ and $(\C^n)^*$ be the standard representation and the standard dual representation of the Lie algebra $\mfrak{g}=\mfrak{sl}(n,\C)$, $n\geq 2$, respectively. Then we have the induced representation of $\mfrak{g}$ on $\C[(\C^n)^* \oplus \C^n]$. Using the canonical linear coordinate functions $y=(y_1,y_2,\dots,y_n)$  on $\C^n$ and the dual linear coordinate functions $x=(x_1,x_2,\dots,x_n)$ on $(\C^n)^*$, we get the isomorphism $\C[(\C^n)^*\oplus \C^n] \simeq \C[x,y]$. The induced representation of $\mfrak{g}$ on $\C[x,y]$ is given by
\begin{align}
  \pi(A)=\sum_{i,j=1}^n a_{ij}(x_i\partial_{x_j}-y_j\partial_{y_i})
\end{align}
for all $A \in \mfrak{g}$. Moreover, if we introduce on $\C[x,y]$ the scalar product by \eqref{eq:scalar product}, then
\begin{align}
  q={\textstyle \sum_{i=1}^n} x_iy_i\quad \text{and} \quad \squares = {\textstyle \sum_{i=1}^n} \partial_{x_i}\partial_{y_i}
\end{align}
form a Fischer pair, hence we have the decomposition
\begin{align}
  \C[x,y]=\C[q] \otimes_\C\! \mcal{H},
\end{align}
where $\mcal{H} = \ker \squares$ and $\C[q]$ is the $\C$-subalgebra of $\C[x,y]$ generated by $q$. Since $q$ is a $\mfrak{g}$-invariant polynomial and $\squares$ is a $\mfrak{g}$-invariant differential operator, we obtain that $\C[q]$ and $\mcal{H}$ are representations of $\mfrak{g}$. Because the Euler homogeneity operators
\begin{align}
  E_x = {\textstyle \sum_{i=1}^n} x_i\partial_{x_i}\quad \text{and} \quad E_y = {\textstyle \sum_{i=1}^n} y_i\partial_{y_i}
\end{align}
satisfy $[\squares,E_x]=\squares$, $[\squares,E_y]=\squares$ and $[E_x,E_y]=0$, we get the direct sum decomposition
\begin{align}
  \mcal{H} = {\textstyle \bigoplus_{(a,b) \in \N_0^2}} \mcal{H}_{a,b}
\end{align}
into the common eigenspaces $\mcal{H}_{a,b}$ of $E_x$ and $E_y$ with the eigenvalues $a$ and $b$, respectively. As $E_x$ and $E_y$ are also $\mfrak{g}$-invariant differential operators, the space $\mcal{H}_{a,b}$ is a representation of $\mfrak{g}$. In fact, it is an irreducible representation of $\mfrak{g}$, since $\mcal{H}_{a,b}$ contains an irreducible representation of $\mfrak{g}$ with the highest weight $a\omega_1+b\omega_{n-1}$ ($x_1^ay_n^b$ is the highest weight vector with respect to the Borel subalgebra of upper triangular matrices) and its dimension is equal to the dimension of $\mcal{H}_{a,b}$.

\end{appendices}


\section*{Acknowledgments}

L.\,Křižka is supported by PRVOUK p47,
P.\,Somberg acknowledges the financial support from the grant GA\,P201/12/G028.
The second author is grateful to T. Kobayashi for sharing his ideas and experience.


\providecommand{\bysame}{\leavevmode\hbox to3em{\hrulefill}\thinspace}
\providecommand{\MR}{\relax\ifhmode\unskip\space\fi MR }
\providecommand{\MRhref}[2]{%
  \href{http://www.ams.org/mathscinet-getitem?mr=#1}{#2}
}
\providecommand{\href}[2]{#2}

\end{document}